\documentclass[12pt]{article}

\usepackage{amsmath,amssymb,amsthm,amscd,a4wide, amscd}
\usepackage{graphicx,color}
\begin{document}

\newcommand{\End}{{\rm{End}\ts}}
\newcommand{\Hom}{{\rm{Hom}}}
\newcommand{\Mat}{{\rm{Mat}}}
\newcommand{\ad}{{\rm{ad}\ts}}
\newcommand{\ch}{{\rm{ch}\ts}}
\newcommand{\chara}{{\rm{char}\ts}} 
\newcommand{\diag}{ {\rm diag}}
\newcommand{\pr}{^{\tss\prime}}
\newcommand{\non}{\nonumber}
\newcommand{\wt}{\widetilde}
\newcommand{\wh}{\widehat}
\newcommand{\ot}{\otimes}
\newcommand{\la}{\lambda}
\newcommand{\ls}{\ts\lambda\ts}
\newcommand{\La}{\Lambda}
\newcommand{\De}{\Delta}
\newcommand{\al}{\alpha}
\newcommand{\be}{\beta}
\newcommand{\ga}{\gamma}
\newcommand{\Ga}{\Gamma}
\newcommand{\ep}{\epsilon}
\newcommand{\ka}{\kappa}
\newcommand{\vk}{\varkappa}
\newcommand{\vt}{\vartheta}
\newcommand{\si}{\sigma}
\newcommand{\vp}{\varphi}
\newcommand{\de}{\delta}
\newcommand{\ze}{\zeta}
\newcommand{\om}{\omega}
\newcommand{\ee}{\epsilon^{}}
\newcommand{\su}{s^{}}
\newcommand{\hra}{\hookrightarrow}
\newcommand{\ve}{\varepsilon}
\newcommand{\ts}{\,}
\newcommand{\vac}{\mathbf{1}}
\newcommand{\di}{\partial}
\newcommand{\qin}{q^{-1}}
\newcommand{\tss}{\hspace{1pt}}
\newcommand{\Sr}{ {\rm S}}
\newcommand{\U}{ {\rm U}}
\newcommand{\BL}{ {\overline L}}
\newcommand{\BE}{ {\overline E}}
\newcommand{\BP}{ {\overline P}}
\newcommand{\AAb}{\mathbb{A}\tss}
\newcommand{\CC}{\mathbb{C}\tss}
\newcommand{\KK}{\mathbb{K}\tss}
\newcommand{\QQ}{\mathbb{Q}\tss}
\newcommand{\SSb}{\mathbb{S}\tss}
\newcommand{\ZZ}{\mathbb{Z}\tss}
\newcommand{\X}{ {\rm X}}
\newcommand{\Y}{ {\rm Y}}
\newcommand{\Z}{{\rm Z}}
\newcommand{\Ac}{\mathcal{A}}
\newcommand{\Lc}{\mathcal{L}}
\newcommand{\Mc}{\mathcal{M}}
\newcommand{\Pc}{\mathcal{P}}
\newcommand{\Qc}{\mathcal{Q}}
\newcommand{\Tc}{\mathcal{T}}
\newcommand{\Sc}{\mathcal{S}}
\newcommand{\Bc}{\mathcal{B}}
\newcommand{\Ec}{\mathcal{E}}
\newcommand{\Fc}{\mathcal{F}}
\newcommand{\Hc}{\mathcal{H}}
\newcommand{\Uc}{\mathcal{U}}
\newcommand{\Vc}{\mathcal{V}}
\newcommand{\Wc}{\mathcal{W}}
\newcommand{\Yc}{\mathcal{Y}}
\newcommand{\Ar}{{\rm A}}
\newcommand{\Br}{{\rm B}}
\newcommand{\Ir}{{\rm I}}
\newcommand{\Fr}{{\rm F}}
\newcommand{\Jr}{{\rm J}}
\newcommand{\Mr}{{\rm M}}
\newcommand{\Or}{{\rm O}}
\newcommand{\GL}{{\rm GL}}
\newcommand{\SL}{{\rm SL}}
\newcommand{\Spr}{{\rm Sp}}
\newcommand{\Rr}{{\rm R}}
\newcommand{\Zr}{{\rm Z}}
\newcommand{\gl}{\mathfrak{gl}}
\newcommand{\middd}{{\rm mid}}
\newcommand{\ev}{{\rm ev}}
\newcommand{\Pf}{{\rm Pf}}
\newcommand{\Norm}{{\rm Norm\tss}}
\newcommand{\oa}{\mathfrak{o}}
\newcommand{\spa}{\mathfrak{sp}}
\newcommand{\osp}{\mathfrak{osp}}
\newcommand{\bgot}{\mathfrak{b}}
\newcommand{\kgot}{\mathfrak{k}}
\newcommand{\g}{\mathfrak{g}}
\newcommand{\h}{\mathfrak h}
\newcommand{\n}{\mathfrak n}
\newcommand{\z}{\mathfrak{z}}
\newcommand{\Zgot}{\mathfrak{Z}}
\newcommand{\p}{\mathfrak{p}}
\newcommand{\sll}{\mathfrak{sl}}
\newcommand{\psl}{\mathfrak{psl}}
\newcommand{\agot}{\mathfrak{a}}
\newcommand{\qdet}{ {\rm qdet}\ts}
\newcommand{\Ber}{ {\rm Ber}\ts}
\newcommand{\HC}{ {\mathcal HC}}
\newcommand{\cdet}{ {\rm cdet}}
\newcommand{\tr}{ {\rm tr}}
\newcommand{\gr}{ {\rm gr}}
\newcommand{\str}{ {\rm str}}
\newcommand{\loc}{{\rm loc}}
\newcommand{\Gr}{{\rm G}}
\newcommand{\sgn}{ {\rm sgn}\ts}
\newcommand{\ba}{\bar{a}}
\newcommand{\bb}{\bar{b}}
\newcommand{\bi}{\bar{\imath}}
\newcommand{\bj}{\bar{\jmath}}
\newcommand{\bk}{\bar{k}}
\newcommand{\bl}{\bar{l}}
\newcommand{\hb}{\mathbf{h}}
\newcommand{\Sym}{\mathfrak S}
\newcommand{\fand}{\quad\text{and}\quad}
\newcommand{\Fand}{\qquad\text{and}\qquad}
\newcommand{\For}{\qquad\text{or}\qquad}
\newcommand{\OR}{\qquad\text{or}\qquad}

\renewcommand{\theequation}{\arabic{section}.\arabic{equation}}

\newtheorem{thm}{Theorem}[section]
\newtheorem{lem}[thm]{Lemma}
\newtheorem{prop}[thm]{Proposition}
\newtheorem{cor}[thm]{Corollary}
\newtheorem{conj}[thm]{Conjecture}
\newtheorem*{mthm}{Main Theorem}
\newtheorem*{mthma}{Theorem A}
\newtheorem*{mthmb}{Theorem B}

\theoremstyle{definition}
\newtheorem{defin}[thm]{Definition}

\theoremstyle{remark}
\newtheorem{remark}[thm]{Remark}
\newtheorem{example}[thm]{Example}

\newcommand{\bth}{\begin{thm}}
\renewcommand{\eth}{\end{thm}}
\newcommand{\bpr}{\begin{prop}}
\newcommand{\epr}{\end{prop}}
\newcommand{\ble}{\begin{lem}}
\newcommand{\ele}{\end{lem}}
\newcommand{\bco}{\begin{cor}}
\newcommand{\eco}{\end{cor}}
\newcommand{\bde}{\begin{defin}}
\newcommand{\ede}{\end{defin}}
\newcommand{\bex}{\begin{example}}
\newcommand{\eex}{\end{example}}
\newcommand{\bre}{\begin{remark}}
\newcommand{\ere}{\end{remark}}
\newcommand{\bcj}{\begin{conj}}
\newcommand{\ecj}{\end{conj}}

\newcommand{\bal}{\begin{aligned}}
\newcommand{\eal}{\end{aligned}}
\newcommand{\beq}{\begin{equation}}
\newcommand{\eeq}{\end{equation}}
\newcommand{\ben}{\begin{equation*}}
\newcommand{\een}{\end{equation*}}

\newcommand{\bpf}{\begin{proof}}
\newcommand{\epf}{\end{proof}}

\def\beql#1{\begin{equation}\label{#1}}

%%%%%%%%%%%%%%  Takuya Matsumoto  %%%%%%%%%%%%%%%%%%%%
\newcommand{\cO}{{\mathcal O}}
\newcommand{\cQ}{{\mathcal Q}}
\newcommand{\sC}{{C\!\!\!\!/\,}}
\newcommand{\fL}{{\mathfrak L}}
\newcommand{\fR}{{\mathfrak R}}
\newcommand{\fQ}{{\mathfrak Q}}
\newcommand{\fS}{{\mathfrak S}}
\newcommand{\fC}{{\mathfrak C}}
\newcommand{\fP}{{\mathfrak P}}
\newcommand{\fK}{{\mathfrak K}}
\newcommand{\fJ}{{\mathfrak J}}
\newcommand{\fB}{{\mathfrak B}}
\newcommand{\fF}{{\mathfrak F}}
\newcommand{\fe}{{\mathfrak e}}
\newcommand{\fh}{{\mathfrak h}}

\newcommand{\anti}{{\rm S}}
\newcommand{\tilh}{{\tilde h}}

\newcommand{\sfR}{{\mathfrak R}\!\!\!\!/\,}
\newcommand{\sfL}{{\mathfrak L}\!\!\!/}
\newcommand{\sfC}{{\mathfrak C}\!\!\!/}
\newcommand{\cT}{\mathcal T}
\newcommand{\Tr}{{\rm Tr}}
\newcommand{\bs}[1]{\boldsymbol{#1}}
\newcommand{\alg}[1]{\mathfrak{#1}}
\newcommand{\el}{\nonumber}
\newcommand{\nln}{\nonumber\\}
\newcommand{\eg}{\widetilde{{\mathfrak g}}}

\newcommand{\ms}{\medskip}

\newtheorem{dfn}{Definition}[section]
\newtheorem{rem}{Remark}[section]
%%%%%%%%%%%%%%%%%%%%%%%%%%%%%%%%%%%%%%%%%%%%%%%%%%%

\title{\Large\bf 
Drinfeld realization of the centrally \\
extended $\mathfrak{psl}(2|2)$ Yangian algebra \\
with the manifest coproducts 
} 

\author{Takuya Matsumoto$^\dagger$ }

\date{} % Start March 2014
\maketitle

\vspace{25 mm}

\begin{abstract}
The Lie superalgebra $\mathfrak{psl}(2|2)$ is recognized as a pretty special one 
in both mathematics and theoretical physics. 
In this paper, we present the Drinfeld realization of the Yangian algebra associated with 
the centrally extended Lie superalgebra $\mathfrak{psl}(2|2)$\,. 
Furthermore, we show that it possesses the Hopf algebra structures, particularly the coproducts. 
The idea to prove the existence of the manifest coproducts is the following. 
Firstly, we shall introduce them to Levendorskii's realization, 
a system of a finite truncation of the Drinfeld generators.
Secondly, we show that Levendorskii's realization is isomorphic to the Drinfeld realization by induction.

% \medskip
%
% Mathematics Subject Classification 2010: 17B37, 17B67

\end{abstract}

%%%\vspace{5 mm}
%%%
%%%{\it Key words:}
%%%

\vfill

\noindent
${}^\dagger$
Department of Applied Physics, 
Faculty of Engineering, 
University of Fukui, \\
$3$-$9$-$1$ Bunkyo, Fukui-shi, Fukui $910$-$8507$, Japan\\
E-mail: {\tt takuyama@u-fukui.ac.jp} 

%\vspace{7 mm}

%\noindent
%School of Mathematics and Statistics\newline
%University of Sydney,
%NSW 2006, Australia\newline
%alexander.molev@sydney.edu.au
%

\setcounter{footnote}{0}
\setcounter{page}{0}
\thispagestyle{empty}

\newpage

\tableofcontents

\section{Introduction and main theorem}
\label{sec:int}
\setcounter{equation}{0}

It is known that the Lie superalgebra $\alg{psl}(2|2)$ has the three-dimensional 
universal central extension \cite{ik:ce}\,, which we denote by $\g$ in this paper, 
\begin{align}
\g:=\alg{psl}(2|2) \oplus \CC^3\,, 
\end{align}
where $\CC^3$ is the three-dimensional abelian Lie algebra. 
As shown in \cite{ik:ce}\,, except the Lie superalgebra $\alg{psl}(2|2)$\,, 
the dimension of the universal central extension of the 
basic classical Lie superalgebra is one at most. 
In this sense, this algebra is an outstanding one. 

\ms

In addition to the mathematical property, this algebra has special interests from physics,
which were discovered by N. Beisert \cite{b:su,b:ab}\,.
He proposed a certain integrable spin-chain model having the symmetry of
the Lie superalgebra $\g$ and succeeded in solving the spectrum problem
of the planar AdS/CFT correspondence in the string theory \cite{b:su}\,.
In this model, the dispersion relation of the excitation on the spin chain, so-called
{\it magnon}, is completely determined by the presentation theory of the Lie superalgebra $\g$\,.
More precisely, the value of {\it energy} is given by the eigenvalue of one of the central
elements $\CC^3$ and the eigenvalues of the other two are the {\it momenta} of
the magnon.
The dispersion relation is mathematically interpreted as the atypical condition for
the Kac module of $\g$\,.
The finite-dimensional irreducible representations of $\alg{gl}(2|2)$ and
$\alg{psl}(2|2)$ are discussed \cite{kkp:fd, ps:fd} and \cite{gqs:tp}\,, respectively,
and that of the Lie superalgebra $\g$ is investigated in detail by \cite{mm:fd}\,.
Another remarkable property of this model is that the S-matrix describing
the scattering of two magnons is completely determined by the symmetry of the Lie algebra $\g$
up to an overall phase factor \cite{b:su}.
This S-matrix satisfies the Yang-Baxter equation.
Thus, this model is integrable.
Furthermore, as pointed out in \cite{b:ab}\,, it is surprising that
this S-matrix also appears in a completely different physical context,
the one-dimensional Hubbard model.
Beisert's S-matrix \cite{b:su} coincides with Shastry's S-matrix for
the one-dimensional Hubbard model \cite{Shastry} up to a particular twist
\footnote{The R-matrix based on the quantum affine algebra
associated with $\alg{sl}(2|2)$ is also proposed by \cite{y:ce}.}.

\ms 

Since it is known that the Hubbard hamiltonian has the infinite-dimensional Yangian
symmetry $\Y(\alg{sl}(2))\oplus \Y(\alg{sl}(2))$ \cite{uk:ys}\,,
it is natural to expect a similar symmetry of the S-matrix.
Indeed, it turned out to be the Yangian symmetry $\Y(\g)$ associated with the centrally
extended Lie superalgebra $\g$ \cite{b:yan}\,,
which includes the Yangian $\Y(\alg{sl}(2))\oplus \Y(\alg{sl}(2))$ as its subalgebra.
The superYangian algebra $\Y(\g)$ is the main object to consider in this article.

\ms 

The Yangian algebra $\Y(\g)$ in \cite{b:yan} is given by the formulation of  
the Drinfeld's original proposal, 
so-called the {\it Drinfleld's first realization} \cite{dri:ybe,dri:qg}\,. 
However, this formulation is not preferable for studying the representation theory. 
For this purpose, the {\it Drinfeld realization}\footnote{
This is also referred as to {\it Drinfleld's second realization}, 
in contrast to {\it Drinfleld's first realization}. 
In Drinfeld's original article,  it is called the {\it new realization} \cite{dri:new}. }
is more appropriate presentation  
since all Cartan generators are explicitly included in the defining relations \cite{dri:new}. 
Another important formulation is so-called {\it RTT-realization}, 
which is the formulation based on the Yang-Baxter equations.  
With this presentation, the coproducts structures are more transparent. 
The RTT-realization of $\Y(\g)$ is proposed in \cite{bd:rtt}\,. 
However, the closure of the algebra is not wholly proved, and consequently, 
the relations to the other realizations are not very obvious.  

\ms

%%%%%%%%%%%%%%%%%%%%%%%%%%%%%%%%%%%%%%%

Thus it is desirable to consider the Drinfeld realization of the Yangian $\Y(\g)$\,.
For this point, there is already a proposal by Spill-Torrielli \cite{ST}. 
However, the Hopf algebra structure of $\Y(\g)$ is not proved completely. 
In particular, it is far from evident that the coproducts 
$\Delta : \Y(\g)\to \Y(\g)\otimes \Y(\g)$
are compatible with the defining relations.  
Regarding the defining relations in \cite{ST}\,, they are associated with  
the Dynkin diagram $\otimes$-$\otimes$-$\otimes$, where all Chevalley generators are odd. 
With these generators, it is hard to see the decomposition of the module
of $\Y(\g)$ into that of $\Y(\alg{sl}(2))\oplus \Y(\alg{sl}(2))$\,. 
Therefore, in this paper, we will formulate the Drinfeld realization of $\Y(\g)$
based on the Dynkin diagram O-$\otimes$-O, where we only have the odd simple root
for the second node and the other two are even. 
In particular, the first and third even roots are associated with the Chevalley 
generators of the two $\Y(\alg{sl}(2))$ subalgebras respectively. 

\ms 

The centrally extended Lie superalgebra $\g$
is also obtained from 
the exceptional Lie superalgebra $D(2,1;\alpha)$ 
by a certain limit for the parameter $\al$\,. 
The Drinfeld realization of the quantum affine algebra associated with 
$D(2,1;\alpha)$ is proposed by \cite{HSTY} using the notation of the 
Weyl groupoid. Hence, it would be interesting to see the relation between 
our Drinfeld realization of the Yangian and that of the quantum affine algebra.

\ms

To introduce the coproduct structures to the Drinfeld realization of $\Y(\g)$\,, 
we shall start with Levendorskii's realization \cite{lev:gen}\,,
which we denote by $\Y_L(\g)$\,. 
In this formulation we only need the finite number of the Drinfeld generators; 
\begin{align}
h_{i,0},~ x^\pm_{i,0}, ~h'_{i,1}, ~x^\pm_{i,1} \quad \text{with}\quad 
i=1,2,3\,, \qquad \text{and}\qquad P_0^\pm, ~P_1^\pm \,,
\end{align}
where the generators $h_{i,0}, x^\pm_{i,0}$ and $P_0^\pm$ are {\it degree} zero, 
and nothing but the Chevalley generators of the Lie superalgebra $\g$\,. 
The extended central elements are $P_0^\pm$\,.
On the other hand, the generators $h'_{i,1}, x^\pm_{i,1}$ and $P_1^\pm$
are degree one Yangian generators in $\Y(\g)$\,, which have the non-trivial 
coproduct structures. 
Therefore, Levendorskii's realization \cite{lev:gen} is a truncation of the 
Drinfeld realization up to degree zero and one generators. 
With this realization, we can show that the coproducts are indeed compatible
with the defining relations. 
Next, we will prove that the Levendorskii's realization $\Y_L(\g)$ is isomorphic to 
the Drinfeld realization$\Y(\g)$\,, where we have the infinite numbers of generators; 
\begin{align}
h_{i,r},~ x^\pm_{i,r}, \quad \text{with}\quad 
i=1,2,3\,, \qquad \text{and}\qquad P_r^\pm \,,
\end{align}
where the degrees are all non-negative integers $r=0,1,2,\cdots$. 
The isomorphism is explicitly proved by the induction with respect to the 
degree $r$\,. 
Now we are ready to state our main theorem as below; 
\begin{mthm}
The Yangian $\Y_D(\g)$ is generated by 
$h_{i,r}, x^\pm_{i,r}$ with $i=1,2,3$  and the central elements $P^\pm_r$
with $r=0,1,2,\cdots$, satisfying the following relations ; 
\begin{align}
[h_{i,r},h_{j,s}]&=0 
\\
[x^+_{i,r},x^-_{j,s}]&=\delta_{ij}h_{i,r+s}
\\
[h_{i,0},x^\pm_{j,r}]&=\pm a_{ij}\xi^\pm_{j,r}
\\
[h_{i,r+1},x^\pm_{j,s}]-[h_{i,r},x^\pm_{j,s+1}]&=\pm\frac{1}{2}a_{ij}\{h_{i,r},x^\pm_{j,s}\}
\\
[x^\pm_{i,r+1}, x^\pm_{j,s}]-[x^\pm_{i,r}, x^\pm_{j,s+1}]&=\pm\frac{1}{2} a_{ij}\{x^\pm_{i,r}, x^\pm_{j,s}\}
\\ 
[x_{1,r}^\pm,x_{3,s}^\pm]=[x_{2,r}^\pm, x_{2,s}^\pm]&=0 
%\qquad \text{for} \qquad i+j=4 
% [x_{i,r}^\pm,x_{j,s}^\pm]&=0 \qquad \text{for} \qquad i+j=4 
\\
\bigl[x^\pm_{j,r},[x^\pm_{j,s},x^\pm_{2,t}]\bigr]+\bigl[x^\pm_{j,s},[x^\pm_{j,r},x^\pm_{2,t}]\bigr]&=0
\qquad \text{for} \qquad j=1,3 
\\
%[x_{2,r}^\pm, x_{2,s}^\pm]&=0
%\\
\bigl[[x^\pm_{1,r},x^\pm_{2,0}],[x^\pm_{3,s},x^\pm_{2,0}]\bigr]&=P^\pm_{r+s} 
% \label{extP}
% \\
% \bigl[[\xi^-_{1,r},\xi^-_{2,0}],[\xi^-_{3,s},\xi^-_{2,0}]\bigr]&=K_{r+s} 
% \label{extK}
\end{align} 
where the indices of subscript $r,s,t$ runs non-negative integers 
and $i,j=1,2,3$. 
The symmetrized Catran matrix is given by 
\begin{align}
(a_{ij})=
%\begin{pmatrix} 2 & -1 & 0 \\ -1& 0 & 1 \\ 0 & 1 & -2   \end{pmatrix}
\left(\!
\begin{array}{rrr}
2 & -1 & 0 \\
-1& 0 & 1 \\
0 & 1 & -2  
\end{array} 
\right)
\label{car-mat}
\end{align}
The generators $x^\pm_{2,r}$ are odd, and the others are even.
All commutators $[~,~]$ (and anti-commutators $\{~,~\}$) are generalized as
super (anti-)commutators, respectively. 
The Yangian $\Y_D(\g)$ has a Hopf algebra structure with the 
coproducts $\Delta : \Y_D(\g)\to \Y_D(\g)\otimes \Y_D(\g)$\,, counits 
$\ep : \Y_D(\g)\to \CC$\,,
and antipodes ${\rm S} : \Y_D(\g)\to\Y_D(\g)$\,,
which are given in {\rm Proposition \ref{prop:copro}}. 
\end{mthm}

\ms 

In particular, the coproducts structure allows us to prove the PBW theorem 
for the Yangian $\Y_D(\g)$ along the idea of \cite{lev:pbw}\,. 
We will report the PBW theorem for $\Y_D(\g)$ in the near future. 
%in \cite{M-PBW}\,.  

\ms 

This paper is organized as follows.
In Section \ref{sec:sl22}\,, we  review the centrally extended Lie superalgebra $\g$\,. 
In Section \ref{sec:YL}\,, we define the Yangian algebra $\Y_L(\g)$ based on Levendorskii's realization
and prove the Hopf algebra structures, in particular, the coproducts. 
In Section \ref{sec:DII}, we will construct Drinfeld realization of the Yangian $\Y_D(\g)$
from Levendorskii's realization $\Y_L(\g)$\,, inductively. 
Consequently, the Hopf algebra structures of $\Y_D(\g)$ are induced from $\Y_L(\g)$\,. 
In Appendix \ref{app:D1}\,, we explain how our Yangian $\Y_L(\g)$ is relating to 
the Drinfeld's first realization.

%Using the coproducts, we shall prove the PBW theorem for the Yangian $\Y(\g)$ in Section \ref{sec:PBW}. 

%%%%%%%%%%%%%%%%%%%%%%%%%%%%%%%%%%%%%%%%%%
\section{Central extensions of the Lie superalgebra $\alg{psl}(2|2)$}
\label{sec:sl22}

We start to define the central extensions of the Lie superalgebra $\alg{psl}(2|2)$\,. 
\begin{defin}
\label{def:lie}
The centrally extended Lie superalgebra 
$\g=\alg{psl}(2|2)\oplus \CC^3$ over $\CC$ has the generators 
%universal enveloping algebra $\U(\g)$ over $\CC$ 
%associated with the Lie superalgebra $\g:=\alg{psl}(2|2)\oplus \CC^3$ 
$h_{i,0}\,, x^\pm_{i,0}$ with $i=1,2,3$ and the central elements $P_0^\pm$\,, 
and they satisfy the following relations; 
\begin{align}
[h_{i,0},h_{j,0}]&=0 \label{lie-1}
\\
[h_{i,0},x^\pm_{j,0}]&=\pm a_{ij}x^\pm_{j,0} \label{lie-2}%\label{s1_hx}
\\
%[h'_{i,1},x^\pm_{j,0}]&=\pm a_{ij}x^\pm_{j,1} %\label{s1_h1x}
%\\
[x^+_{i,0},x^-_{j,0}]&=\delta_{ij}h_{i,0} \label{lie-3}%\label{s1_xpxm0}
\\
%[x^+_{i,1},x^-_{j,0}]&=\delta_{ij}h_{i,1}\equiv \delta_{ij}(h'_{i,1}+\tfrac{1}{2}h^2_{i,0}) %\label{s1_xpxm1}
%\\
%[x^\pm_{i,1},x^\pm_{j,0}]-[x^\pm_{i,0},x^\pm_{j,1}]&=\pm\frac{a_{ij}}{2}\{x^\pm_{i,0},x^\pm_{j,0}\} %\label{s1_xxxx}
%\\
%[h'_{i,1},[x^+_{j,1},x^-_{j,1}]]&=0 
%\\
[x^\pm_{2,0},x^\pm_{2,0}]=
[x^\pm_{1,0},x^\pm_{3,0}]&=0 \label{lie-4}
\\
[x^\pm_{i,0},[x^\pm_{i,0},x^\pm_{2,0}]]&=0 \qquad \text{for}\qquad i=1,3 %\label{s1_serre2} 
\label{lie-5}
\\
[[x^\pm_{1,0},x^\pm_{2,0}],[x^\pm_{3,0},x^\pm_{2,0}]]&=P^\pm_0\,.  \label{lie-last}
%\\
%[[x^\pm_{1,1},x^\pm_{2,0}],[x^\pm_{3,0},x^\pm_{2,0}]]&=P^\pm_1 
\end{align}
The $\ZZ_2$-grading $p : \g\to\ZZ_2$ is defined by setting 
$p(x^\pm_{2,0})=1$ and $p(\text{others})=0$\,. 
The symmetrized Cartan matrix $a_{ij}$ is given in \eqref{car-mat}.  
\end{defin}

\begin{remark}
The Lie superalgebra $\g$ is $17$-dimension over $\CC$\,. 
The abelian Lie subalgebra $\CC^3$ is spanned by the two 
central elements $P_0^+\tss,\,P^-_0$ and 
\begin{align}
C_0:=-\frac{1}{2}h_{1,0}-h_{2,0}-\frac{1}{2}h_{3,0} \,. 
\label{C0}
\end{align}
Noting that the generator $C_0$ is the central element of the Lie superalgebra 
$\alg{sl}(2|2)$\,. 
\end{remark}

\begin{remark}
Definition \ref{def:lie} is associated with the Dynkin diagram 
O-$\otimes$-O rather than 
$\otimes$-$\otimes$-$\otimes$ in \cite{ST}\,, 
where the nodes O and $\otimes$ denote even and odd roots,
respectively. 
In particular, the sets of generators $\{x_{1,0}^\pm\tss,\,h_{1,0}\}$ and 
$\{x_{3,0}^\pm\tss,\,h_{3,0}\}$ span  
% are the Chevalley generators of 
two $\alg{sl}(2)$ subalgebras of $\g$\,. 
\end{remark}

\ms

The universal enveloping algebra $\U(\g)$ of $\g$ has the Hopf algebra structures 
with the coproducts $\Delta : \U(\g)\to \U(\g)\otimes \U(\g)$\,, 
counits $\ep : \U(\g)\to \CC$ 
and antipodes ${\rm S} : \U(\g)\to \U(\g)$ defined by\footnote{
Here we have omitted the abelian braiding factor introduced in \cite{b:su,b:ab}. 
} 
\begin{align}
\Delta(a)=a\otimes 1+1\otimes a\,, \quad 
\ep(a)=0\,, \quad {\rm S}(a)=-a
\end{align}
for $a\in \U(\g)$\,. 
Since $\U(\g)$ is $\ZZ_2$-graded algebra, to see the homomorphism 
of the coproducts, we need to take into account the following rule
\begin{align}
(a\otimes b)(c\otimes d)=(-1)^{p(b)p(c)}\,ac\otimes bd\,.  
\end{align}
for $a,b,c,d\in \U(\g)$\,. 

\ms 

For the later arguments, it would be convenient to introduce the matrix  
presentation of the Lie superalgebra $\g$ to express the non-simple root generators.  
The matrix presentation of $\g$ is a quotient of the 
centrally extended general linear Lie superalgebra $\alg{gl}(2|2)$\,.  
Let $\alg{gl}(2|2)\oplus\CC^2$ be a Lie superalgebra generated by 
the standard basis $E_{ij}$ with $1\leq i,j\leq4$ and the central elements $P,K$, 
satisfying the relations; 
\begin{align}
[E_{ij},E_{kl}]=\de_{kj}E_{il}-\de_{il}E_{kj}(-1)^{(p(i)+p(j))(p(k)+p(l))}
+\bar{\ep}_{ik}\ep_{jl}P+\ep_{ik}\bar{\ep}_{jl}K\,,  
\label{lie-mat}
\end{align}
where the square bracket is a super-commutator and the constants 
$\ep_{ij}$ and $\bar{\ep}_{ij}$ are zero except for the values 
\begin{align}
\ep_{12}=-\ep_{21}=1 \qquad\text{and}\qquad \bar\ep_{34}=-\bar\ep_{43}=1\,. 
\end{align}
The $\ZZ_2$-grading of the indices is defined by $p(1)=p(2)=0$\,,  
$p(3)=p(4)=1$ and that of the generators are $p(E_{ij})=p(i)+p(j)$\,. 
Consider a quotient of the Lie superalgebra $\alg{gl}(2|2)\oplus\CC^2$ by the relation
\begin{align}
E_{11}+E_{22}-E_{33}-E_{44}=0\,, 
\end{align}
and denote this algebra by $\tilde\g=\alg{sl}(2|2)\oplus\CC^2$\,. 
This notation is consistent with the following direct sum decomposition 
\begin{align}
\alg{gl}(2|2)=\alg{sl}(2|2)\oplus \CC(E_{11}+E_{22}-E_{33}-E_{44})\,. 
\end{align}

\ms 

\begin{prop} The centrally extended Lie super algebra $\g$ in Definition \ref{def:lie}
is isomorphic to the Lie superalgebra $\tilde\g$\,. 
The isomorphism $\g\to\tilde\g$ is given by 
\begin{align}
x^+_{1,0}&\mapsto E_{21}  & x^-_{1,0}&\mapsto E_{12}  
& h_{1,0}&\mapsto E_{22} - E_{11} 
\nonumber \\ 
x^+_{2,0}&\mapsto  -E_{32}    & x^-_{2,0} &\mapsto E_{23}   
& h_{2,0}&\mapsto -E_{33} - E_{22}  
\el \\
x^+_{3,0}&\mapsto  -E_{43}  &  x^-_{3,0}&\mapsto  E_{34}  
& h_{3,0}&\mapsto  -E_{44} + E_{33} 
\el \\
% E_{31}&=-x^+_{4,0} = [x^+_{1,0}, x^+_{2,0}]  & E_{13} &=x^-_{4,0} =[x^-_{1,0}, x^-_{2,0}] \el \\
% E_{42}&=-x^+_{5,0} = -[x^+_{2,0}, x^+_{3,0}]  & E_{24} &=x^-_{5,0} =[x^-_{2,0}, x^-_{3,0}] \el \\ 
% E_{41}&=-x^+_{6,0} = [[x^+_{1,0}, x^+_{2,0}],x^+_{3,0}]   & E_{14} &=x^-_{6,0} =[[x^-_{1,0}, x^-_{2,0}],x^-_{3,0}] \el \\
P^+_0 &\mapsto P  &  P^-_0 &\mapsto   -K
%& C_0&=\tfrac{1}{2}(E_{ll}+E_{\mu\mu}) ~. 
\label{lie-iso}
\end{align}
\end{prop}

\begin{proof}It is straightforward to see that the map in \eqref{lie-iso} 
preserves the defining relations \eqref{lie-1}--\eqref{lie-last} by using the relations
\eqref{lie-mat}. Thus, it is a homomorphism. 
The surjectivity is also obvious. Precisely, the non-simple roots are given by 
\begin{align}
[x_{1,0}^+,x_{2,0}^+]&\mapsto  E_{31}  & [x_{1,0}^-,x_{2,0}^-]&\mapsto  E_{13}   
\nonumber \\ 
-[x_{2,0}^+,x_{3,0}^+]&\mapsto  E_{42}  & [x_{2,0}^-,x_{3,0}^-]&\mapsto  E_{24} 
\nonumber \\ 
[[x_{1,0}^+,x_{2,0}^+],x_{3,0}^+]&\mapsto  E_{41}  
& [[x_{1,0}^-,x_{2,0}^-],x_{3,0}^-]&\mapsto  E_{14} \,. 
\end{align}
To prove the injectivety, we need to see that $\dim\g\leq \dim\tilde\g=17$\,.  
Setting 
\begin{align}
&\Gamma=\{\tss x_{i,0}^\pm\tss,\,h_{i,0}\tss,\,P\tss,\,K\tss\}_{i=1,2,3}\,,  \el \\
&\Gamma'=\{\tss [x_{1,0}^\pm,x_{2,0}^\pm,]\tss,\,
[x_{2,0}^\pm,x_{3,0}^\pm]\tss,\,[x_{1,0}^\pm,[x_{2,0}^\pm,x_{3,0}^\pm]]
\tss \}\cup\Gamma\,, \el \\
\text{and}\quad &V={\rm span}_\CC\{X|X\in \Gamma' \}\subset \g \,,   
\end{align}
%$\Gamma=\{\tss x_{i,0}^\pm\tss,\,h_{i,0}\tss,\,P\tss,\,K\tss\}_{i=1,2,3}$\,,  
%$\Gamma'=\{\tss [x_{1,0}^\pm,x_{2,0}^\pm,]\tss,\,
%[x_{2,0}^\pm,x_{3,0}^\pm]\tss,\,[x_{1,0}^\pm,[x_{2,0}^\pm,x_{3,0}^\pm]]
%\tss \}\cup\Gamma$
%and $V={\rm span}_\CC\{X|X\in \Gamma' \}\subset \g$ \,.  
it is sufficient to show $\g=V$\,. 
To confirm this, it is enough to check $[\Gamma,V]\subset V$\,. 
Since it is obvious that $[\Gamma,\Gamma]\subset V$\,, 
all we need to see are the relations 
$[\Gamma,Y]\subset V$ for $Y\in 
\{\tss [x_{1,0}^\pm,x_{2,0}^\pm,]\tss,\,
[x_{2,0}^\pm,x_{3,0}^\pm]\tss,\,[x_{1,0}^\pm,[x_{2,0}^\pm,x_{3,0}^\pm]] \}$\,. 
We can verify this by direct computations 
using the relations \eqref{lie-1}--\eqref{lie-last}. 
This completes the proof. 
\end{proof}

\begin{remark} The defining relation \eqref{lie-mat} is almost same with that
of $\alg{gl}(2|2)$ except for the following nontrivial relations 
\begin{align}
&[E_{13},E_{24}]= -[E_{23},E_{14}]=K\,, \\
&[E_{31},E_{42}]= -[E_{32},E_{41}]=P\,. 
\end{align}
Through the isomorphism in \eqref{lie-iso}, the above relations correspond to 
the extended Serre relations in \eqref{lie-last}.   
See \cite{dob:serre} for the relating arguments.  
\end{remark}

%%%%%%%%%%%%%%%%%%%%%%%%%%%%%%%%%%
\section{Levendorskii's realization of the Yangian $\Y_L(\g)$ }
\label{sec:YL}
\setcounter{equation}{0}

This section defines the Yangian algebra $\Y_L(\g)$ associated with the Lie superalgebra $\g$ by Levendorskii's realization.
In Subsection \ref{subsec:DII}\,, we define $\Y_L(\g)$ (Definition \ref{def:yan})\,.
In Subsection \ref{subsec:hopf}\,, we prove the hopf algebra structures, especially the coproducts (Proposition \ref{prop:copro})\,.

\subsection{Definition of the Yangian $\Y_L(\g)$}
\label{subsec:DII}

Let us now define the Yangian algebra $\Y_L(\g)$ associated with $\g$\,. 
\begin{defin}
\label{def:yan}
The Yangian $\Y_L(\g)$ associated with the Lie superalgebra $\g$ is 
generated by $h_{i,0}, ~x^\pm_{i,0}, ~\tilh_{i,1},~ x^\pm_{i,1}$ 
with $i=1,2,3$ and the central elements $P_0^\pm, ~P_1^\pm$\,.
They satisfy the relations \eqref{lie-1}--\eqref{lie-last} and 
\begin{align}
%[h_{i,0},h_{j,0}]&=
[\tilh_{i,1},h_{j,0}]&=0 
\label{yan-1}\\
[\tilh_{i,1}, \tilh_{j,1}]&=0 \label{deg-2}
\\
%[h_{i,0},x^\pm_{j,0}]&=\pm a_{ij}x^\pm_{j,0} %\label{s1_hx}
%\\
[\tilh_{i,1},x^\pm_{j,0}]&=\pm a_{ij}x^\pm_{j,1} \label{h1x0}
\\
%[x^+_{i,0},x^-_{j,0}]&=\delta_{ij}h_{i,0} %\label{s1_xpxm0}
%\\
[x^+_{i,1},x^-_{j,0}]&=\delta_{ij}h_{i,1} %\equiv \delta_{ij}(\tilh_{i,1}+\tfrac{1}{2}h^2_{i,0})
\label{yan-tilh}
\\
[x^\pm_{i,1},x^\pm_{j,0}]-[x^\pm_{i,0},x^\pm_{j,1}]
&=\pm\frac{1}{2}a_{ij} \{x^\pm_{i,0},x^\pm_{j,0}\} 
\label{x1x0}
\\
[x^\pm_{2,1},x^\pm_{2,0}]
&=0
\label{xx20} \\
%[\tilh_{i,1},[x^+_{j,1},x^-_{j,1}]]&=0 
%\label{deg-3}
%\\
[\tilh_{j,1},[x^+_{j,1},x^-_{j,1}]]&=0 \qquad\text{for}\qquad  j=1,3
\label{deg-31} \\ %\nln
[\tilh_{1,1},[x^+_{2,1},x^-_{2,1}]]&=0 
\label{deg-32}  \\ %\nln 
%[x^\pm_{2,0},x^\pm_{2,0}]=[x^\pm_{1,0},x^\pm_{3,0}]&=0 
%\\
%[x^\pm_{i,0},[x^\pm_{i,0},x^\pm_{2,0}]]&=0 \qquad \text{for}\qquad i=1,3 %\label{s1_serre2} 
%\\
%[[x^\pm_{1,0},x^\pm_{2,0}],[x^\pm_{3,0},x^\pm_{2,0}]]&=P^\pm_0 
%\\
[[x^\pm_{1,1},x^\pm_{2,0}],[x^\pm_{3,0},x^\pm_{2,0}]]&=P^\pm_1 
\label{yan-serre}
\end{align}
where the generator $h_{i,1}$ in \eqref{yan-tilh} is defined by  
\begin{align}
h_{i,1}= \tilh_{i,1}+\frac{1}{2}(h_{i,0})^2\,. 
\end{align}
The $\ZZ_2$-grading $p : \Y_L(\g)\to\ZZ_2$ is defined by setting 
$p(x^\pm_{2,0})=p(x^\pm_{2,1})=1$ and $p(\text{others})=0$\,. 
The symmetrized Cartan matrix ($a_{ij})$ is given in \eqref{car-mat}.    
\end{defin}

\begin{remark}
This realization of Yangian, which is associated with a finite-dimensional complex simple Lie algebra, was originally proposed by Levendorskii \cite{lev:gen}.
Its generalization to the Yangians of the Lie superalgebras $A(m,n)=\alg{sl}(m+1|n+1)$ is
given by \cite{stu:ann}. 
Compared with the Yangian of the Lie superalgebra
$\alg{sl}(2|2)$\,, the Serre relations \eqref{yan-serre} are extended.
\end{remark}

Here, let us introduce a terminology {\it degree}.
The degree of the generators of $\Y(\g)$ is the second index of the subscript.
For example, $\deg(x_{2,0}^+)=0$ and $\deg(x_{2,1}^+)=1$\,.
The degree of a monomial is the sum of the degree of each generator.
The degree of a polynomial is the highest degree of the monomials
included in the polynomial.
The degree of tensor product of monomials is defined by the sum of the degree of each tensor factor.
The degree of a polynomial of tensor products is given by the highest degree of the tensor products included in the polynomial.

\ms

Set 
\begin{align}
\Y_L(\g)_{r}:=\{~x\in \Y_L(\g) ~|~\deg(x)\leq r~\}\,. 
\label{def_filt}
\end{align}
Then, the Yangian $\Y(\g)$ has a natural filtration with respect to the degree, 
\begin{align}
\{0\}=\Y_L(\g)_{-1}\subset 
\Y_L(\g)_{0}\subset \Y_L(\g)_{1}\subset  \Y_L(\g)_{2}\subset\cdots\,. 
\label{filt}
\end{align}
%This filtration satisfies 
%\begin{align}
%&\bigcup_{r=-1}^\infty Y_L(\g)_r=Y_L(\g)\,, \qquad \bigcap_{r=-1}^\infty Y_L(\g)_r=\CC\,, 
%\nln
%\end{align}

\begin{remark}
The universal enveloping algebra $\U(\g)$ in Definition \ref{def:lie} 
is a subalgebra of $\Y_L(\g)$ and identified with  
$\Y_L(\g)_{0}$\,. 
\end{remark}

\subsection{Hopf algebra structures of  $\Y_L(\g)$}
\label{subsec:hopf}

We show the Hopf algebra structures of $\Y_L(\g)$ in this subsection. 

\begin{prop} 
\label{prop:copro}
The Yangian $\Y_L(\g)$ has the Hopf algebra structures with the coproducts 
$\Delta : \Y_L(\g)\to \Y_L(\g)\otimes \Y_L(\g)$ given by 
\begin{align}
\Delta(x_{i,0}^\pm)&= x_{i,0}^\pm\otimes 1 + 1\otimes x_{i,0}^\pm \el \\
\Delta(h_{i,0})&= h_{i,0}\otimes 1 + 1\otimes h_{i,0} \qquad (\,i=1,2,3\,) \el \\
\Delta(P_{0}^\pm)&=P_{0}^\pm\otimes1+1\otimes P_{0}^\pm \el \\  
\Delta(x_{1,1}^+)&= x_{1,1}^+\otimes 1 + 1\otimes x_{1,1}^+ + h_{1,0}\otimes x_{1,0}^+ 
- \sum_{\mu=3,4} E_{2\mu}\otimes E_{\mu1} \el\\ 
\Delta(x_{2,1}^+)&= x_{2,1}^+\otimes 1 + 1\otimes x_{2,1}^+ + h_{2,0}\otimes x_{2,0}^+ 
+ E_{12}\otimes E_{31} + E_{34}\otimes E_{42}- E_{14}\otimes P_0^+  \el\\ 
\Delta(x_{3,1}^+)&= x_{3,1}^+\otimes 1 + 1\otimes x_{3,1}^+ + h_{3,0}\otimes x_{3,0}^+ 
- %{\textstyle \sum_{l=1,2}}\,
\sum_{l=1,2}E_{l3}\otimes E_{4l} \el\\
\Delta(\tilh_{1,1})&= \tilh_{1,1}\otimes 1+1\otimes \tilh_{1,1} -2 x_{1,0}^-\otimes x_{1,0}^+ 
+ \sum_{\mu=3,4} \bigl(E_{1\mu}\otimes E_{\mu1}- E_{2\mu}\otimes E_{\mu2}\bigr) \el\\ 
\Delta(\tilh_{2,1})&= \tilh_{2,1}\otimes 1+1\otimes \tilh_{2,1} 
+E_{12}\otimes E_{21}- E_{13}\otimes E_{31}+E_{34}\otimes E_{43}+ E_{24}\otimes E_{42} 
-P_0^-\otimes P_0^+ \el\\ 
\Delta(\tilh_{3,1})&= \tilh_{3,1}\otimes 1+1\otimes \tilh_{3,1} +2 x_{3,0}^-\otimes x_{3,0}^+ 
+\sum_{l=1,2} \bigl( E_{l3}\otimes E_{3l}- E_{l4}\otimes E_{4l} \bigr)\el\\ 
\Delta(x_{1,1}^-)&= x_{1,1}^-\otimes 1 + 1\otimes x_{1,1}^- + x_{1,0}^-\otimes h_{1,0} 
-\sum_{\mu=3,4}  E_{1\mu}\otimes E_{\mu2} \el\\ 
\Delta(x_{2,1}^-)&= x_{2,1}^-\otimes 1 + 1\otimes x_{2,1}^- + x_{2,0}^-\otimes h_{2,0} 
- E_{13}\otimes E_{21} - E_{24}\otimes E_{43} -P_0^-\otimes E_{41}   \el\\ 
\Delta(x_{3,1}^-)&= x_{3,1}^- \otimes 1 + 1\otimes x_{3,1}^- + x_{3,0}^- \otimes h_{3,0} 
+\sum_{l=1,2}  E_{l4}\otimes E_{3l} \el \\
\Delta(P_{1}^+)&=P_{1}^+\otimes1+1\otimes P_{1}^+ -2 C_0\otimes P_0^+ \el \\
\Delta(P_{1}^-)&=P_{1}^-\otimes1+1\otimes P_{1}^- -2 P_0^-\otimes C_0 ~,  
\label{coprod}
\end{align}
the counits $\ep:\Y_L(\g)\to\CC$\,,  
\begin{align}
\ep(X)=0 \qquad \text{for }\qquad X\in \Y_L(\g)\,,   
\end{align}
and the antipodes ${\rm S} : \Y_L(\g)\to\Y_L(\g)$\,,   
\begin{align}
{\rm S}(x_{i,0}^\pm)&= -x_{i,0}^\pm\,, \quad %\el \\
{\rm S}(h_{i,0})= - h_{i,0}\,, \quad (\,i=1,2,3\,)\,, \quad %\el \\
{\rm S}(P_{0}^\pm)= -P_{0}^\pm \el \\  
{\rm S}(x_{1,1}^+)&= -x_{1,1}^+ + h_{1,0} x_{1,0}^+ 
- \sum_{\mu=3,4} E_{2\mu} E_{\mu1} \el\\ 
{\rm S}(x_{2,1}^+)&= -x_{2,1}^+ + h_{2,0} x_{2,0}^+ 
+ E_{12} E_{31} + E_{34} E_{42}- E_{14} P_0^+  \el\\ 
{\rm S}(x_{3,1}^+)&= -x_{3,1}^+ + h_{3,0} x_{3,0}^+ 
-\sum_{l=1,2}E_{l3} E_{4l} \el\\
{\rm S}(\tilh_{1,1})&= -\tilh_{1,1}  -2 x_{1,0}^- x_{1,0}^+ 
+ \sum_{\mu=3,4} \bigl(E_{1\mu} E_{\mu1}- E_{2\mu} E_{\mu2}\bigr) \el\\ 
{\rm S}(\tilh_{2,1})&= -\tilh_{2,1}  
+E_{12} E_{21}- E_{13} E_{31}+E_{34} E_{43}+ E_{24} E_{42} 
-P_0^- P_0^+ \el\\ 
{\rm S}(\tilh_{3,1})&= -\tilh_{3,1} +2 x_{3,0}^- x_{3,0}^+ 
+\sum_{l=1,2} \bigl( E_{l3} E_{3l}- E_{l4} E_{4l} \bigr)\el\\ 
{\rm S}(x_{1,1}^-)&= -x_{1,1}^-  + x_{1,0}^- h_{1,0} 
-\sum_{\mu=3,4}  E_{1\mu} E_{\mu2} \el\\ 
{\rm S}(x_{2,1}^-)&= -x_{2,1}^- + x_{2,0}^- h_{2,0} 
- E_{13} E_{21} - E_{24} E_{43} -P_0^- E_{41}   \el\\ 
{\rm S}(x_{3,1}^-)&= -x_{3,1}^- + x_{3,0}^-  h_{3,0} 
+\sum_{l=1,2}  E_{l4} E_{3l} \el \\
{\rm S}(P_{1}^+)&= -P_{1}^+ -2 C_0 P_0^+ \el \\
{\rm S}(P_{1}^-)&=-P_{1}^- -2 P_0^- C_0 \,.   
\label{anti}
\end{align} 
\end{prop}

\begin{remark}
The coproducts listed in \eqref{coprod} are essentially translated from \cite{b:yan}, where Drinfeld's first realization is adopted.
\end{remark}

% restart from here 

To prove Proposition \ref{prop:copro}, the following lemmas are useful. 

\begin{lem}
A map $\rho : \Y_L(\g)\to\Y_L(\g)$ defined by 
\begin{align}
x^+_{1,r}&\mapsto x^-_{3,r}   & x^-_{1,r}&\mapsto -x^+_{3,r}   
& h_{1,r}& \mapsto h_{3,r} \el \\
x^+_{2,r}&\mapsto -x^-_{2,r}   & x^-_{2,r}&\mapsto -x^+_{2,r}  & 
h_{2,r}& \mapsto h_{2,r}  \el \\
x^+_{3,r}&\mapsto -x^-_{1,r}  & x^-_{3,r}&\mapsto x^+_{1,r}  & 
h_{3,r}& \mapsto h_{1,r}  \el \\
P^+_r &\mapsto -P^-_r & P^-_r &\mapsto -P^+_r &
%(r&=0,1\,. )
\label{auto_pm}
\end{align}
with $r=0,1$ is an automorphism of the Yangian $\Y_L(\g)$\,. 
The map $\rho$ preserve the $\ZZ_2$-grading 
\begin{align}
p(\rho(x)) = p(x)\qquad \text{for}\qquad x\in \Y_L(\g)\,,  
\end{align}
and satisfies $\rho^2={\rm id}$\,. 
\end{lem} 

\begin{proof}
By direct computations, it is verified that the map $\rho$ preserves 
all defining relations in \eqref{lie-1}--\eqref{lie-last} 
and \eqref{yan-1}--\eqref{yan-serre}. 
%Since the odd generators are only $x^\pm_{2,r}$ ($r=0,1$)\,, it is obvious that the $\ZZ_2$-grading is invariant under $\rho$\,. 
The properties $\rho^2(x)=x$ and $p(\rho(x)) = p(x)$ for any  $ x\in \Y_L(\g)$ immediately follow from the definition \eqref{auto_pm}. 
\end{proof}

\begin{remark}
It is easy to see that 
the Lie superalgebra $\g$ defined in a matrix form \eqref{lie-mat} enjoys 
an automorphism 
\begin{align}
&E_{ij}\mapsto E_{5-i,5-j} \qquad \text{for} \qquad i,j=1,\cdots,4\,,
\nln
&P\mapsto K\,, \qquad K\mapsto P\,.  
\end{align}
The automorphism $\rho$ in \eqref{auto_pm} 
for the generators with $r=0$ agrees with this map. 
Then, $\rho$  is a natural extension of this to Yangian $\Y_L(\g)$\,.  
\end{remark}

%This map may simplify the proof of the homomorphism of coproducts 
%because of the following proposition. 

\begin{lem}
\label{lem:auto}
The automorphism $\rho : \Y_L(\g) \to \Y_L(\g)$ 
% $\rho\in {\rm Aut}(\Y(\g))$ 
defined in \eqref{auto_pm} is compatible with the coproducts in \eqref{coprod}, 
\begin{align}
(\rho\otimes \rho) \circ \widetilde{\Delta} = \Delta \circ \rho ~,  
\label{compati}
\end{align}
where we denoted the opposite coproducts by 
$\widetilde{\Delta}=\sigma \circ \Delta$ with the graded permutation $\sigma$
defined by $\sigma(x\otimes y)=(-1)^{p(x)p(y)}\tss y\otimes x$ for $x,y\in \Y_L(\g)$. 
\end{lem}

\begin{proof} Set ${\rm G}_r=\{\tss x^\pm_{i,r}\tss,\, h_{i,r}\tss,\, P^\pm_r~|~
i=1,2,3\tss\}\subset \Y_L(\g)$ with $r=0,1$\,. 
It is obvious that \eqref{compati} holds 
for $x\in {\rm G}_0$ since their coproducts are primitive, 
$\De(x)=x\otimes 1 +1\otimes x$ for $x\in {\rm G}_0$\,. 
Then, to prove the lemma, it is sufficient to show 
\begin{align}
(\rho\otimes \rho) \circ \widetilde{\Delta} (x)= \Delta \circ \rho(x)
\qquad \text{for}\qquad x\in {\rm G}_1\,. 
\end{align}
This is confirmed by direct calculations. For example, 
\begin{align}
(\rho\otimes \rho) \circ \widetilde{\Delta} (\tilh_{1,1})
&= (\rho\otimes \rho) \circ\sigma \Bigl(
 \tilh_{1,1}\otimes 1+1\otimes \tilh_{1,1} -2 x_{1,0}^-\otimes x_{1,0}^+ 
\el\\
&\qquad + \sum_{\mu=3,4} 
\bigl(E_{1\mu}\otimes E_{\mu1}- E_{2\mu}\otimes E_{\mu2}\bigr) \Bigr) 
\el \\
&= (\rho\otimes \rho) \Bigl(
1\otimes \tilh_{1,1}+\tilh_{1,1}\otimes1 -2 x_{1,0}^+\otimes x_{1,0}^- 
\el\\
&\qquad - \sum_{\mu=3,4} 
\bigl(E_{\mu1}\otimes E_{1\mu}- E_{\mu2}\otimes E_{2\mu}\bigr) \Bigr) 
\el \\
&=
1\otimes \tilh_{3,1}+\tilh_{3,1}\otimes1 +2 x_{3,0}^-\otimes x_{3,0}^+ 
\el\\
&\qquad - \sum_{\mu=1,2} 
\bigl(E_{l 4}\otimes E_{4l}- E_{l3}\otimes E_{3l}\bigr)  
\el \\
&= \Delta \circ \rho(\tilh_{1,1})\,. 
\end{align}
The computations are similar for the other generators in ${\rm G}_1$\,. 
This completes the proof. 
\end{proof}

We are ready to prove Proposition \ref{prop:copro}. 

\begin{proof}[Proof of Proposition \ref{prop:copro}]
Let $\mu : \Y_L(\g) \otimes \Y_L(\g)\to \Y_L(\g)$ be the product of $\Y_L(\g)$ defined by $\mu(x\otimes y)=xy$ for $x,y\in \Y_L(\g)$\,, and $u : \CC\to \Y_L(\g)$ the unit define by $u(1)=1\in \Y_L(\g)$\,.  
Then, it is straightforward to check that the antipodes in \eqref{anti} satisfy 
the antipode relations, 
\begin{align}
\mu\circ({\rm S}\otimes 1)\circ \Delta (x)
=\mu\circ(1\otimes {\rm S})\circ \Delta (x)
= u\circ \ep(x)=0
\end{align}
for $x\in \Y_L(\g)$\,. 

\ms

The non-trivial part is to show that the coproducts given in \eqref{coprod} is a homomorphism of Yangian $\Y_L(\g)$\,. We shall prove this with respect to each degree of the relations.

\ms

For the degree zero relations \eqref{lie-1}--\eqref{lie-last}, 
since the products of the degree zero generators are primitive, 
\begin{align}
\Delta(x)=x\otimes1+1\otimes x  \qquad \text{for}\qquad x\in \Y_L(\g)_0\,, 
\end{align}
it is easy to see that they are consistent with the relations
\begin{align}
[\Delta(x),\Delta(y)]=\Delta([x,y]) \qquad \text{for}\qquad x,y\in \Y_L(\g)_0\,.  
\end{align}

\ms

The degree one relations are \eqref{yan-1}, \eqref{h1x0} ,\eqref{yan-tilh}, \eqref{x1x0}, \eqref{xx20}, and \eqref{yan-serre}. 
It is not so difficult to see the compatibility of the coproducts with the relations by direct computations.

\ms 

\noindent$\bullet$ Relation \eqref{yan-1}. 
\begin{align}
[\Delta(\tilh_{i,1}),\Delta(h_{j,0})]&=0 \qquad \text{for}\qquad i,j=1,2,3\,. 
\end{align}

\noindent$\bullet$ Relations \eqref{h1x0}.  
After some computations, we can verify 
\begin{align}
[\Delta(\tilh_{i,1}),\Delta(x^+_{j,0})]&=a_{ij}\Delta(x^+_{j,1})
\label{h1xp}
\end{align}
for $i,j=1,2,3$\,. Then, the following relations are obtained from the above 
by using the automorphism $\rho$ in \eqref{auto_pm} and Lemma \ref{lem:auto}\,, 
\begin{align}
[\Delta(\tilh_{i,1}),\Delta(x^-_{j,0})]&=-a_{ij}\Delta(x^-_{j,1})\,.  
\label{h1xm}
\end{align}

\noindent$\bullet$ Relations \eqref{yan-tilh}. 
At first, we compute 
\begin{align}
[\Delta(x^+_{1,1}),\Delta(x^-_{1,0})]=\Delta(h_{1,1})\,, \el \\
[\Delta(x^+_{2,1}),\Delta(x^-_{2,0})]=\Delta(h_{2,1})\,, 
\end{align}
where $h_{i,1}=\tilh_{i,1}+\frac{1}{2}h^2_{i,0}$\,. 
Applying the automorphism $\rho$ in \eqref{auto_pm} for the first relation, 
we have 
\begin{align}
[\Delta(x^+_{3,1}),\Delta(x^-_{3,0})]=\Delta(h_{3,1})\,. 
\end{align}
The remaining task here is to prove 
\begin{align}
[\Delta(x^+_{i,1}),\Delta(x^-_{j,0})]=0 \qquad\text{for}\qquad i\neq j\,. 
\label{xpxm-proof}
\end{align}
Let us assume that $i\neq j$\,. 
Acting $\Delta(\tilh_{i,1})$ and $\Delta(\tilh_{j,1})$ on the relations 
\begin{align}
[\Delta(x^+_{i,0}),\Delta(x^-_{j,0})]=0\,, 
\end{align}
and using \eqref{h1xp}\,, \eqref{h1xm}\,, we have
\begin{align}
a_{ii}[\Delta(x^+_{i,1}),\Delta(x^-_{j,0})]-a_{ij}[\Delta(x^+_{i,0}),\Delta(x^-_{j,1})]&=0 \el \\
a_{ji}[\Delta(x^+_{i,1}),\Delta(x^-_{j,0})]-a_{jj}[\Delta(x^+_{i,0}),\Delta(x^-_{j,1})]&=0 \,. 
\end{align}
Noting that $a_{ii}a_{ij}-a_{ji}a_{jj}\neq 0$ if $i\neq j$\,, we obtain 
\begin{align}
[\Delta(x^+_{i,1}),\Delta(x^-_{j,0})]=[\Delta(x^+_{i,0}),\Delta(x^-_{j,1})]=0 \,. 
%\el \\
%a_{ji}[\Delta(x^+_{i,1}),\Delta(x^-_{j,0})]-a_{jj}[\Delta(x^+_{i,0}),\Delta(x^-_{j,1})]&=0 \,. 
\end{align}
This proves \eqref{xpxm-proof}.

\noindent$\bullet$ Relations \eqref{x1x0}. 
Using the automorphism $\rho$ in \eqref{auto_pm}\,, the relations 
\begin{align}
[\Delta(x^+_{i,1}),\Delta(x^+_{j,0})]-[\Delta(x^+_{i,0}),\Delta(x^+_{j,1})]
&= \frac{1}{2}a_{ij} \{\Delta(x^+_{i,0}),\Delta(x^+_{j,0})\}  \label{x1x0-proof}
\end{align}
immediately imply that 
\begin{align}
[\Delta(x^-_{i,1}),\Delta(x^-_{j,0})]-[\Delta(x^-_{i,0}),\Delta(x^-_{j,1})]
&=-\frac{1}{2}a_{ij} \{\Delta(x^-_{i,0}),\Delta(x^-_{j,0})\}\,.   %\label{x1x0}
\end{align}
In addition, the relations \eqref{x1x0-proof} are symmetric under exchanging 
$i$ and $j$\,.  
Thus, it is sufficient to check the relations with respect to the following six cases ; 
\begin{align}
(i,j)=(1,1), (2,2),(3,3),(1,2),(1,3), ~\text{and}~ (2,3)\,.  
\end{align}
By computing all cases, we prove \eqref{x1x0-proof}. 

\noindent$\bullet$ Relations \eqref{xx20}. 
We can show that 
\begin{align}
[\De(x^+_{2,1}), \De(x^+_{2,0})]&=0\,. 
\end{align}
Applying the automorphism $\rho$ for this, we obtain 
$[\De(x^-_{2,1}), \De(x^-_{2,0})]=0$\,.

\noindent$\bullet$ Relations \eqref{yan-serre}. 
By calculations, we get 
\begin{align}
&\quad \bigl[[\Delta(x^+_{1,1}),\Delta(x^+_{2,0})],
[\Delta(x^+_{3,0}),\Delta(x^+_{2,0})]\bigr] \el \\
&=\bigl[[x_{1,1}^+,x_{2,0}^+]\otimes1 +1\otimes [x_{1,1}^+,x_{2,0}^+]+E_{32}\otimes E_{21}
-E_{11}\otimes E_{31} \el \\
&\quad -\sum_{\mu=3,4}E_{3\mu}\otimes E_{\mu1} -E_{24}\otimes P_0^+ \,,~
E_{42}\otimes 1 +1\otimes E_{42} \bigr]
\el \\
&=P_1^+ \otimes1 +1\otimes P_1^+ -E_{32}\otimes E_{41}
-E_{11}\otimes P_0^+ \el \\
&\quad +E_{32}\otimes E_{41} -E_{33}\otimes P_0^+ -(E_{22}+E_{44})\otimes P_0^+ 
\el\\
&= P_1^+ \otimes1 +1\otimes P_1^+ -(E_{11}+E_{22}+E_{33}+E_{44})\otimes P_0^+ 
\el \\
&= P_1^+ \otimes1 +1\otimes P_1^+ -2C_0 \otimes P_0^+
\el \\
&=\Delta(P^+_1)\,. %\label{yan-serre}
\end{align}
Here $C_0$ is defined in \eqref{C0} and we have used the matrix notation 
\begin{align}
C_0=\frac{1}{2}(E_{11}+E_{22}+E_{33}+E_{44})\,, 
\end{align}
which follows from the isomorphism \eqref{lie-iso}. 
By automorphism $\rho$ in \eqref{auto_pm}, we immediately obtain 
another relation,  
\begin{align}
\bigl[[\Delta(x^-_{1,1}),\Delta(x^-_{2,0})],
[\Delta(x^-_{3,0}),\Delta(x^-_{2,0})]\bigr] =\Delta(P^-_1)\,. 
\end{align}
Thus, the copoducts in \eqref{coprod} are compatible with the defining relations 
of $\Y_L(\g)$ up to degree one. 

\ms

Next, let us move to the degree two relations, \eqref{deg-2}. 

\ms

\noindent $\bullet$ Relations \eqref{deg-2}. 
It is trivial that 
\begin{align}
[\Delta(\tilh_{i,1}), \Delta(\tilh_{i,1})]=0 \qquad\text{for}\qquad i=1,2,3\,. 
\end{align}
Due to the automorphism $\rho$\,, it is sufficient to show that 
\begin{align}
[\Delta(\tilh_{1,1}), \Delta(\tilh_{2,1})]&=0 \label{hh1}\,,\\
[\Delta(\tilh_{1,1}), \Delta(\tilh_{3,1})]&=0 \label{hh2}\,. 
\end{align}

\ms

For convenience, we introduce a notation $\De^{(k)}(x)$ for $0\leq k\leq l\,,x\in \Y_L(\g)_l$ as follows.
Let $\De(x)$ be a coproduct for an element $x\in \Y_L(\g)_l$ with degree $l$\,.
% For the definition, see \eqref{filt}\,.
Assume that it is symbolically  given by  
\begin{align}
\De(x)=\sum_i a_i\otimes b_i \qquad \text{with} \qquad 
a_i\,,b_i \in \Y_L(\g)_l\,. 
\end{align}
For $0\leq k\leq l$\,, we set  
\begin{align}
\De^{(k)}(x)=\sum_{i, {\rm deg}=k} a_i\otimes b_i
~\in \Y_L(\g)_k\otimes \Y_L(\g)_k
\,, %\qquad \text{where} \qquad 
%{\rm deg}(a_i)+{\rm deg}(b_i)=k\,,  
\end{align}
where the summation runs over tensor products of degree $k$\,, ${\rm deg}(a_i)+{\rm deg}(b_i)=k$\,. 
Then, a coproduct $\De(x)$ for $x\in \Y_L(\g)_l$ is uniquely written as 
\begin{align}
\De(x)=\De^{(l)}(x)+\De^{(l-1)}(x)+\cdots +\De^{(1)}(x)+\De^{(0)}(x)\,. 
%\sum_{k=0}^l\De^{(k)}(x)\,. 
\end{align}
%
%Consider  the graded algebra 
%\begin{align}
%{\rm gr}\Y_L(\g)=\bigoplus_{l=1}^\infty {\rm gr}_l\Y_L(\g)
%\qquad \text{with}\qquad 
%{\rm gr}_l\Y_L(\g)=\Y_L(\g)_l/\Y_L(\g)_{l-1}
%\end{align}
%associated with the filtration \eqref{filt}, and  let $\pi_l$ ($l\geq0$) be a natural projection of $\Y_L(\g)_l$ to ${\rm gr}_l\Y_L(\g)$\,. 
%For $0\leq k\leq l$\,, introduce the composition map by 
%\begin{align}
%\De^{(k)}=\sum_{i=0}^k(\pi_{k-i}\otimes \pi_i)\circ \De :
%\Y_L(\g)_l \longrightarrow \sum_{i=0}^k {\rm gr}_{k-i}\Y_L(\g)\otimes {\rm gr}_{i}\Y_L(\g)\,. 
%\end{align}
%Then, we have the direct decomposition 
%\begin{align}
%\De(\Y_L(\g)_l)=\bigoplus_{k=0}^l\De^{(k)}(\Y_L(\g)_l)\,. 
%\end{align}
%
%Introduce the notation for degree $r$ part of the coproducts of $X\in \Y(\g)$ by 
%\begin{align}
%\Delta^{(r)}(X)\equiv \Delta(X) \cap \Y(\g)_r -\Delta(X) \cap \Y(\g)_{r-1} 
%\end{align}
%with $r=0,1$\,. 
%Any coproducts of degree one generators $X_1\in \Y(\g)$ are written as  
%\begin{align}
%\Delta(X_1)=\Delta^{(1)}(X_1)+\Delta^{(0)}(X_1)\,. 
%\end{align}
For example, 
\begin{align}
\Delta(\tilh_{1,1})
&=\Delta^{(1)}(\tilh_{1,1})+\Delta^{(0)}(\tilh_{1,1}) \el \\
\text{where}\quad 
\Delta^{(1)}(\tilh_{1,1})
&= \tilh_{1,1}\otimes 1+1\otimes \tilh_{1,1} \el \\
\Delta^{(0)}(\tilh_{1,1})
&= -2 x_{1,0}^-\otimes x_{1,0}^+ 
+ \sum_{\mu=3,4} \bigl(E_{1\mu}\otimes E_{\mu1}- E_{2\mu}\otimes E_{\mu2}\bigr)\,.
\end{align}

\ms

With this notation, the relation \eqref{hh1} is written as 
\begin{align}
[\Delta(\tilh_{1,1}), \Delta(\tilh_{2,1})]
&=[\Delta^{(1)}(\tilh_{1,1}), \Delta^{(1)}(\tilh_{2,1})] \el \\
&+[\Delta^{(1)}(\tilh_{1,1}), \Delta^{(0)}(\tilh_{2,1})]
+[\Delta^{(0)}(\tilh_{1,1}), \Delta^{(1)}(\tilh_{2,1})] \el \\
&+[\Delta^{(0)}(\tilh_{1,1}), \Delta^{(0)}(\tilh_{2,1})]\,. 
\end{align}
The degree two part is trivial, 
\begin{align}
[\Delta^{(1)}(\tilh_{1,1}), \Delta^{(1)}(\tilh_{2,1})]
=[\tilh_{1,1}, \tilh_{2,1}]\otimes 1+1\otimes [\tilh_{1,1}, \tilh_{2,1}] 
=0\,. 
\end{align}
The degree one terms are computed respectively as 
\begin{align}
[\Delta^{(1)}(\tilh_{1,1}), \Delta^{(0)}(\tilh_{2,1})]
&=[\tilh_{1,1}\otimes 1+1\otimes \tilh_{1,1} \,, \el \\
&\qquad E_{12}\otimes E_{21}- E_{13}\otimes E_{31}+E_{34}\otimes E_{43} 
%\el \\
%&\hspace{40mm}
+ E_{24}\otimes E_{42} -P_0^-\otimes P_0^+ ]
\el \\
&=-2(x_{1,1}^-\otimes x_{1,0}^+-x_{1,0}^-\otimes x_{1,1}^+) \el \\
&\quad +[x_{1,1}^-,x_{2,0}^-]\otimes [x_{1,0}^+,x_{2,0}^+]
-[x_{1,0}^-,x_{2,0}^-]\otimes [x_{1,1}^+,x_{2,0}^+] 
\el \\
&\quad -[x_{2,1}^-,x_{3,0}^-]\otimes [x_{2,0}^+,x_{3,0}^+]
+[x_{2,0}^-,x_{3,0}^-]\otimes [x_{2,1}^+,x_{3,0}^+] 
\el \\
&\quad +\frac{1}{2}\{x_{1,0}^-,x_{2,0}^-\}\otimes [x_{1,0}^+,x_{2,0}^+]
+\frac{1}{2}[x_{1,0}^-,x_{2,0}^-]\otimes \{x_{1,0}^+,x_{2,0}^+\} \,, 
\el \\
[\Delta^{(0)}(\tilh_{1,1}), \Delta^{(1)}(\tilh_{2,1})]
&=[-2 x_{1,0}^-\otimes x_{1,0}^+ 
+ \sum_{\mu=3,4} \bigl(E_{1\mu}\otimes E_{\mu1}- E_{2\mu}\otimes E_{\mu2}\bigr) \,, 
\el \\
&\qquad \tilh_{2,1}\otimes 1+1\otimes \tilh_{2,1} ]
\el \\
&=2(x_{1,1}^-\otimes x_{1,0}^+-x_{1,0}^-\otimes x_{1,1}^+) \el \\
&\quad -[x_{1,1}^-,x_{2,0}^-]\otimes [x_{1,0}^+,x_{2,0}^+]
+[x_{1,0}^-,x_{2,0}^-]\otimes [x_{1,1}^+,x_{2,0}^+] 
\el \\
&\quad +[x_{2,0}^-,x_{3,1}^-]\otimes [x_{2,0}^+,x_{3,0}^+]
-[x_{2,0}^-,x_{3,0}^-]\otimes [x_{2,0}^+,x_{3,1}^+] 
\el \\
&\quad +\frac{1}{2}\{E_{13},E_{34} \}\otimes E_{41} 
-\frac{1}{2}E_{14}\otimes \{E_{31},E_{43} \} 
\el \\
&\quad -\frac{1}{2}\{E_{12},E_{24} \}\otimes E_{41} 
+\frac{1}{2}E_{14}\otimes \{E_{21},E_{42} \} \,. 
\end{align}
Combining them, we have 
\begin{align}
&\quad [\Delta^{(1)}(\tilh_{1,1}), \Delta^{(0)}(\tilh_{2,1})]
+[\Delta^{(0)}(\tilh_{1,1}), \Delta^{(1)}(\tilh_{2,1})] 
\el \\ 
&=-([x_{2,1}^-,x_{3,0}^-]-[x_{2,0}^-,x_{3,1}^-])\otimes [x_{2,0}^+,x_{3,0}^+]
+[x_{2,0}^-,x_{3,0}^-]\otimes([x_{2,1}^+,x_{3,0}^+]- [x_{2,0}^+,x_{3,1}^+] )
\el \\
&\quad +\frac{1}{2}\{x_{1,0}^-,x_{2,0}^-\}\otimes [x_{1,0}^+,x_{2,0}^+]
+\frac{1}{2}[x_{1,0}^-,x_{2,0}^-]\otimes \{x_{1,0}^+,x_{2,0}^+\} \,, 
\el \\
&\quad +\frac{1}{2}\{E_{13},E_{34} \}\otimes E_{41} 
-\frac{1}{2}E_{14}\otimes \{E_{31},E_{43} \} 
\el \\
&\quad -\frac{1}{2}\{E_{12},E_{24} \}\otimes E_{41} 
+\frac{1}{2}E_{14}\otimes \{E_{21},E_{42} \} 
\el \\
&=-\frac{1}{2}\{E_{23},E_{34}\}\otimes E_{42} 
+\frac{1}{2} E_{24}\otimes  \{E_{32},E_{43}\} 
\el \\
%
%-([x_{2,1}^-,x_{3,0}^-]-[x_{2,0}^-,x_{3,1}^-])\otimes [x_{2,0}^+,x_{3,0}^+]
%+[x_{2,0}^-,x_{3,0}^-]\otimes([x_{2,1}^+,x_{3,0}^+]- [x_{2,0}^+,x_{3,1}^+] )
%\el \\
&\quad +\frac{1}{2}\{E_{12},E_{23}\}\otimes E_{31} 
-\frac{1}{2} E_{13}\otimes  \{E_{21},E_{32}\} 
\el \\
%&\quad +\frac{1}{2}\{x_{1,0}^-,x_{2,0}^-\}\otimes [x_{1,0}^+,x_{2,0}^+]
%+\frac{1}{2}[x_{1,0}^-,x_{2,0}^-]\otimes \{x_{1,0}^+,x_{2,0}^+\} \,, 
%\el \\
&\quad +\frac{1}{2}\{E_{13},E_{34} \}\otimes E_{41} 
-\frac{1}{2}E_{14}\otimes \{E_{31},E_{43} \} 
\el \\
&\quad -\frac{1}{2}\{E_{12},E_{24} \}\otimes E_{41} 
+\frac{1}{2}E_{14}\otimes \{E_{21},E_{42} \} 
~\in \Y_L(\g)_0\otimes \Y_L(\g)_0
\,, 
\end{align}
where we have used the relations \eqref{x1x0} in the last equality.  
Here we observe that the above quantity is degree zero. 
Furthermore, it is equal to the {\it minus} sign of 
\begin{align}
[\Delta^{(0)}(\tilh_{1,1}), \Delta^{(0)}(\tilh_{2,1})]
&=\bigl[-2 x_{1,0}^-\otimes x_{1,0}^+ 
+ \sum_{\mu=3,4} \bigl(E_{1\mu}\otimes E_{\mu1}- E_{2\mu}\otimes E_{\mu2}\bigr) 
\,, \el \\
&
E_{12}\otimes E_{21}- E_{13}\otimes E_{31}+E_{34}\otimes E_{43}+ E_{24}\otimes E_{42} 
-P_0^-\otimes P_0^+
\bigr]\,.  \el 
\end{align}
Hence, we prove \eqref{hh1}. 
Due to the automorphism $\rho$, this also implies that 
\begin{align}
[\Delta(\tilh_{3,1}), \Delta(\tilh_{2,1})]=0\,. 
\end{align}

\ms

Another relation \eqref{hh2} is proved similarly. 
The degree decomposition is 
\begin{align}
[\Delta(\tilh_{1,1}), \Delta(\tilh_{3,1})]
&=[\Delta^{(1)}(\tilh_{1,1}), \Delta^{(1)}(\tilh_{3,1})] \el \\
&+[\Delta^{(1)}(\tilh_{1,1}), \Delta^{(0)}(\tilh_{3,1})]
+[\Delta^{(0)}(\tilh_{1,1}), \Delta^{(1)}(\tilh_{3,1})] \el \\
&+[\Delta^{(0)}(\tilh_{1,1}), \Delta^{(0)}(\tilh_{3,1})]\,. 
\end{align}
The degree two term vanishes trivially, 
\begin{align}
[\Delta^{(1)}(\tilh_{1,1}), \Delta^{(1)}(\tilh_{3,1})]
=[\tilh_{1,1}, \tilh_{3,1}]\otimes 1 
+1\otimes [\tilh_{1,1}, \tilh_{3,1}] =0\,. 
\end{align}
The summation of the degree one terms turns out to be 
\begin{align}
&\quad [\Delta^{(1)}(\tilh_{1,1}), \Delta^{(0)}(\tilh_{3,1})]
+[\Delta^{(0)}(\tilh_{1,1}), \Delta^{(1)}(\tilh_{3,1})] \el \\
&=-\{E_{12},E_{23}\}\otimes E_{31} +E_{13}\otimes \{E_{21},E_{32}\} \el \\
&\quad -\{E_{13},E_{34}\}\otimes E_{41} +E_{14}\otimes \{E_{31},E_{43}\} \el \\
&\quad+\{E_{12},E_{24}\}\otimes E_{41} -E_{14}\otimes \{E_{21},E_{42}\} \el \\
&\quad+\{E_{23},E_{34}\}\otimes E_{42} -E_{24}\otimes \{E_{32},E_{43}\} 
~\in \Y_L(\g)_0\otimes \Y_L(\g)_0
\,. 
\end{align}
This is again equal to the opposite sign of 
\begin{align}
[\Delta^{(0)}(\tilh_{1,1}), \Delta^{(0)}(\tilh_{3,1})]
&= [-2 x_{1,0}^-\otimes x_{1,0}^+ 
+ \sum_{\mu=3,4} \bigl(E_{1\mu}\otimes E_{\mu1}- E_{2\mu}\otimes E_{\mu2}\bigr)\,, 
\el \\
&\quad 2 x_{3,0}^-\otimes x_{3,0}^+ 
+\sum_{l=1,2} \bigl( E_{l3}\otimes E_{3l}- E_{l4}\otimes E_{4l} \bigr)]\,. 
\end{align}
Thus, we prove \eqref{hh2}. 
Therefore, we complete the proof of the compatibility of the coproducts with the degree two relations \eqref{deg-2}.

%\Delta(\tilh_{1,1})&= \tilh_{1,1}\otimes 1+1\otimes \tilh_{1,1} 
%-2 x_{1,0}^-\otimes x_{1,0}^+ 
%+ \sum_{\mu=3,4} \bigl(E_{1\mu}\otimes E_{\mu1}- E_{2\mu}\otimes E_{\mu2}\bigr) \el\\ 
%
%\Delta(\tilh_{2,1})&= \tilh_{2,1}\otimes 1+1\otimes \tilh_{2,1} 
%+E_{12}\otimes E_{21}- E_{13}\otimes E_{31}+E_{34}\otimes E_{43}+ E_{24}\otimes E_{42} 
%-P_0^-\otimes P_0^+ \el\\ 
%
%\Delta(\tilh_{3,1})&= \tilh_{3,1}\otimes 1+1\otimes \tilh_{3,1} 
%+2 x_{3,0}^-\otimes x_{3,0}^+ 
%+\sum_{l=1,2} \bigl( E_{l3}\otimes E_{3l}- E_{l4}\otimes E_{4l} \bigr)\el\\ 

\ms 

Finally, we show the compatibility with the degree three relations \eqref{deg-31} and \eqref{deg-32}, the most challenging part of this proof.

\ms 

\noindent $\bullet$ 
Relation \eqref{deg-31}. 
%\noindent $\diamond$ Proof of $\bigl[\Delta(\tilh_{1,1}),[\Delta(x^+_{1,1}),\Delta(x^-_{1,1})]\bigr]=0$\,. 
Set the degree two element by $h_{1,2}=[x^+_{1,1}, x^-_{1,1}] $\,. 
The degree decomposition is 
\begin{align}
\Delta(h_{1,2})=[\Delta(x^+_{1,1}),\Delta(x^-_{1,1})] 
=\Delta^{(2)}(h_{1,2})+\Delta^{(1)}(h_{1,2})+\Delta^{(0)}(h_{1,2})\,, 
\end{align}
and each term is explicitly calculated by 
\begin{align}
\Delta^{(2)}(h_{1,2})&=[x^+_{1,1},x^-_{1,1}]\otimes1+1\otimes[x^+_{1,1},x^-_{1,1}] \el \\
\Delta^{(1)}(h_{1,2})&=h_{1,1}\otimes h_{1,0}+h_{1,0}\otimes h_{1,1} \el \\
&\quad -2(x_{1,1}^-\otimes x_{1,0}^++x_{1,0}^-\otimes x_{1,1}^+ )\el \\ 
&\quad +[h_{1,1},x_{2,0}^-]\otimes x_{2,0}^+ - x_{2,0}^-\otimes [h_{1,1},x_{2,0}^+ ] \el \\
&\quad +[x_{1,1}^-,x_{2,0}^-]\otimes E_{31} + E_{13} \otimes [x_{1,1}^+,x_{2,0}^+ ] \el \\
&\quad -[[h_{1,1},x_{2,0}^-],x_{3,0}^-]\otimes E_{42}
-E_{24}\otimes[[h_{1,1},x_{2,0}^+],x_{3,0}^+] \el\\
&\quad+ [x_{1,1}^-[x_{2,0}^-,x_{3,0}^-]]\otimes E_{41}
+E_{14}\otimes [x_{1,1}^+,[x_{2,0}^+,x_{3,0}^+]] \el \\ 
\Delta^{(0)}(h_{1,2})
&= -\{h_{1,0},x_{1,0}^-\}\otimes x_{1,0}^+ -x_{1,0}^-\otimes\{h_{1,0},x_{1,0}^+\} \el \\
&\quad +\frac{1}{2}\bigl(
(\{h_{1,0},E_{1\mu}\}+\{x_{1,0}^-,E_{2\mu}\})\otimes E_{\mu1}
+E_{1\mu}\otimes (\{h_{1,0},E_{\mu1}\}+\{x_{1,0}^+,E_{\mu2}\}) \el\\
%&\qquad 
%+\{x_{1,0}^-,E_{2\mu}\}\otimes E_{\mu1}+E_{1\mu}\otimes \{x_{1,0}^+,E_{\mu2}\} \el\\
&\qquad 
-(\{E_{13},E_{24}\}+\{E_{23},E_{14}\})\otimes P_0^+
-P_0^+ \otimes (\{E_{31},E_{42}\} +\{E_{32},E_{41}\})\,. 
\label{h12-deg}
%&\qquad 
%-\{E_{23},E_{14}\}\otimes P_0^+-P_0^+ \otimes \{E_{32},E_{41}\} \bigr)
%\el \\
%\Delta^{(2)}(h_{2,2})&=[x^+_{2,1},x^-_{2,1}]\otimes1+1\otimes[x^+_{2,1},x^-_{2,1}]\el \\
%\Delta^{(1)}(h_{2,2})&=
%h_{2,1}\otimes h_{2,0} + h_{2,0}\otimes h_{2,1} \el \\
%&\quad + x_{1,1}^-\otimes x_{1,0}^++x_{1,0}^-\otimes x_{1,1}^+ \el \\
%&\quad - x_{3,1}^-\otimes x_{3,0}^+ -x_{3,0}^-\otimes x_{3,1}^+ \el \\
%&\quad -[x_{1,0}^-,x_{2,1}^-]\otimes E_{31} - E_{13} \otimes [x_{1,0}^+,x_{2,1}^+] \el \\
%&\quad +[x_{2,1}^-,x_{3,0}^-]\otimes E_{42} - E_{24} \otimes [x_{2,1}^+,x_{3,0}^+] \el \\
%&\quad +\frac{1}{2} \bigl(
%\{h_{2,0},x_{1,0}^-\} \otimes x_{1,0}^+ +x_{1,0}^-\otimes \{h_{2,0},x_{1,0}^+ \} \bigr) \el \\ 
%&\quad - \frac{1}{2} \bigl(
%\{h_{2,0},x_{3,0}^-\} \otimes x_{3,0}^+ +x_{3,0}^-\otimes \{h_{2,0},x_{3,0}^+ \} \bigr) 
%\el \\ 
%\Delta^{(0)}(h_{2,2})&= 
%-\frac{1}{2}\bigl(
%(\{h_{2,0}, E_{13}\}+\{E_{12},E_{23}\})\otimes E_{31} 
%+E_{13} \otimes (\{h_{2,0}, E_{31}\}+\{E_{21},E_{32}\}) \el \\
%&\qquad +(\{E_{12}, E_{24}\}-\{E_{34},E_{13}\}) \otimes E_{41} 
%+E_{14}\otimes (\{E_{21}, E_{42}\}-\{E_{43},E_{31}\}) \el \\
%&\qquad -(\{h_{2,0}, E_{24}\}+\{E_{34},E_{23}\}) \otimes E_{42} 
%-E_{24}\otimes (\{h_{2,0}, E_{42}\}+\{E_{43},E_{32}\}) \el \\
%&\qquad +\{h_{2,0},P_0^-\}\otimes P_0^+ + P_0^-\otimes \{h_{2,0},P_0^+\} \bigr)\,. 
%\\~
%\el \\ 
%\Delta^{(2)}(\tilh_{3,2})&=[x^+_{3,1},x^-_{3,1}]\otimes1+1\otimes[x^+_{3,1},x^-_{3,1}] \el\\
%\Delta^{(1)}(\tilh_{3,2})&=\el \\
%\Delta^{(0)}(\tilh_{3,2})&= 
\end{align}
Then, the decomposition of the left hand side of \eqref{deg-31} with $j=1$ becomes 
\begin{align}
\bigl[\Delta(\tilh_{1,1}),[\Delta(x^+_{1,1}),\Delta(x^-_{1,1})]\bigr]
%&=\bigl[\Delta^{(1)}(\tilh_{1,1})+\Delta^{(0)}(\tilh_{1,1}),
%\Delta^{(2)}(\tilh_{1,2})+\Delta^{(1)}(\tilh_{1,2})+\Delta^{(0)}(\tilh_{1,2})\bigr] \el \\
&= \bigl[\Delta^{(1)}(\tilh_{1,1}),\Delta^{(2)}(\tilh_{1,2})\bigr] \el \\
& +\bigl[\Delta^{(1)}(\tilh_{1,1}),\Delta^{(1)}(\tilh_{1,2})\bigr] 
+\bigl[\Delta^{(0)}(\tilh_{1,1}),\Delta^{(2)}(\tilh_{1,2})\bigr] \el \\
& +\bigl[\Delta^{(1)}(\tilh_{1,1}),\Delta^{(0)}(\tilh_{1,2})\bigr] 
+\bigl[\Delta^{(0)}(\tilh_{1,1}),\Delta^{(1)}(\tilh_{1,2})\bigr] \el \\
& +\bigl[\Delta^{(0)}(\tilh_{1,1}),\Delta^{(0)}(\tilh_{1,2})\bigr] \,. 
\end{align}

\ms 

It is obvious that the degree three term vanishes by using the relation \eqref{deg-31}, 
\begin{align}
\bigl[\Delta^{(1)}(\tilh_{1,1}),\Delta^{(2)}(h_{1,2})\bigr]
=\bigl[\tilh_{1,1},[ x^+_{1,1}, x^-_{1,1}]\bigr]\otimes 1 
+1\otimes \bigl[\tilh_{1,1},[ x^+_{1,1}, x^-_{1,1}]\bigr] =0\,. 
\end{align}
The degree two terms are somehow complicated, but after some calculations,  
\begin{align}
&\quad \bigl[\Delta^{(1)}(\tilh_{1,1}),\Delta^{(1)}(h_{1,2})\bigr] 
+\bigl[\Delta^{(0)}(\tilh_{1,1}),\Delta^{(2)}(h_{1,2})\bigr] \nln
&=\frac{1}{2} \Bigl(
-\{h_{1,0},x_{2,0}^-\}\otimes x_{2,1}^+ +x_{2,1}^-\otimes \{h_{1,0},x_{2,0}^+\} \nln 
&\quad 
-\{x_{1,0}^-,x_{2,0}^-\}\otimes [x_{1,0}^+,x_{2,1}^+] 
-[x_{1,0}^-,x_{2,1}^-]\otimes \{x_{1,0}^+,x_{2,0}^+\} \nln 
&\quad 
-\{h_{1,0},[x_{2,0}^-,x_{3,0}^-]\}\otimes [x_{2,1}^+,x_{3,0}^+] 
+[x_{2,1}^-,x_{3,0}^-]\otimes \{h_{1,0},[x_{2,0}^+,x_{3,0}^+]\} \nln 
&\quad 
-\{x_{1,0}^-,[x_{2,0}^-,x_{3,0}^-]\}\otimes [[x_{1,0}^+,x_{2,1}^+],x_{3,0}^+] 
-[[x_{1,0}^-,x_{2,1}^-],x_{3,0}^-]\otimes \{x_{1,0}^+,[x_{2,0}^+,x_{3,0}^+]\} \Bigr) \nln 
&\quad -2 (\{h_{1,0},x_{1,1}^-\}+\{h_{1,1},x_{1,0}^-\})\otimes x_{1,0}^+ 
+2x_{1,0}^-\otimes (\{h_{1,0},x_{1,1}^+\}+\{h_{1,1},x_{1,0}^+\}) \nln
&\quad 
+\frac{1}{2}(\{x_{1,0}^-,x_{2,1}^-\}-[\tilh_{1,1},\{h_{1,0},[x_{1,0}^-,x_{2,0}^-]\}]
+\{h_{1,1},[x_{1,0}^-,x_{2,0}^-]\})\otimes E_{31} \nln
&\quad 
+\frac{1}{2}E_{13}\otimes (\{x_{1,0}^+,x_{2,1}^+\}-[\tilh_{1,1},\{h_{1,0},[x_{1,0}^+,x_{2,0}^+]\}]
+\{h_{1,1},[x_{1,0}^+,x_{2,0}^+]\})  \nln
&\quad 
+\frac{1}{2} \bigl(\{x_{1,0}^-,[x_{2,1}^-,x_{3,0}^-]\}
+\{h_{1,0},[[x_{1,1}^-,x_{2,0}^-],x_{3,0}^-]\} +\{h_{1,1},[[x_{1,0}^-,x_{2,0}^-],x_{3,0}^-]\} 
\bigr)\otimes E_{41} \nln
&\quad 
+\frac{1}{2} E_{14}\otimes \bigl(\{x_{1,0}^+,[x_{2,1}^+,x_{3,0}^+]\}
-\{h_{1,0},[[x_{1,1}^+,x_{2,0}^+],x_{3,0}^+]\} -\{h_{1,1},[[x_{1,0}^+,x_{2,0}^+],x_{3,0}^+]\} 
\bigr)\nln
&\quad 
+\frac{1}{4}\bigl(
\{h_{1,0},[[x_{1,0}^-,x_{2,0}^-],x_{3,0}^-]\} \otimes E_{41} 
+E_{14}\otimes \{h_{1,0},[[x_{1,0}^-,x_{2,0}^-],x_{3,0}^-]\} \bigr) \nln
&\quad +\frac{1}{2} \bigl(-\{h_{1,1},x_{2,0}^-\}\otimes x_{2,0}^+
+ x_{2,0}^-\otimes \{h_{1,1},x_{2,0}^+\} \bigr) \nln 
&\quad +\frac{1}{2} \bigl(\{h_{1,1},E_{24} \}\otimes E_{42} 
-E_{24} \otimes \{h_{1,1},E_{42}\} \bigr)\,,  
\end{align}
we observe that these elements are essentially in degree one, 
\begin{align}
\bigl[\Delta^{(1)}(\tilh_{1,1}),\Delta^{(1)}(h_{1,2})\bigr] 
+\bigl[\Delta^{(0)}(\tilh_{1,1}),\Delta^{(2)}(h_{1,2})\bigr]
\nln
\in %\sum_{i=0,1}\Y_L(\g)_{1-i} \otimes \Y_L(\g)_i \,. 
\Y_L(\g)_{1} \otimes \Y_L(\g)_0+\Y_L(\g)_{0} \otimes \Y_L(\g)_1
\end{align}
Furthermore, adding the degree one terms to the above, we have  
\begin{align}
&\quad [\Delta^{(1)}(\tilh_{1,1}),\Delta^{(1)}(h_{1,2})\bigr] 
+\bigl[\Delta^{(0)}(\tilh_{1,1}),\Delta^{(2)}(h_{1,2})\bigr] \nln
&\quad + \bigl[\Delta^{(1)}(\tilh_{1,1}),\Delta^{(0)}(h_{1,2})\bigr] 
+\bigl[\Delta^{(0)}(\tilh_{1,1}),\Delta^{(1)}(h_{1,2})\bigr] 
\nln
&= \frac{1}{2}\sum_{\mu=3,4} \bigl(-\{\{h_{1,0},E_{2\mu}\},x_{1,0}^-\}\otimes E_{\mu1} 
+E_{1\mu}\otimes \{\{h_{1,0},E_{\mu2}\},x_{1,0}^+\} \bigr)\nln
&+\frac{1}{2}\sum_{\mu=3,4} \bigl(
\{h_{1,0},E_{1\mu}\}\otimes \{x_{1,0}^+,E_{\mu2} \}
-\{x_{1,0}^+,E_{2\mu} \}\otimes \{h_{1,0},E_{\mu1}\}\bigr) \nln
&+\frac{1}{4}
\bigl(\{\{E_{12},E_{23}\},E_{24}\} - \{\{E_{12},E_{24}\},E_{23}\}
+\{\{h_{1,0},E_{13}\},E_{24}\}+\{\{h_{1,0},E_{23}\},E_{14}\} \bigr)\otimes P_0^+ \nln
&-\frac{1}{4} P_0^- \otimes 
\bigl(\{\{E_{21},E_{32}\},E_{42}\} - \{\{E_{21},E_{42}\},E_{32}\}
+\{\{h_{1,0},E_{31}\},E_{42}\}+\{\{h_{1,0},E_{32}\},E_{41}\} \bigr)\nln
&+\frac{1}{4}\bigl(
(\{E_{13},E_{24}\}+\{E_{23},E_{14}\})\otimes \{h_{1,0},P_0^+\}
-\{h_{1,0},P_0^-\}\otimes (\{E_{31},E_{42}\}+\{E_{32},E_{41}\}) \bigr) \nln 
&+\frac{1}{2}\bigl(
\{x_{1,0}^-,P_0^-\}\otimes \{E_{31},E_{41} \}
-\{E_{13},E_{14} \} \otimes\{x_{1,0}^+,P_0^+\} \bigr) \,. 
\label{hh11-zero}
\end{align}
We see that there is no degree one element and 
remaining terms are all in degree zero, 
\begin{align}
&\Bigl(\bigl[\Delta^{(1)}(\tilh_{1,1}),\Delta^{(1)}(h_{1,2})\bigr] 
+\bigl[\Delta^{(0)}(\tilh_{1,1}),\Delta^{(2)}(h_{1,2})\bigr] \nln
&+\bigl[\Delta^{(1)}(\tilh_{1,1}),\Delta^{(0)}(h_{1,2})\bigr] 
+\bigl[\Delta^{(0)}(\tilh_{1,1}),\Delta^{(1)}(h_{1,2})\bigr] \Bigr)
\in 
\Y_L(\g)_0 \otimes  \Y_L(\g)_0 \,. 
%\subset \Y(\g)_0\,.  
\end{align} 
Finally, we notice that the resulting terms in \eqref{hh11-zero}
exactly coincide with the minus sign of 
$\bigl[\Delta^{(0)}(\tilh_{1,1}),\Delta^{(0)}(h_{1,2})\bigr] $\,.
%the degree zero term. 
Thus, summing up all, we have 
\begin{align}
&\quad \bigl[\Delta^{(1)}(\tilh_{1,1}),\Delta^{(1)}(h_{1,2})\bigr] 
+\bigl[\Delta^{(0)}(\tilh_{1,1}),\Delta^{(2)}(h_{1,2})\bigr] \nln
&+\bigl[\Delta^{(1)}(\tilh_{1,1}),\Delta^{(0)}(h_{1,2})\bigr] 
+\bigl[\Delta^{(0)}(\tilh_{1,1}),\Delta^{(1)}(h_{1,2})\bigr] \el \\
& +\bigl[\Delta^{(0)}(\tilh_{1,1}),\Delta^{(0)}(h_{1,2})\bigr] \nln
&=0\,. 
\end{align}
This proves the desired result 
\begin{align}
\bigl[\Delta(\tilh_{1,1}),\Delta(h_{1,2})\bigr] 
=\bigl[\Delta(\tilh_{1,1}),[\Delta(x_{1,1}^+),\Delta(x_{1,1}^-)]\bigr] 
=0\,. 
\end{align}

\ms

\noindent $\bullet$ 
Relation \eqref{deg-32}. 
We shall repeat similar computations. 
%\noindent $\diamond$ Proof of 
%$\bigl[\Delta(\tilh_{1,1}),[\Delta(x^+_{2,1}),\Delta(x^-_{2,1})]\bigr]=0$\,. 
%
Set $h_{2,2}=[x^+_{2,1}, x^-_{2,1}]$\,.
The degree decomposition is 
\begin{align}
\bigl[\Delta(\tilh_{1,1}),[\Delta(x^+_{2,1}),\Delta(x^-_{2,1})]\bigr]
%&=[\Delta(\tilh_{1,1}),\Delta(h_{2,2})] \nln
&=\bigl[\Delta^{(2)}(\tilh_{1,1}),\Delta^{(1)}(h_{2,2})\bigr] \nln
&
+\bigl[\Delta^{(1)}(\tilh_{1,1}),\Delta^{(1)}(h_{2,2})\bigr] 
+\bigl[\Delta^{(0)}(\tilh_{1,1}),\Delta^{(2)}(h_{2,2})\bigr] \nln
&+\bigl[\Delta^{(1)}(\tilh_{1,1}),\Delta^{(0)}(h_{2,2})\bigr] 
+\bigl[\Delta^{(0)}(\tilh_{1,1}),\Delta^{(1)}(h_{2,2})\bigr] \el \\
& +\bigl[\Delta^{(0)}(\tilh_{1,1}),\Delta^{(0)}(h_{2,2})\bigr]\,. 
\end{align}
%where 
%$\Delta^{(r)}(h_{2,2})\equiv [\Delta(x^+_{2,1}),\Delta(x^-_{2,1})]$ 
%with $r=2,1,0$ are given by 
%The degree decompositions 
Each degree component of 
\begin{align}
\Delta(h_{2,2})=[\Delta(x^+_{2,1}),\Delta(x^-_{2,1})]\bigr]
=\Delta^{(2)}(h_{2,2})+\Delta^{(1)}(h_{2,2})+\Delta^{(0)}(h_{2,2})\,, 
\end{align}
is calculated by 
\begin{align}
\Delta^{(2)}(h_{2,2})&=[x^+_{2,1},x^-_{2,1}]\otimes1+1\otimes[x^+_{2,1},x^-_{2,1}]\el \\
\Delta^{(1)}(h_{2,2})&=
h_{2,1}\otimes h_{2,0} + h_{2,0}\otimes h_{2,1} \el \\
&\quad + x_{1,1}^-\otimes x_{1,0}^++x_{1,0}^-\otimes x_{1,1}^+ \el \\
&\quad - x_{3,1}^-\otimes x_{3,0}^+ -x_{3,0}^-\otimes x_{3,1}^+ \el \\
&\quad -[x_{1,0}^-,x_{2,1}^-]\otimes E_{31} - E_{13} \otimes [x_{1,0}^+,x_{2,1}^+] \el \\
&\quad +[x_{2,1}^-,x_{3,0}^-]\otimes E_{42} - E_{24} \otimes [x_{2,1}^+,x_{3,0}^+] \el \\
&\quad 
-[E_{14},x_{2,1}^+]\otimes P_0^+- P_0^-\otimes [E_{41},x_{2,1}^-]\nln
&\quad +\frac{1}{2} \bigl(
\{h_{2,0},x_{1,0}^-\} \otimes x_{1,0}^+ +x_{1,0}^-\otimes \{h_{2,0},x_{1,0}^+ \} \bigr) \el \\ 
&\quad - \frac{1}{2} \bigl(
\{h_{2,0},x_{3,0}^-\} \otimes x_{3,0}^+ +x_{3,0}^-\otimes \{h_{2,0},x_{3,0}^+ \} \bigr) 
\el \\ 
\Delta^{(0)}(h_{2,2})&= 
-\frac{1}{2}\Bigl(
(\{h_{2,0}, E_{13}\}+\{E_{12},E_{23}\})\otimes E_{31} 
+E_{13} \otimes (\{h_{2,0}, E_{31}\}+\{E_{21},E_{32}\}) \el \\
&\qquad +(\{E_{12}, E_{24}\}-\{E_{34},E_{13}\}) \otimes E_{41} 
+E_{14}\otimes (\{E_{21}, E_{42}\}-\{E_{43},E_{31}\}) \el \\
&\qquad -(\{h_{2,0}, E_{24}\}+\{E_{34},E_{23}\}) \otimes E_{42} 
-E_{24}\otimes (\{h_{2,0}, E_{42}\}+\{E_{43},E_{32}\}) \el \\
&\qquad +\{h_{2,0},P_0^-\}\otimes P_0^+ + P_0^-\otimes \{h_{2,0},P_0^+\} \Bigr)\,. 
\end{align}
Firstly, 
it is obvious that the degree three term vanishes by the relation \eqref{deg-32}, 
\begin{align}
\bigl[\Delta^{(1)}(\tilh_{1,1}),\Delta^{(2)}(h_{2,2})\bigr]
=\bigl[\tilh_{1,1},[ x^+_{2,1}, x^-_{2,1}]\bigr]\otimes 1 
+1\otimes \bigl[\tilh_{1,1},[ x^+_{2,1}, x^-_{2,1}]\bigr] =0\,. 
\end{align}
Though the degree two terms look quite complicated, 
\begin{align}
&\quad \bigl[\Delta^{(1)}(\tilh_{1,1}),\Delta^{(1)}(h_{2,2})\bigr] 
+\bigl[\Delta^{(0)}(\tilh_{1,1}),\Delta^{(2)}(h_{2,2})\bigr]  \nln
&=\bigl(\{h_{2,1},x_{1,0}^-\} + \{h_{2,0},x_{1,1}^-\}\bigr)\otimes x_{1,0}^+
-x_{1,0}^-\otimes \bigl(\{h_{2,1},x_{1,0}^+\} + \{h_{2,0},x_{1,1}^+\}\bigr) \nln
&\quad 
+\frac{1}{2}\Bigl(\bigl(-\{x_{1,1}^-,x_{2,0}^-\} + \{x_{1,0}^-,x_{2,1}^-\}
\bigr)\otimes E_{31} 
+E_{13}\otimes \bigl(-\{x_{1,1}^+, x_{2,0}^+\} + \{x_{1,0}^+, x_{2,1}^+\} \bigr) \Bigr)
\nln 
&\quad 
+\frac{1}{2}\bigl(
-\{h_{2,1},[x_{1,0}^-,x_{2,0}^-]\} -\{h_{2,0}, [x_{1,1}^-,x_{2,0}^-]\}\bigr)\otimes E_{31} \nln
&\quad 
+\frac{1}{2}\tss E_{13}\otimes\bigl(
\{h_{2,1},[x_{1,0}^+,x_{2,0}^+]\} + \{h_{2,0}, [x_{1,1}^+,x_{2,0}^+]\}\bigr) \nln
&\quad 
+\{x_{1,0}^-,x_{2,0}^-\} \otimes [x_{1,1}^+,x_{2,0}^+] 
+[x_{1,1}^-,x_{2,0}^-] \otimes \{x_{1,0}^+,x_{2,0}^+\} \nln
&\quad   
+\frac{1}{2}\Bigl(-\bigl(\{x_{2,1}^-,x_{3,0}^-\} + \{x_{2,0}^-,x_{3,1}^-\}
\bigr)\otimes E_{42} 
+E_{24}\otimes \bigl( \{x_{2,1}^+, x_{3,0}^+\} + \{x_{2,0}^+, x_{3,1}^+\} \bigr) \Bigr)\nln 
&\quad 
+\frac{1}{2}\bigl(
-\{h_{2,1},[x_{2,0}^-,x_{3,0}^-]\} -\{h_{2,0}, [x_{2,0}^-,x_{3,1}^-]\}\bigr)\otimes E_{42} \nln
&\quad 
+\frac{1}{2}\tss E_{24}\otimes\bigl(
-\{h_{2,1},[x_{2,0}^+,x_{3,0}^+]\} - \{h_{2,0}, [x_{2,0}^+,x_{3,1}^+]\} \bigr) 
\nln
&\quad 
+\frac{1}{2}\bigl(\{E_{13},[x_{2,1}^-,x_{3,0}^- ]\} - \{[x_{1,0}^-,x_{2,1}^- ],E_{24}\}
\bigr)\otimes P_0^+ \nln 
&\quad 
+\frac{1}{4} \bigl(\{\{E_{13},E_{23}\},E_{34} \}+\{\{E_{13},E_{34}\},E_{23} \}\bigr) 
\otimes P_0^+ \nln 
&\quad 
+\frac{1}{4} \bigl( \{\{E_{24},E_{12}\},E_{23} \}+\{\{E_{24},E_{23}\},E_{12} \}
\bigr)\otimes P_0^+ \nln 
&\quad 
+\frac{1}{2}\tss P_0^-\otimes 
\bigl(\{E_{31},[x_{2,1}^+,x_{3,0}^+ ]\} + \{[x_{1,0}^+,x_{2,1}^+ ],E_{42}\} \bigr) \nln 
&\quad 
-\frac{1}{4} \tss P_0^-\otimes 
\bigl(\{\{E_{31},E_{32}\},E_{43} \} - \{\{E_{31},E_{43}\},E_{32} \}\bigr) \nln 
&\quad 
-\frac{1}{4}\tss P_0^-\otimes 
\bigl( \{\{E_{42},E_{21}\},E_{32} \} - \{\{E_{42},E_{32}\},E_{21} \} \bigr)\nln 
&\quad 
+\frac{1}{2}\Bigl(\bigl(\{x_{1,1}^-,E_{24}\}+\{E_{13},x_{3,1}^-\}\bigr)\otimes E_{41} 
+E_{14}\otimes \bigl( \{x_{1,1}^+,E_{42}\}+\{E_{31},x_{3,1}^+\}\bigr) \Bigr)  \nln
&\quad  
+\frac{1}{2}\bigl(
-\{x_{1,0}^-,[x_{2,1}^-,x_{3,0}^-]\} + \{[x_{1,0}^-,x_{2,1}^-],x_{3,0}^-\}
\bigr)\otimes E_{41} \nln
&\quad  
+\frac{1}{2}\tss E_{14}\otimes \bigl(
-\{x_{1,0}^+,[x_{2,1}^+,x_{3,0}^+]\} + \{[x_{1,0}^+,x_{2,1}^+],x_{3,0}^+\} \bigr)\nln
&\quad  
+\frac{1}{4}\bigl(
\{h_{2,0} ,\{E_{13},E_{34}\} \} - \{h_{2,0} ,\{E_{12},E_{24}\} \} 
\bigr)\otimes E_{41} \nln
&\quad  
+\frac{1}{4}\tss E_{14}\otimes \bigl(
-\{h_{2,0} ,\{E_{31},E_{43}\} \} + \{h_{2,0} ,\{E_{21},E_{42}\} \} \bigr)\,,  
\end{align}
we observe that the above elements are in degree one, 
\begin{align}
\bigl[\Delta^{(1)}(\tilh_{1,1}),\Delta^{(1)}(h_{2,2})\bigr] 
+\bigl[\Delta^{(0)}(\tilh_{1,1}),\Delta^{(2)}(h_{2,2})\bigr]
\nln
\in \Y_L(\g)_1\otimes \Y_L(\g)_0+\Y_L(\g)_0\otimes \Y_L(\g)_1 \,. 
\end{align}
Furthermore, adding the degree one terms to the above, we have  
\begin{align}
&\quad [\Delta^{(1)}(\tilh_{1,1}),\Delta^{(1)}(h_{2,2})\bigr] 
+\bigl[\Delta^{(0)}(\tilh_{1,1}),\Delta^{(2)}(h_{2,2})\bigr] \nln
&\quad + \bigl[\Delta^{(1)}(\tilh_{1,1}),\Delta^{(0)}(h_{2,2})\bigr] 
+\bigl[\Delta^{(0)}(\tilh_{1,1}),\Delta^{(1)}(h_{2,2})\bigr] \nln
&
=\frac{1}{4}\Bigl(
\{h_{2,0},\{x_{1,0}^-,x_{2,0}^-\}\} \otimes E_{31} 
+E_{13}\otimes \{h_{2,0},\{x_{1,0}^+,x_{2,0}^+\}\} \nln
&\qquad 
+\{h_{2,0},E_{13}\} \otimes \{x_{1,0}^+,x_{2,0}^+\}
+\{x_{1,0}^-,x_{2,0}^-\} \otimes \{h_{2,0},E_{31}\} \Bigr) \nln 
&
+ \frac{1}{4}\Bigl(
-\{h_{2,0},\{x_{2,0}^-,x_{3,0}^-\}\} \otimes E_{42}  
+E_{24} \otimes \{h_{2,0},\{x_{2,0}^+,x_{3,0}^+\}\} \nln
&\qquad 
+\{h_{2,0},E_{24}\} \otimes \{x_{2,0}^+,x_{3,0}^+\}
-\{x_{2,0}^-,x_{3,0}^-\} \otimes \{h_{2,0},E_{42}\} \Bigr) \nln 
&
+\frac{1}{2}\bigl(-\{\{x_{1,0}^-,x_{2,0}^-\},x_{3,0}^-\} \otimes E_{41} 
+E_{14}\otimes \{\{x_{1,0}^+,x_{2,0}^+\},x_{3,0}^+\} 
\bigr) \nln
&
+\frac{1}{4}\Bigl(\bigl(\{h_{2,0},\{E_{13},E_{34}\}\} - \{h_{2,0},\{E_{12},E_{24}\}\}
\bigr)\otimes E_{41} \nln
&\quad 
+E_{14}\otimes \bigl(-\{h_{2,0},\{E_{31},E_{43}\}\} + \{h_{2,0},\{E_{21},E_{42}\}\} \bigr)
\Bigr) \nln 
&
+\frac{1}{2} \bigl(
\{E_{34},E_{13}\}\otimes \{E_{21},E_{42}\} - \{E_{12},E_{24}\}\otimes \{E_{43},E_{31}\} 
\bigr)  \nln
&
+\frac{1}{4} \Bigl(\bigl(\{E_{13},E_{34}\} - \{E_{12},E_{24}\}\bigr)\otimes \{E_{41},h_{2,0} \}
\nln
&\quad 
+\{E_{14},h_{2,0} \}\otimes \bigl(-\{E_{31},E_{43}\} + \{E_{21},E_{42}\} \bigr) \Bigr)
\nln 
&
+\frac{1}{4}\bigl(\{\{E_{13},E_{23}\},E_{34} \} + \{\{E_{13},E_{34}\},E_{23} \} \bigr)
\otimes P_0^+ \nln 
&
+\frac{1}{4} \bigl(\{\{E_{24},E_{12}\},E_{23} \} + \{\{E_{24},E_{23}\},E_{12} \} \bigr)
\otimes P_0^+ \nln
&
-\frac{1}{4}\tss P_0^-\otimes 
\bigl(\{\{E_{31},E_{32}\},E_{43} \} - \{\{E_{31},E_{43}\},E_{32} \} \bigr) \nln 
&
-\frac{1}{4}\tss P_0^-\otimes 
\bigl(\{\{E_{42},E_{21}\},E_{32} \} + \{\{E_{42},E_{32}\},E_{21} \} \bigr) \nln
&
+\frac{1}{2}\bigl(
-\{E_{23},E_{14}\} \otimes \{E_{31},E_{42}\} 
+ \{E_{13},E_{24}\} \otimes \{E_{32},E_{41}\} 
\bigr) \nln
&
+\frac{1}{2}\bigl(
-\{E_{13},E_{14}\} \otimes \{x_{1,0}^+,P_0^+ \} 
+ \{x_{1,0}^-,P_0^- \} \otimes \{E_{31},E_{41}\} 
\bigr)\,. 
\label{hh12-zero}
\end{align}
We see that there is no degree one element and 
remaining terms are all degree zero, 
\begin{align}
&\Bigl(\bigl[\Delta^{(1)}(\tilh_{1,1}),\Delta^{(1)}(h_{2,2})\bigr] 
+\bigl[\Delta^{(0)}(\tilh_{1,1}),\Delta^{(2)}(h_{2,2})\bigr] \nln
&\quad +\bigl[\Delta^{(1)}(\tilh_{1,1}),\Delta^{(0)}(h_{2,2})\bigr] 
+\bigl[\Delta^{(0)}(\tilh_{1,1}),\Delta^{(1)}(h_{2,2})\bigr] \Bigr)
\in \Y_L(\g)_0\otimes \Y_L(\g)_0\,.  
\end{align} 
Finally, we observe that the terms in \eqref{hh12-zero}
exactly coincide with the minus sign of the degree zero term, 
$\bigl[\Delta^{(0)}(\tilh_{1,1}),\Delta^{(0)}(h_{2,2})\bigr]$\,. 
Thus, summing up all, we have 
\begin{align}
&\quad \bigl[\Delta^{(1)}(\tilh_{1,1}),\Delta^{(1)}(h_{2,2})\bigr] 
+\bigl[\Delta^{(0)}(\tilh_{1,1}),\Delta^{(2)}(h_{2,2})\bigr] \nln
&+\bigl[\Delta^{(1)}(\tilh_{1,1}),\Delta^{(0)}(h_{2,2})\bigr] 
+\bigl[\Delta^{(0)}(\tilh_{1,1}),\Delta^{(1)}(h_{2,2})\bigr] \el \\
& +\bigl[\Delta^{(0)}(\tilh_{1,1}),\Delta^{(0)}(h_{2,2})\bigr] \nln
&=0\,. 
\end{align}
This proves the desired result 
\begin{align}
\bigl[\Delta(\tilh_{1,1}),\Delta(h_{2,2})\bigr] 
=\bigl[\Delta(\tilh_{1,1}),[\Delta(x_{2,1}^+),\Delta(x_{2,1}^-)]\bigr] 
=0\,. 
\end{align}
Therefore, we have proved that the coproducts defined in \eqref{coprod}
are indeed compatible with the defining relations \eqref{yan-1}--\eqref{yan-serre} of $\Y_L(\g)$\,. 

\ms 

This completes the proof of Proposition \ref{prop:copro}.  
\end{proof}

\section{Drinfeld realization of the Yangian $\Y_D(\g)$}
\label{sec:DII}
\setcounter{equation}{0}

We define the  Drinfeld realization of the Yangian $\Y_D(\g)$ in Definition \ref{def:D2}, 
and we show that it is isomorphic to $\Y_L(\g)$\,, Theorem \ref{thm:D2}. 

\begin{defin}
\label{def:D2}
The Yangian $\Y_D(\g)$ is generated by   the generators  
$h_{i,r}, x^\pm_{i,r}$ with $i=1,2,3$ 
and the central elements $P^\pm_r$ with $r=0,1,2,\cdots$\,. 
They satisfy the following relations, 
\begin{align}
[h_{i,r},h_{j,s}]&=0 
\label{D2-1}
\\
[x^+_{i,r},x^-_{j,s}]&=\delta_{ij}h_{i,r+s}
\label{D2-2}
\\
[h_{i,0},x^\pm_{j,r}]&=\pm a_{ij} x^\pm_{j,r}
\label{D2-3}
\\
[h_{i,r+1},x^\pm_{j,s}]-[h_{i,r},x^\pm_{j,s+1}]&=\pm\frac{1}{2}a_{ij}\{h_{i,r},x^\pm_{j,s}\}
\quad \text{for $i,j$ not both $2$}
\label{D2-4}
\\
[h_{2,r},x^\pm_{2,s}]&=0 
\label{D2-4-1}
\\
[x^\pm_{i,r+1}, x^\pm_{j,s}]-[x^\pm_{i,r}, x^\pm_{j,s+1}]
&=\pm\frac{1}{2} a_{ij}\{x^\pm_{i,r}, x^\pm_{j,s}\}
\quad \text{for $i,j$ not both $2$}
\label{D2-5}
\\ 
%[x_{1,r}^\pm,x_{3,s}^\pm]&=0 
%\label{D2-6}
%\\
[x_{2,r}^\pm, x_{2,s}^\pm]&=0 
\label{D2-7}
\\
\bigl[x^\pm_{j,r},[x^\pm_{j,s},x^\pm_{2,t}]\bigr]
+\bigl[x^\pm_{j,s},[x^\pm_{j,r},x^\pm_{2,t}]\bigr]&=0
\qquad \text{for} \qquad j=1,3 
\label{D2-8}
\\
%[x_{2,r}^\pm, x_{2,s}^\pm]&=0
%\\
\bigl[[x^\pm_{1,r},x^\pm_{2,0}],[x^\pm_{3,s},x^\pm_{2,0}]\bigr]&=P^\pm_{r+s} \,.
\label{D2-9}
% \label{extP}
% \\
% \bigl[[\xi^-_{1,r},\xi^-_{2,0}],[\xi^-_{3,s},\xi^-_{2,0}]\bigr]&=K_{r+s} 
% \label{extK}
\end{align} 
The $\ZZ_2$-grading $p : \Y_D(\g)\to\ZZ_2$ is defined by setting 
$p(x^\pm_{2,r})=1$  for $r\in \ZZ_{\geq0}$ and $p(\text{others})=0$\,. 
The symmetrized Cartan matrix $a_{ij}$ is given in \eqref{car-mat}.  
\end{defin}

\begin{remark}
The Yangian $\Y_D(\g)$ includes the subalgebra $\Y(\alg{sl}(2))\oplus\Y(\alg{sl}(2))$, 
which is spanned by the generators $\{x_{i,r}^\pm,h_{i,r}\}$ with $i=1,3$\,, and $r\in \ZZ_{\geq0}$\,.
\end{remark}

We introduce the generators in $\Y_L(\g)$ inductively by,  for $r\in \ZZ_{\geq0}$\,, 
\begin{align}
&
%x_{i,r+1}^\pm \equiv \pm\frac{1}{a_{ii}}\tss[\tilh_{i,1},x_{i,r}^\pm] &
%&(~ i=1,3 ~) \nln
x_{1,r+1}^\pm = \pm\frac{1}{2}\tss[\tilh_{1,1},x_{1,r}^\pm] \,,\qquad 
x_{2,r+1}^\pm = \mp [\tilh_{1,1},x_{2,r}^\pm] \,,\qquad 
x_{3,r+1}^\pm = \mp \frac{1}{2}\tss[\tilh_{3,1},x_{3,r}^\pm] \,, \nln
%&
%x_{2,r+1}^\pm \equiv \mp [\tilh_{1,1},x_{2,r}^\pm] & & \nln
&
h_{i,r}= [x^+_{i,r},x_{i,0}^-] \qquad  (~ i=1,2,3 ~) \,, \nln
&
P^\pm_r= [[x_{1,r}^\pm,x_{2,0}^\pm],[x_{3,0}^\pm,x_{2,0}^\pm]] \,. 
\label{iso-map}
\end{align}
%for $r\geq\ZZ_+$\,. 

\begin{thm}
\label{thm:D2}
The Yangian $\Y_D(\g)$ is isomorphic to $\Y_L(\g)$\,. 
The isomorphism $\phi: \Y_D(\g)\to \Y_L(\g)$ is given by 
\begin{align}
h_{i,r}\mapsto h_{i,r} \,, \qquad 
x_{i,r}^\pm\mapsto x_{i,r}^\pm \,, \qquad  
P^\pm_r\mapsto P^\pm_r\,, 
\label{iso-phi}
\end{align}
where the image of $\phi$ is defined in \eqref{iso-map}\,. 
\end{thm}

\begin{remark} 
The original idea of Theorem \ref{thm:D2} is based on Levendorskii's construction of the Yangian associated with the finite-dimensional complex Lie algebra \cite{lev:gen}.
Its generalization to the super-Yangian associated with Lie superalgebra
$\alg{sl}(m|n)$ is obtained by Stukopin \cite{stu:ann}.
Compared to the defining relations of the Yangian $\Y(\alg{sl}(2|2))$,
the only differences are the extended Serre relations \eqref{D2-9}.
\end{remark}

\ms 
Theorem \ref{thm:D2} allows us to induce the Hopf algebra structure to $\Y_D(\g)$ from $\Y_L(\g)$ by the following commutative diagrams,
\begin{center}
$
\begin{CD}
\Y_D(\g) @>{\phi}>> \Y_L(\g)  \\
@V{\De_D}VV    @VV{\De}V \\
\Y_D(\g) \otimes \Y_D(\g)    @>>{\phi\otimes \phi}>  \Y_L(\g)\otimes \Y_L(\g)  
\end{CD}
$
\end{center}

\ms 

\begin{center}
$\begin{CD}
\Y_D(\g) @>{\phi}>> \Y_L(\g)  \\
@V{{\rm S}_D}VV    @VV{\rm S}V \\
\Y_D(\g)    @>>{\phi}>  \Y_L(\g)
\end{CD}$
\qquad\qquad  
\raise1ex\hbox{% 
$\begin{CD}
\Y_D(\g) @>{\phi}>> \Y_L(\g)  \\
@V{\ep_D}VV    @VV{\ep}V \\
\CC    @=  \CC  
\end{CD}$ 
}%
\end{center}
%\begin{align}
%&
%\begin{CD}
%\Y_D(\g) @>{\phi}>> \Y_L(\g)  \\
%@V{\De_D}VV    @VV{\De}V \\
%\Y_D(\g) \otimes \Y_D(\g)    @>>{\phi\otimes \phi}>  \Y_L(\g)\otimes \Y_L(\g)  
%\end{CD}
%\nln ~\nln
%& 
%\begin{CD}
%\Y_D(\g) @>{\phi}>> \Y_L(\g)  \\
%@V{{\rm S}_D}VV    @VV{\rm S}V \\
%\Y_D(\g)    @>>{\phi}>  \Y_L(\g)
%\end{CD}
%\qquad\qquad  
%\raise1ex\hbox{% 
%$\begin{CD}
%\Y_D(\g) @>{\phi}>> \Y_L(\g)  \\
%@V{\ep_D}VV    @VV{\ep}V \\
%\CC    @=  \CC  
%\end{CD}$ 
%}%
%\notag 
%\end{align}
Here $\De_D\,, {\rm S}_D$\,, and $\ep_D$ are induced coproduct, antipode, and 
counit in $\Y_D(\g)$, respectively. 
In particular, the manifest expressions of the Drinfeld coproduct $\De_D$ are obtained 
from \eqref{coprod} by the induction \eqref{iso-map}\,. 

\ms 

\begin{cor}
The Yangian $\Y_D(\g)$ in Definition \ref{def:D2} has the Hopf algebra structures with the coproducts given by 
\[
\Delta_D=(\phi^{-1}\otimes \phi^{-1}) 
\circ\De \circ \phi : \Y_D(\g)\to \Y_D(\g)\otimes \Y_D(\g)\,, 
\]
the antipodes 
\[
{\rm S}_D=\phi^{-1}\circ {\rm S} \circ \phi : \Y_D(\g)\to \Y_D(\g)\,, 
\]
and the counits
\[
\ep_D=\ep \circ\phi : \Y_D(\g)\to \CC\,. 
\] 
\end{cor}

\begin{proof}
Since the map $\phi$ is isomorphism by Theorem \ref{thm:D2}\,,  so is  $\phi^{-1}$\,. 
The maps $\De\,, {\rm S}\,, \ep$ define the Hopf algebra structures by Proposition \ref{prop:copro}\,. 
Thus, the composite maps  $\De_D\,, {\rm S}_D\,, \ep_D$ are homomorphism of the Yangian $\Y_D(\g)$\,. 
Hence, these maps define the Hopf algebra structures of $\Y_D(\g)$\,. 
\end{proof}

\begin{proof}[Proof of Theorem \ref{thm:D2}]

First, we prove that the map $\psi: \Y_L(\g)\to \Y_D(\g)$ defined by 
\begin{align}
&
h_{i,0}\mapsto h_{i,0} \,, 
&& x_{i,0}^\pm\mapsto x_{i,0}^\pm \,, 
&& P^\pm_0\mapsto P^\pm_0\,, 
\nln
&
\tilh_{i,1}\mapsto h_{i,1}-\frac{1}{2}h_{i,0}^2\,, 
&& x_{i,1}^\pm\mapsto x_{i,1}^\pm \,,   
&& P^\pm_1\mapsto P^\pm_1\,. 
\end{align}
%\begin{align}
%h_{i,r}\mapsto h_{i,r} \,, \qquad 
%x_{i,r}^\pm\mapsto x_{i,r}^\pm \,, \qquad  
%P^\pm_r\mapsto P^\pm_r\,, 
%\end{align}
is a homomorphism.  
%where the generators in $\Y_L(\g)$ are given in \eqref{iso-map}. 
%
Second, we show that the map $\phi: \Y_D(\g)\to \Y_L(\g)$ in \eqref{iso-phi} is a homomorphism. 
Since $\phi\circ\psi=\psi\circ\phi={\rm id.}$\,, we can conclude that $\phi$ in \eqref{iso-phi} is the isomorphism. 

\ms 
First, it is easy to see that the map $\psi: \Y_L(\g)\to \Y_D(\g)$ is a homomorphism. 
Indeed, the relations \eqref{yan-1}--\eqref{yan-serre} immediately follow from 
%\eqref{D2-1}, \eqref{D2-2}, \eqref{D2-5}, \eqref{D2-7}, \eqref{D2-9} 
\eqref{D2-1}--\eqref{D2-9} 
as truncation of the degree. 
In particular, by \eqref{D2-3} and \eqref{D2-4}, 
the relation \eqref{h1x0} is derived as 
\begin{align}
[\tilh_{i,1},x_{j,0}]
&\mapsto [h_{i,1},x_{j,0}]-\frac{1}{2}\bigl\{h_{i,0},[h_{i,0},x_{i,0}^\pm]\bigr\}
\nln 
&=[h_{i,0},x_{j,1}]\pm \frac{1}{2}a_{ij} \{h_{i,0}, x_{i,0}^\pm\}
\mp  \frac{1}{2}a_{ij} \{h_{i,0}, x_{i,0}^\pm\}
\nln 
&= \pm a_{ij} x_{j,1}^\pm\,. 
\end{align}
The relations \eqref{deg-31} and \eqref{deg-32} follow from \eqref{D2-1}\,, 
\begin{align}
[ \tilh_{j,1}, h_{j,2}]& \mapsto [ h_{j,1}-\frac{1}{2}h_{j,0}^2, h_{j,2}]=0\qquad(j=1,3)\,, 
\notag \\
[\tilh_{1,1}, h_{2,2}] & \mapsto
[ h_{1,1}-\frac{1}{2}h_{1,0}^2, h_{2,2}]= 0\,. 
\end{align}

\ms

Second, we show that the map $\phi: \Y_D(\g)\to \Y_L(\g)$ in \eqref{iso-phi} is a homomorphism
for the relations  \eqref{D2-1}--\eqref{D2-9}, respectively. 
Hereafter, all calculations are understood in $\Y_L(\g)$\,.

\item 
{\bf Proof of \eqref{D2-3}.}
Recall the definition of $x^\pm_{i,r}$ in \eqref{iso-map} and the relation \eqref{yan-1}. 
Commuting $\tilh_{j,1}$ ($j=1,3 $) by $r$-times with 
\begin{align}
[h_{i,0},x^\pm_{j,0}]&=\pm a_{ij} x^\pm_{j,0} \el 
%\label{D2-3}
\end{align} 
we have the relation \eqref{D2-3} with $j=1, 3$\,.  
When $j=2$, we can prove the relation \eqref{D2-3} by using $\tilh_{1,1}$ 
instead of  $\tilh_{2,1}$\,.  
\qed

\ms 

To prove \eqref{D2-5}, we prepare Lemma \ref{lem:deg-2} and Proposition \ref{prop:deg2-boost}. 
\begin{lem}
\label{lem:deg-2}
The following relations hold in $\Y(\g)$\,, for $r\in \ZZ_{\geq 0}$\,,
$i,j\in \{1,2,3\}\,,$
\begin{align}
[\tilh_{i,1},x_{j,r}^\pm]&=\pm a_{ij} x_{j,r+1}^\pm 
\label{eq1}
\\
[x^\pm_{i,2}, x^\pm_{j,1}]-[x^\pm_{i,1}, x^\pm_{j,1}]
&=\pm\frac{1}{2} a_{ij}\{x^\pm_{i,1}, x^\pm_{j,0}\} 
\label{eq2}
\\
h_{i,2}= [x_{i,2}^+,x_{i,0}^-]=[x_{i,1}^+,x_{i,1}^-]&=[x_{i,0}^+,x_{i,2}^-]
\label{eq3}
\\
[h_{i,2},\tilh_{j,1}]&=0 
\label{eq4}
\\
[h_{i,2}, x^\pm_{j,r}]-[h_{i,1}, x^\pm_{j,r+1}]
&=\pm\frac{1}{2} a_{ij}\{h_{i,1}, x^\pm_{j,r}\}\,.
\label{eq5}
\end{align}
\end{lem}

\begin{proof}[Proof of Lemma \ref{lem:deg-2}]
We shall prove these relations in order. 

\ms 

The proof of \eqref{eq1} is parallel to that of \eqref{D2-3}. 
For $j=1,3$\,, commuting $\tilh_{j,1}$ with the relation \eqref{h1x0} by $r$-times, we have \eqref{eq1}. 
For $j=2$\,, we use $\tilh_{1,1}$\,, and  obtain \eqref{eq1}. 

\ms 

To prove \eqref{eq2}, we start with the degree one relations \eqref{x1x0}\,, 
\begin{align}
[x^\pm_{i,1}, x^\pm_{j,0}]-[x^\pm_{i,0}, x^\pm_{j,1}]
&=\pm\frac{1}{2} a_{ij}\{x^\pm_{i,0}, x^\pm_{j,0}\}\,.  \el 
\end{align}
Acting $\tilh_{k,1}$ on both hand sides and using the relation \eqref{eq1}, 
we have 
\begin{align}
a_{ki} {\rm R}_{ij}(1,0)+a_{kj} {\rm R}_{ij}(0,1)=0\,,
\end{align}
where we have introduced the notations 
\begin{align}
{\rm R}_{ij}(r,s)= [x^\pm_{i,r+1}, x^\pm_{j,s}]-[x^\pm_{i,r}, x^\pm_{j,s+1}]
\mp\frac{1}{2} a_{ij}\{x^\pm_{i,r}, x^\pm_{j,s}\} \,, 
\qquad (r,s\geq0)\,. 
\label{R-rs}
\end{align}
When $i\neq j$\,, by setting $k=i$ and $k=j$, we have 
\begin{align}
a_{ii} {\rm R}_{ij}(1,0)+a_{ij} {\rm R}_{ij}(0,1)&=0\,, \\
a_{ji} {\rm R}_{ij}(1,0)+a_{jj} {\rm R}_{ij}(0,1)&=0\,.  
\end{align}
Since $a_{ii}a_{jj}-a_{ij}a_{ji}\neq0$ for the Cartan matrix \eqref{car-mat}, 
we obtain  
\begin{align}
{\rm R}_{ij}(1,0)={\rm R}_{ij}(0,1)=0\,. 
\end{align}
When $i=j=1, 3$\,, by setting $k=i$\,, 
%When $i=j=1$ or $2$\,, setting $k=1$
%(when $i=j=3$\,, setting $k=3$, respectively), 
we have 
\begin{align}
{\rm R}_{ii}(1,0)+{\rm R}_{ii}(0,1)=0\,. 
\end{align}
Due to the symmetry ${\rm R}_{ii}(1,0)={\rm R}_{ii}(0,1)$\,, we obtain
\begin{align}
{\rm R}_{ii}(1,0)={\rm R}_{ii}(0,1)=0\,. 
\end{align}
When $i=2$, as $a_{22}=0$\,, we have 
\begin{align}
{\rm R}_{22}(1,0)= [x^\pm_{2, 2}, x^\pm_{2, 0}]-[x^\pm_{2,1}, x^\pm_{2,1}]\,.  
\end{align}
On the other hand, commuting $\tilh_{1,1}$ with \eqref{xx20}, we obtain  
\begin{align}
[x^\pm_{2, 2}, x^\pm_{2, 0}]+ [x^\pm_{2,1}, x^\pm_{2,1}]=0 \,. 
\end{align}
Hence, it holds that 
\begin{align}
[x^\pm_{2, 2}, x^\pm_{2, 0}]= [x^\pm_{2,1}, x^\pm_{2,1}]=0 \,. 
\end{align}
This means that, for any $i,j=1,2,3$, we have   
\begin{align}
{\rm R}_{ij}(1,0)=0\,. 
\end{align}
Thus, \eqref{eq2} is proved.  

\ms 

Now we move to \eqref{eq3}. 
Acting $\tilh_{k,1}$ on $[x_{i,0}^+,x_{i,0}^-]=h_{i,0}$\,, we have 
\begin{align}
a_{ki} ([x_{i,1}^+,x_{i,0}^-]-[x_{i,0}^+,x_{i,1}^-])=0\,. 
\end{align}
Chose $k$ such as $k\neq 4-i$\,. Because of $a_{ik}\neq 0$\,, we have 
\begin{align}
h_{i,1}= [x_{i,1}^+,x_{i,0}^-]=[x_{i,0}^+,x_{i,1}^-]\,. 
\end{align}
Repeating this procedure for the above relation, we obtain 
\begin{align}
0&= [x_{i,2}^+,x_{i,0}^-]-[x_{i,1}^+,x_{i,1}^-]\,,\nln 
0&= [x_{i,1}^+,x_{i,1}^-]- [x_{i,0}^+,x_{i,2}^-]\,. \el 
\end{align}
Combining these two relations and taking into account 
the definition of $h_{i,2}$\,, we have \eqref{eq3}. 

\ms 

The left hand side of \eqref{eq4} is rewritten as 
\begin{align}
[h_{i,2},\tilh_{j,1}]
=[[x_{i,2}^+,x_{i,0}^-],\tilh_{j,1}] 
=[[x_{i,1}^+,x_{i,1}^-],\tilh_{j,1}]\,,   
\end{align}
where the last equality is due to \eqref{eq3}. 
Furthermore, by \eqref{eq1}, we have 
\begin{align}
%[h_{i,2},\tilh_{j,1}]
%&=[[x_{i,2}^+,x_{i,0}^-],\tilh_{j,1}] \nln
%&=
[[x_{i,1}^+,x_{i,1}^-],\tilh_{j,1}]
= -a_{ji}([x_{i,2}^+,x_{i,1}^-]-[x_{i,1}^+,x_{i,2}^-]) \,. 
\end{align}
Here, it is noted that the right hand side vanishes due to the relations 
\eqref{deg-31}, \eqref{deg-32} and \eqref{eq1}. 
Thus, we prove the relation \eqref{eq4}. 

\ms

Finally, we shall prove the relation \eqref{eq5}.  
Commuting $x^\mp_{i,0}$ with the relation \eqref{eq2} 
and taking \eqref{eq3} into account, we have 
\begin{align}
&[h_{i,2}, x^\pm_{j,0}]-[h_{i,1}, x^\pm_{j,1}]
\mp\frac{1}{2} a_{ij}\{h_{i,1}, x^\pm_{j,0}\} \nln 
&+\de_{ij}(-1)^{\de_{i,2}}\bigl([x^\pm_{i,2}, h_{j,0}]-[x^\pm_{i,1}, h_{j,1}]
\mp\frac{1}{2} a_{ij}\{x^\pm_{i,1}, h_{j,0}\}\bigr)=0\,. 
\end{align}
Since the second line on the left hand side %proportional to $\delta_{ij}$ 
is zero due to 
\begin{align}
[x^\pm_{i,2}, h_{j,0}]-[x^\pm_{i,1}, h_{j,1}]
\mp\frac{1}{2} a_{ij}\{x^\pm_{i,1}, h_{j,0}\}  
&=\mp a_{ij}x^\pm_{i,2}+[h_{i,1}-\frac{1}{2}h_{i,0}^2\tss,x_{i,1}^\pm] \nln
&=\mp a_{ij}x^\pm_{i,2}+[\tilh_{i,1}\tss,x_{i,1}^\pm] \nln
&=0\,, 
\end{align}
it becomes 
\begin{align}
[h_{i,2}, x^\pm_{j,0}]-[h_{i,1}, x^\pm_{j,1}]=\pm\frac{1}{2} a_{ij}\{h_{i,1}, x^\pm_{j,0}\} \,. 
\end{align}
Then, noting that $\tilh_{j,1}$ commutes with $h_{i,2}$ by \eqref{eq4}
and acting $\tilh_{j,1}$ on the above relation by $r$-times, 
we obtain the last relation \eqref{eq5}; 
\begin{align}
[h_{i,2}, x^\pm_{j,r}]-[h_{i,1}, x^\pm_{j,r+1}]=\pm\frac{1}{2} a_{ij}\{h_{i,1}, x^\pm_{j,r}\} \,. 
\end{align}
This completes the proof of the lemma. 
\end{proof}

Lemma \ref{lem:deg-2} allows us to prove the following proposition,  
which will play an important role later. 
\begin{prop}
\label{prop:deg2-boost} 
The operators defined by 
\begin{align}
B_{ij} 
%= h_{i,2} -h_{i,0} \tilh_{i,1} -\frac{1}{6}h_{i,0}\left(h_{i,0}^2+\frac{a_{ij}^2}{2} \right)
=h_{i,2} -h_{i,0} \tilh_{i,1} -\frac{1}{6}h_{i,0}^3
-\frac{1}{12}a_{ij}^2\tss h_{i,0}
\label{deg2-boost} 
\end{align}
satisfy the relations, for $i,j=1,2,3$ and $r\in \ZZ_{\geq 0}$\,,
\begin{align}
[B_{ij},x_{j,r}^\pm]=\pm a_{ij} x_{j,r+2}^\pm\,. 
%\qquad (\tss r\geq0 \tss)\,. 
\end{align}
\end{prop}

\begin{proof}
%By \eqref{eq1}\,, \eqref{eq5} in Lemma \ref{lem:deg-2}, and 
%the definition $h_{i,1}=\tilh_{i,1}+\frac{1}{2}h_{i,0}^2$
Adding the following two relations, 
\begin{align}
[h_{i,2}, x^\pm_{j,r}]-[h_{i,1}, x^\pm_{j,r+1}]
&=\pm \frac{1}{2} a_{ij}\{h_{i,1}, x^\pm_{j,r}\}\,, 
\notag\\
[h_{i,1}, x^\pm_{j,r+1}]-[h_{i,0}, x^\pm_{j,r+2}]
&=\pm \frac{1}{2} a_{ij}\{h_{i,0}, x^\pm_{j,r+1}\}\,, 
\end{align}
we have 
\begin{align}
[h_{i,2}, x^\pm_{j,r}]-[h_{i,0}, x^\pm_{j,r+2}]
&=\pm \frac{1}{2} a_{ij}\{h_{i,1}, x^\pm_{j,r}\}
\pm \frac{1}{2} a_{ij}\{h_{i,0}, x^\pm_{j,r+1}\}\,.
\end{align} 
By the definition $h_{i,1}=\tilh_{i,1}+\frac{1}{2}h_{i,0}^2$\,, 
it becomes 
\begin{align}
[h_{i,2}, x^\pm_{j,r}]-[h_{i,0}, x^\pm_{j,r+2}]
&=[ \frac{1}{2} \{ \tilh_{i,1}, h_{j,0}\},  x^\pm_{j,r}] 
\pm \frac{1}{4} a_{ij}\{h_{i,0}^2, x^\pm_{j,r}\}
\nln
\Leftrightarrow 
\quad 
[h_{i,2} - \tilh_{i,1} h_{j,0}\, , x^\pm_{j,r}]
&=\pm a_{ij} x^\pm_{j,r+2}
\pm \frac{1}{4} a_{ij}\{h_{i,0}^2, x^\pm_{j,r}\}\,.
\end{align} 
%The second term on the right hand side can be expressed by commutators 
Using  the following identity, 
\begin{align}
[h^3_{i,0}, x^\pm_{j,r}]=\pm\frac{3}{2}a_{ij} \{h_{i,0}^2, x^\pm_{j,r}\}
\mp\frac{1}{2}a_{ij} [h_{i,0}, [h_{i,0},x^\pm_{j,r} ]]\,,  
\end{align}
we obtain 
\begin{align}
[h_{i,2} - \tilh_{i,1} h_{j,0}\, , x^\pm_{j,r}]
&=\pm a_{ij} x^\pm_{j,r+2}
+\frac{1}{6}[h^3_{i,0},  x^\pm_{j,r}]
\pm \frac{1}{12} a_{ij}[h_{i,0}, [h_{i,0},x^\pm_{j,r} ]]
\nln
\Leftrightarrow 
\quad 
[h_{i,2} - \tilh_{i,1} h_{j,0}-\frac{1}{6}h^3_{i,0}\,,  x^\pm_{j,r}]
&=\pm a_{ij} x^\pm_{j,r+2}
\pm \frac{1}{12} a_{ij}[h_{i,0}, [h_{i,0},x^\pm_{j,r} ]]
\nln
&=\pm a_{ij} x^\pm_{j,r+2}
+ \frac{1}{12} a^2_{ij} [h_{i,0},x^\pm_{j,r} ] \,. 
\end{align}
Thus, we get the desired relations, 
\begin{align}
[h_{i,2} - \tilh_{i,1} h_{j,0}-\frac{1}{6}h^3_{i,0}
- \frac{1}{12} a^2_{ij} h_{i,0} \,,  x^\pm_{j,r}]
=[B_{ij}\,,  x^\pm_{j,r}]
=\pm a_{ij} x^\pm_{j,r+2} \,. 
\end{align}
%This is proved by the direct computations,  
%\begin{align}
%[B_{ij},x_{j,r}]
%&=[h_{i,2} -h_{i,0} \tilh_{i,1} -\frac{1}{6}h_{i,0}^3
%-\frac{a_{ij}^2}{12}\tss h_{i,0}\,,x_{j,r}^\pm] \nln 
%&=[h_{i,1}\tss,x_{j,r+1}^\pm ]\pm \frac{a_{ij}}{2}\{h_{i,1},x_{j,r}^\pm\}
%\mp \frac{a_{ij}}{2}\{\tilh_{i,1},x_{j,r}^\pm\} 
%\mp \frac{a_{ij}}{2}\{h_{i,0},x_{j,r+1}^\pm\}
%\nln 
%&\quad \mp \frac{a_{ij}}{6}
%\bigl(x_{j,r}^\pm h_{i,0}^2+h_{i,0}x_{j,r}^\pm h_{i,0}+h_{i,0}^2x_{j,r}^\pm\bigr) 
%\mp \frac{a_{ij}^3}{12} x_{j,r}^\pm 
%\nln 
%&= [\tilh_{i,1}+\frac{1}{2}h_{i,0}^2\tss,x_{j,r+1}^\pm ]
%\pm \frac{a_{ij}}{4}\{h_{i,0}^2,x_{j,r}^\pm\}
%\mp \frac{a_{ij}}{2}\{h_{i,0},x_{j,r+1}^\pm\}
%\nln 
%&\quad \mp \frac{a_{ij}}{6}
%\Bigl(\frac{3}{2}\{h_{i,0}^2\tss, x_{j,r}^\pm\}
%-\frac{1}{2}[h_{i,0}\tss,[h_{i,0}\tss,x_{j,r}^\pm]]\Bigr) 
%\mp \frac{a_{ij}^3}{12} x_{j,r}^\pm 
%\nln 
%&= [\tilh_{i,1}\tss,x_{j,r+1}^\pm ]
%\pm  \frac{a_{ij}}{2}\{h_{i,0},x_{j,r+1}^\pm\}
%\mp \frac{a_{ij}}{2}\{h_{i,0},x_{j,r+1}^\pm\}
%\nln 
%&=\pm a_{ij} x_{j,r+2}^\pm\,.  
%\end{align}
%Here in the second equality we have used the relations \eqref{eq5} and \eqref{eq1}
%in Lemma \ref{lem:deg-2}. 
%The third equality is due to the definition $h_{i,1}=\tilh_{i,1}+\frac{1}{2}h_{i,0}^2$
%and the following identity; 
%\begin{align}
%x_{j,r}^\pm h_{i,0}^2+h_{i,0}x_{j,r}^\pm h_{i,0}+h_{i,0}^2x_{j,r}^\pm
%&=\frac{3}{2}\{h_{i,0}^2\tss, x_{j,r}^\pm\}
%-\frac{1}{2}[h_{i,0}\tss,[h_{i,0}\tss,x_{j,r}^\pm]]\,. 
%\end{align}
Hence, the proposition is proved. 
\end{proof}

We shall use the {\it degree-two raising operator} $B_{ij}$ in \eqref{deg2-boost} 
in the following proof.
% of \eqref{D2-5}. 

\item {\bf Proof of \eqref{D2-5}.}
Let us use the notation \eqref{R-rs}
with 
%
%To prove the relation \eqref{D2-5}, we need the degree two raising operator
%in addition to $\tilh_{i,1}$\,, which is explicitly given by the lemma�@.... 
%
%We define 
%\begin{align}
%{\rm R}_{ij}(r,s)\equiv [x^\pm_{i,r+1}, x^\pm_{j,s}]-[x^\pm_{i,r}, x^\pm_{j,s+1}]
%\mp\frac{1}{2} a_{ij}\{x^\pm_{i,r}, x^\pm_{j,s}\}
%\qquad (\tss r,s\geq0\tss)\,,  
%\end{align}
%
$i,j=1,2,3$ which are not both $2$\,. 
Note that ${\rm R}_{ij}(r,s)$ has the following symmetry; 
\begin{align}
{\rm R}_{ij}(r,s)={\rm R}_{ji}(s,r)
\qquad\text{if} \qquad (i,j)\neq (2,2) \,. 
\label{R-sym}
\end{align}
The relation \eqref{D2-5} is now expressed as  
%\begin{align}
${\rm R}_{ij}(r,s)=0$\,. 
%\end{align}
We will prove this by induction with respect to $K=r+s$. 

\ms 

When $K=0$\,, the relation ${\rm R}_{ij}(0,0)=0$ is nothing but \eqref{x1x0}. 
When $K=1$\,, we have already proved the relations 
${\rm R}_{ij}(1,0)={\rm R}_{ij}(0,1)=0$ in Lemma \ref{lem:deg-2}. 

\ms 

%
%Acting $\tilh_{k,1}$ on the both hand sides of ${\rm R}_{ij}(0,0)=0$\,, we have 
%\begin{align}
%a_{ki}{\rm R}_{ij}(1,0)+a_{kj}{\rm R}_{ij}(0,1)=0\,. 
%\label{aR+aR}
%\end{align}
%When $i\neq j$\,, setting $k=i$ and $j$\,, the above relations become 
%\begin{align}
%a_{ii}{\rm R}_{ij}(1,0)+a_{ij}{\rm R}_{ij}(0,1)&=0 \nln 
%a_{ji}{\rm R}_{ij}(1,0)+a_{jj}{\rm R}_{ij}(0,1)&=0\,. 
%\end{align}
%Since $a_{ii}a_{jj}-a_{ij}a_{ji}\neq0$ from the Cartan matrix in \eqref{car-mat} 
%in this case, we obtain
%\begin{align}
%{\rm R}_{ij}(1,0)={\rm R}_{ij}(0,1)=0
%\qquad \text{for}\qquad i\neq j\,. 
%\end{align}
%When $i=j=1,3$\,, setting $k=i$ in \eqref{aR+aR}, it becomes 
%\begin{align}
%{\rm R}_{ii}(1,0)+{\rm R}_{ii}(0,1)=0\,. 
%\end{align}
%Due to the symmetry \eqref{R-sym}, this relation implies 
%\begin{align}
%{\rm R}_{ii}(1,0)={\rm R}_{ii}(0,1)=0\,. 
%\end{align}
%
%Thus, we obtain 
%\begin{align}
%{\rm R}_{ij}(r,s)=0
%\qquad \text{for}\qquad 0\leq r+s\leq 1\,. 
%\end{align}

Next, we suppose that, for a fixed $K\geq1$\,, it holds that 
\begin{align}
{\rm R}_{ij}(l,K-l)=0
\qquad \text{for any}\qquad 0\leq l\leq K\,,  
\label{R-hypo}
\end{align}
and show that the above relation holds for $K+1$\,. 
% We have already shown that this is true for $K=0,1$\,. 

\ms 

As we did for $K=1$ case, 
commuting $\tilh_{k,1}$ with \eqref{R-hypo}, we have  
\begin{align}
a_{ki}{\rm R}_{ij}(l+1,K-l)+a_{kj}{\rm R}_{ij}(l,K-l+1)=0\,. 
\label{aR-general}
\end{align}
When $i\neq j$\,, setting $k=i$ and $j$\,, the above relation becomes
\begin{align}
a_{ii}{\rm R}_{ij}(l+1,K-l)+a_{ij}{\rm R}_{ij}(l,K-l+1)&=0\,, \nln
a_{ji}{\rm R}_{ij}(l+1,K-l)+a_{jj}{\rm R}_{ij}(l,K-l+1)&=0\,. 
\end{align}
Noting that $a_{ii}a_{jj}-a_{ij}a_{ji}\neq0$ ($i\neq j$) for the Cartan matrix \eqref{car-mat}, 
we obtain
\begin{align}
{\rm R}_{ij}(l+1,K-l)={\rm R}_{ij}(l,K-l+1)=0
\qquad \text{for}\qquad i\neq j\,. 
\end{align}
Thus, we prove that 
\begin{align}
{\rm R}_{ij}(l,K+1-l)=0
\qquad \text{for}\qquad i\neq j\,,\quad 0\leq l \leq K+1\,.  
\end{align}

When $i=j=1,3$ and $K$ is even, we denote $K=2m$ with $m\geq1$\,.  
In this case, by setting $k=i$\,, the relation \eqref{aR-general} becomes 
\begin{align}
{\rm R}_{ii}(l+1,2m-l)+{\rm R}_{ii}(l,2m-l+1)=0 
\qquad \text{for}\qquad 0\leq l\leq 2m\,. 
\label{RmRm}
\end{align}
In particular, when $l=m$\,, this relation is 
\begin{align}
{\rm R}_{ii}(m+1,m)+{\rm R}_{ii}(m,m+1)=0\,.  
\end{align}
Taking into account the symmetry \eqref{R-sym}, we obtain 
\begin{align}
{\rm R}_{ii}(m+1,m)={\rm R}_{ii}(m,m+1)=0\,.  
\end{align}
Next, setting $l=m+1$ in \eqref{RmRm} and substituting 
${\rm R}_{ii}(m+1,m)=0$\,, we get 
\begin{align}
{\rm R}_{ii}(m+2,m-1)={\rm R}_{ii}(m-1,m+2)=0\,. 
\end{align}
Repeating this procedure for $l=m+2,\cdots,2m$, we find that 
\begin{align}
{\rm R}_{ii}(l,2m+1-l)=0 
\qquad \text{for}\qquad 0\leq l\leq 2m+1\,. 
\end{align}

\ms 

Finally, when $i=j=1,3$ and $K$ is odd, we denote $K=2m'+1$ with $m'\geq0$\,.  
In addition to the induction hypothesis \eqref{R-hypo}, 
we also suppose that 
\begin{align}
R_{ii}(l,2m'-l)=0 
\qquad \text{for}\qquad 0\leq l\leq 2m'\,. 
\label{R-hypo2}
\end{align}
Commuting the generator $B_{ii}$ in \eqref{deg2-boost} with this relation,  
%which raises the degree by two, 
we have 
\begin{align}
R_{ii}(l+2,2m'-l)+R_{ii}(l,2m'-l+2)=0\,.  
%\qquad \text{for}\qquad 0\leq l\leq 2m\,. 
\end{align}
In particular, setting $l=m'$ and using the symmetry \eqref{R-sym}, 
we find that 
\begin{align}
R_{ii}(m'+2,m')=R_{ii}(m',m'+2)=0\,. 
\label{Rm2}
\end{align}
On the other hand, commuting $\tilh_{i,1}$ with the induction hypothesis 
\eqref{R-hypo} %with $i=j=1,3$ and $K=2m'+1$\,, 
it turns out that 
\begin{align}
{\rm R}_{ii}(l+1,2m'+1-l)+{\rm R}_{ii}(l,2m'+2-l)=0
\qquad \text{for}\qquad 0\leq l\leq 2m'+1\,. 
\label{rm-odd}
\end{align}
Setting $l=m'$\,, it reduces to 
\begin{align}
{\rm R}_{ii}(m'+1,m'+1)+{\rm R}_{ii}(m',m'+2)=0\,. 
\end{align}
Remembering \eqref{Rm2}, we find 
\begin{align}
{\rm R}_{ii}(m'+1,m'+1)=0\,. 
\end{align}
For \eqref{rm-odd}, 
repeating this procedure for $l=m'+1,\cdots, 2m'+1$ in order, we prove that 
\begin{align}
{\rm R}_{ii}(l,2m'+2-l)=0
\qquad \text{for}\qquad 0\leq l\leq 2m'+2\,. 
\end{align}

Therefore, the induction hypothesis \eqref{R-hypo} holds for $K+1$\,. 
This completes the proof of the relation \eqref{D2-5}. 
\qed

\item {\bf Proof of \eqref{D2-7}. } 
%\begin{align}
%[x_{2,r}^\pm, x_{2,s}^\pm]&=0 
%%\label{D2-7}
%\end{align} 
%
Introduce the notation, 
\begin{align}
{\rm Q}(r,s)=[x_{2,r}^\pm, x_{2,s}^\pm]\,. 
\end{align}
Since $x_{2,r}^\pm$ are odd and the square bracket is the super-commutator, 
we have the symmetry;
\begin{align}
{\rm Q}(r,s)={\rm Q}(s,r)\,. 
\label{Q-sym}
\end{align}
The relation \eqref{D2-7} is now expressed as ${\rm Q}(r,s)=0$ for any $r,s\geq0$\,. 
We will prove this by induction with respect to $K=r+s$. 

\ms

For $K=0, 1$\,, due to the relations \eqref{lie-2} and \eqref{xx20}, the following relations hold,
\begin{align}
{\rm Q}(0,0)=0 \qquad \text{and} \qquad {\rm Q}(1,0)={\rm Q}(0,1)=0\,. 
%\\
%{\rm Q}(r,s)=0 \qquad \text{for} \qquad 0\leq r+s\leq 1\,. 
\end{align}

\ms 

Next, we suppose that for a fixed $K\geq1$ the following relations hold, 
\begin{align}
{\rm Q}(l,K-l)=0 \qquad \text{for} \qquad 0\leq l\leq K\,. 
%{\rm Q}(r,s)=0 \qquad \text{for} \qquad 0\leq r+s\leq K\,. 
\end{align}
Commuting $\tilh_{1,1}$ and noting 
$[\tilh_{1,1}\tss,x_{2,r}^\pm]=\mp x_{2,r+1}^\pm$\,, we have
\begin{align}
{\rm Q}(l+1,K-l)+{\rm Q}(l,K-l+1)=0 \qquad \text{for} \qquad 0\leq l\leq K\,.
\label{QQ-gen}
\end{align}

\ms 

When $K$ is even, set $K=2m$ ($m\geq0$)\,. Then, the above relations become 
\begin{align}
{\rm Q}(l+1,2m-l)+{\rm Q}(l,2m-l+1)=0 \qquad \text{for} \qquad 0\leq l\leq 2m\,.
\label{QQ-even}
\end{align}
In particular, when $l=m$\,, we have 
\begin{align}
{\rm Q}(m+1,m)+{\rm Q}(m,m+1)=0\,. 
% \qquad \text{for} \qquad 0\leq l\leq 2m\,.
\end{align}
Taking into account the symmetry \eqref{Q-sym}, this gives us 
\begin{align}
{\rm Q}(m+1,m)={\rm Q}(m,m+1)=0\,. 
% \qquad \text{for} \qquad 0\leq l\leq 2m\,.
\end{align}
Then, substituting $l=m+1,m+2,\cdots,2m$ for the relation \eqref{QQ-even} in order, 
we obtain
\begin{align}
{\rm Q}(l,2m+1-l)=0 \qquad \text{for} \qquad 0\leq l\leq 2m+1\,.
\end{align}

\ms

When $K$ is odd, setting $K=2m'+1$ ($m'\geq0$), the relations \eqref{QQ-gen} become 
\begin{align}
{\rm Q}(l+1,2m'+1-l)+{\rm Q}(l,2m'+2-l)=0 \qquad \text{for} \qquad 0\leq l\leq 2m'+1\,.
\label{QQ-odd}
\end{align}
On the other hand, by the induction hypothesis, the following relations hold;  
\begin{align}
{\rm Q}(l,2m'-l)=0 \qquad \text{for} \qquad 0\leq l\leq 2m'\,.
\end{align}
Then, acting $B_{12}$ in \eqref{deg2-boost} on these relations 
and noting $[B_{12},x_{2,r}^\pm]=\mp x_{2,r+2}^\pm$\,, 
we get
\begin{align}
{\rm Q}(l+2,2m'-l)+{\rm Q}(l,2m'+2-l)=0 \qquad \text{for} \qquad 0\leq l\leq 2m'\,.
\label{QQ-odd2}
\end{align}
In particular, setting $l=m'$ in \eqref{QQ-odd} and \eqref{QQ-odd2} respectively, 
we obtain
\begin{align}
{\rm Q}(m'+1,n+1)+{\rm Q}(m',m'+2)&=0\,, \nln
{\rm Q}(m'+2,n)+{\rm Q}(m',m'+2)&=0\,. 
\end{align}
Because of  the symmetry \eqref{Q-sym}, this gives us 
\begin{align}
{\rm Q}(m'+2,m')={\rm Q}(m',m'+2)={\rm Q}(m'+1,m'+1)=0\,. 
\end{align}
Then, substituting $l=m'+1,m'+2,\cdots,2m'+1$ for \eqref{QQ-odd} in this order, 
we obtain
\begin{align}
{\rm Q}(l,2m'+2-l)=0 \qquad \text{for} \qquad 0\leq l\leq 2m'+2\,.
\end{align}

\ms 

Thus, for any fixed $K\geq1$\,, we have proved 
\begin{align}
{\rm Q}(l,K+1-l)=0 \qquad \text{for} \qquad 0\leq l\leq K+1\,.
\end{align}
This completes the induction procedure. 
\qed 

\ms 

In order to prove the relations \eqref{D2-1}, \eqref{D2-2}, \eqref{D2-4}, and  \eqref{D2-4-1}, 
we prepare Lemma \ref{lem:hh}. 

\begin{lem}
\label{lem:hh}
For a fixed $K\geq1$\,,  suppose that the following relations hold
\begin{align}
[h_{i,r},h_{j,s}]&=0 \\
[x^+_{i,r},x^-_{j,s}]&=\delta_{ij}h_{i,r+s} \\
[h_{i,r+1},x^\pm_{j,s}]-[h_{i,r},x^\pm_{j,s+1}]&=\pm\frac{1}{2}a_{ij}\{h_{i,r},x^\pm_{j,s}\}\,,  
%\quad \text{for $i,j$ not both $2$}
\end{align} 
for $r,s\in \ZZ_{\geq0}$ such that $r+s \leq K$\,. 
Then, the following relations hold, 
\begin{align}
[h_{i,l},h_{j,K+1-l}]=0 
\qquad \text{for}\qquad 0\leq l \leq K+1\,. 
\label{hh-toprove}
\end{align} 
\end{lem}

\begin{proof}[Proof of Lemma \ref{lem:hh}]
Set 
\begin{align}
{\rm S}_{ij}(r,s)=[h_{i,r+1},x^\pm_{j,s}]-[h_{i,r},x^\pm_{j,s+1}]
\mp\frac{a_{ij}}{2}\{h_{i,r},x^\pm_{j,s}\}\,. 
\label{Sij}
\end{align}
By the induction hypothesis, we have for, $1\leq l \leq K $\,, 
\begin{align}
0&={\rm S}_{ij}(l-1,K-l) 
%\qquad \text{for} \qquad %2\leq l \leq K 
%1\leq l \leq K 
\nln
&=[h_{i,l},x^\pm_{j,K-l}]-[h_{i,l-1},x^\pm_{j,K-l+1}]
\mp\frac{a_{ij}}{2}\{h_{i,l-1},x^\pm_{j,K-l}\}\,.
\end{align}
Commuting $x^\mp_{j,1}$ with the above relation from right, it becomes 
\begin{align}
0=[[h_{i,l},x^\pm_{j,K-l}],x^\mp_{j,1}]-[[h_{i,l-1},x^\pm_{j,K-l+1}],x^\mp_{j,1}]
\mp\frac{a_{ij}}{2}[\{h_{i,l-1},x^\pm_{j,K-l}\},x^\mp_{j,1}]\,.
\end{align}
The first term is computed as 
\begin{align}
[[h_{i,l},x^\pm_{j,K-l}],x^\mp_{j,1}]
&=-[[x^\pm_{j,K-l},x^\mp_{j,1}], h_{i,l}] -[x^\pm_{j,K-l},[h_{i,l}, x^\mp_{j,1}]]
\nln
&=\mp [h_{j, K-l+1}, h_{i,l}]
-[x^\pm_{j,K-l},[h_{i,l+1}, x^\mp_{j,0}]\pm\frac{a_{ij}}{2}\{h_{i,l}, x^\mp_{j,0}\}]
\nln
&=\mp [h_{j, K-l+1}, h_{i,l}]
-[[x^\pm_{j,K-l}, h_{i,l+1}], x^\mp_{j,0}]
-[h_{i,l+1}, [x^\pm_{j,K-l}, x^\mp_{j,0}]]
\nln
&\quad 
\mp\frac{a_{ij}}{2}\{[x^\pm_{j,K-l},h_{i,l}], x^\mp_{j,0}\}
\mp\frac{a_{ij}}{2}\{h_{i,l},[x^\pm_{j,K-l}, x^\mp_{j,0}]\}
\nln
&=\mp [h_{j, K-l+1}, h_{i,l}]
-[[x^\pm_{j,K-l}, h_{i,l+1}], x^\mp_{j,0}]
\mp [h_{i,l+1}, h_{j,K-l}]
\nln
&\quad 
\pm\frac{a_{ij}}{2}\{[h_{i,l}, x^\pm_{j,K-l}], x^\mp_{j,0}\}
-\frac{a_{ij}}{2}\{h_{i,l}, h_{j,K-l} \}
\nln
&=\pm( [h_{i,l}, h_{j, K-l+1}]-[h_{i,l+1}, h_{j,K-l}]) 
+[[ h_{i,l+1}, x^\pm_{j,K-l}], x^\mp_{j,0}]
\nln
&\quad 
\pm\frac{a_{ij}}{2}\{[h_{i,l}, x^\pm_{j,K-l}], x^\mp_{j,0}\}
-\frac{a_{ij}}{2}\{h_{i,l}, h_{j,K-l} \}\,. 
\end{align}
%\begin{align}
%[[h_{i,l},x^\pm_{j,K-l}],x^\mp_{j,1}]
%&= \pm[h_{i,l+1},h_{j,K-l}] \mp [h_{i,l},h_{j,K-l+1}]
%-[[h_{i,l+1},x^\pm_{j,K-l}],x^\mp_{j,0}] 
%\nln
%&\quad 
%\mp\frac{a_{ij}}{2} \{[h_{i,l},x^\pm_{j,K-l}],x^\mp_{j,0}\}
%+\frac{a_{ij}}{2}\{h_{i,l},h_{j,K-l}\}\,. 
%\end{align}
%
%
Since the second term is obtained from the first term 
by shifting the index $l\to l-1$\,, 
the sum of the first and second term reads 
\begin{align}
&\quad [[h_{i,l},x^\pm_{j,K-l}],x^\mp_{j,1}]-[[h_{i,l-1},x^\pm_{j,K-l+1}],x^\mp_{j,1}]
\nln 
&=
\mp( [h_{i,l+1}, h_{j,K-l}]-2 [h_{i,l}, h_{j, K-l+1}]+[h_{i,l-1}, h_{j, K-l+2}]) 
\nln
&\quad +[[ h_{i,l+1}, x^\pm_{j,K-l}]-[ h_{i,l}, x^\pm_{j,K-l+1}], x^\mp_{j,0}]
\nln
&\quad 
\pm\frac{a_{ij}}{2}\{[h_{i,l}, x^\pm_{j,K-l}]- [h_{i,l-1}, x^\pm_{j,K-l+1}], x^\mp_{j,0}\}
\nln
&\quad 
-\frac{a_{ij}}{2}\{h_{i,l}, h_{j,K-l} \}+\frac{a_{ij}}{2}\{h_{i,l-1}, h_{j,K-l+1} \}
\nln
&=
\mp( [h_{i,l+1}, h_{j,K-l}]-2 [h_{i,l}, h_{j, K-l+1}]+[h_{i,l-1}, h_{j, K-l+2}]) 
\nln
&\quad \pm \frac{a_{ij}}{2}[\{ h_{i,l}, x^\pm_{j,K-l}\}, x^\mp_{j,0}]
\nln
&\quad 
+\frac{a_{ij}^2}{4}\{\{h_{i,l-1}, x^\pm_{j,K-l}\}, x^\mp_{j,0}\}
\nln
&\quad 
-\frac{a_{ij}}{2}\{h_{i,l}, h_{j,K-l} \}+\frac{a_{ij}}{2}\{h_{i,l-1}, h_{j,K-l+1} \}
\nln
&=
\mp( [h_{i,l+1}, h_{j,K-l}]-2 [h_{i,l}, h_{j, K-l+1}]+[h_{i,l-1}, h_{j, K-l+2}]) 
\nln
&\quad 
+\frac{a_{ij}^2}{4}\{\{h_{i,l-1}, x^\pm_{j,K-l}\}, x^\mp_{j,0}\}
\nln
&\quad 
+\frac{a_{ij}}{2}\{h_{i,l-1}, h_{j,K-l+1} \}
\pm \frac{a_{ij}}{2}\{x^\pm_{j,K-l}, [h_{i,l}, x^\mp_{j,0}] \}\,. 
\end{align}
%\begin{align}
%[[h_{i,l-1},x^\pm_{j,K-l+1}],x^\mp_{j,1}]
%&=\pm( [h_{i,l-1}, h_{j, K-l}]-[h_{i,l}, h_{j,K-l+1}]) 
%+[[ h_{i,l}, x^\pm_{j,K-l+1}], x^\mp_{j,0}]
%\nln
%&\quad 
%\mp\frac{a_{ij}}{2}\{[h_{i,l-1}, x^\pm_{j,K-l+1}], x^\mp_{j,0}\}
%-\frac{a_{ij}}{2}\{h_{i,l-1}, h_{j,K-l+1} \}\,. 
%\end{align}
%the sum of the first and second term reads
%\begin{align}
%&\quad [[h_{i,l},x^\pm_{j,K-l}],x^\mp_{j,1}]-[[h_{i,l-1},x^\pm_{j,K-l+1}],x^\mp_{j,1}] \nln
%&=\pm\bigl([h_{i,l+1},h_{j,K-l}] -2[h_{i,l},h_{j,K-l+1}] + [h_{i,l-1},h_{j,K-l+2}]\bigr)\nln
%&\quad 
%-\frac{a_{ij}}{2}\{h_{j,K-l},h_{i,l}\} \mp\frac{a_{ij}}{2} \{x^\pm_{j,K-l},[h_{i,l},x^\mp_{j,0}]\} 
%%\nln 
%%&\quad 
%-\left(\frac{a_{ij}}{2}\right)^2\{\{h_{i,l-1},x^\pm_{j,K-l}\},x^\mp_{j,0} \}\,.
%\end{align}
Then, adding the last term
\begin{align}
\mp\frac{a_{ij}}{2}[\{h_{i,l-1},x^\pm_{j,K-l}\},x^\mp_{j,1}]
=-\frac{a_{ij}}{2} \{h_{i,l-1} , h_{j,K-l+1}\} 
\mp \frac{a_{ij}}{2}\{x^\pm_{j,K-l},[h_{i,l-1},x^\mp_{j,1}]\}
\end{align}
to above, we obtain the relation
\begin{align}
0
&=
[[h_{i,l},x^\pm_{j,K-l}],x^\mp_{j,1}]-[[h_{i,l-1},x^\pm_{j,K-l+1}],x^\mp_{j,1}]
\mp\frac{a_{ij}}{2}[\{h_{i,l-1},x^\pm_{j,K-l}\},x^\mp_{j,1}]\,.
\nln
&=
\mp( [h_{i,l+1}, h_{j,K-l}]-2 [h_{i,l}, h_{j, K-l+1}]+[h_{i,l-1}, h_{j, K-l+2}]) 
\nln
&\quad 
+\frac{a_{ij}^2}{4}\{\{h_{i,l-1}, x^\pm_{j,K-l}\}, x^\mp_{j,0}\}
\nln
&\quad 
\pm \frac{a_{ij}}{2}\{x^\pm_{j,K-l}, [h_{i,l}, x^\mp_{j,0}] - [h_{i,l-1}, x^\mp_{j,1}]\}
\nln
&=
\mp( [h_{i,l+1}, h_{j,K-l}]-2 [h_{i,l}, h_{j, K-l+1}]+[h_{i,l-1}, h_{j, K-l+2}]) 
\nln
&\quad 
+\frac{a_{ij}^2}{4}\left(\{\{h_{i,l-1}, x^\pm_{j,K-l}\}, x^\mp_{j,0}\}
-(-1)^{\de_{j, 2}}\{\{h_{i,l-1},  x^\mp_{j,0} \}, x^\pm_{j,K-l}\}
\right)\,. 
\end{align}
%\begin{align}
%0
%&=\pm\bigl([h_{i,l+1},h_{j,K-l}] -2[h_{i,l},h_{j,K-l+1}] + [h_{i,l-1},h_{j,K-l+2}]\bigr)
%\nln
%&\quad 
%+\left(\frac{a_{ij}}{2}\right)^2
%\left(\{x^\pm_{j,K-l},\{x^\mp_{j,0},h_{i,l-1}\}\} 
%-(-1)^{\de_{j,2}} \{x^\mp_{j,0},\{x^\pm_{j,K-l},h_{i,l-1}\}\} \right)\,. 
%\end{align}
Noting that the second line of the last equality is further simplified and it actually vanishes, 
\begin{align}
&\quad 
\{\{h_{i,l-1}, x^\pm_{j,K-l}\}, x^\mp_{j,0}\}
-(-1)^{\de_{j, 2}}\{\{h_{i,l-1},  x^\mp_{j,0} \}, x^\pm_{j,K-l}\}
\nln
&=[h_{i,l-1}, [x^\pm_{j,K-l}, x^\mp_{j,0}]] 
\nln
&=\pm [h_{i,l-1}, h_{j,K-l}] 
\nln
&=0\,. 
\end{align}
%\begin{align}
%&\quad 
%\{x^\pm_{j,K-l},\{x^\mp_{j,0},h_{i,l-1}\}\} 
%-(-1)^{\de_{j,2}} \{x^\mp_{j,0},\{x^\pm_{j,K-l},h_{i,l-1}\}\} 
%\nln
%&=[[x^\pm_{j,K-l},x^\mp_{j,0}],h_{i,l-1}]
%\nln
%&=\pm [h_{j,K-l},h_{i,l-1}]
%\nln
%&=0\,. 
%\end{align}
As a result, we have the following recursion relation
\begin{align}
[h_{i,l+1},h_{j,K-l}] -2[h_{i,l},h_{j,K-l+1}] + [h_{i,l-1},h_{j,K-l+2}]=0
\qquad \text{for} \qquad 1\leq l\leq K\,. 
\label{diff-eq}
\end{align}

\ms 

The recursion relation \eqref{diff-eq} implies that 
\begin{align}
[h_{i,l+1},h_{j,K-l}] - [h_{i,l},h_{j,K-l+1}] 
&= [h_{i,1},h_{j,K}]-[h_{i,0},h_{j,K+1}]  
%\nln
%&\textcolor{red}{=0}
\qquad \text{for} \qquad 1\leq l\leq K\,. 
\end{align}
Since $[h_{i,0},h_{j,K+1}] =0$\,, we have 
\begin{align}
[h_{i,l+1},h_{j,K-l}]=(l+1)[h_{i,1},h_{j,K}]
\qquad \text{for} \qquad 1\leq l\leq K\,. 
\end{align}
When $l=K$\,, it becomes 
\begin{align}
[h_{i,K+1},h_{j,0}]=(K+1)[h_{i,1},h_{j,K}]
\qquad \text{for} \qquad 1\leq l\leq K\,. 
\end{align}
This gives us $[h_{i,1},h_{j,K}]=0$ since the left hand side is zero. 
Therefore, we obtain 
\begin{align}
[h_{i,l},h_{j,K-l+1}] =0
\qquad \text{for} \qquad 0\leq l\leq K+1\,. 
\end{align}
This completes the proof of lemma. 
\end{proof}

We are now ready to prove \eqref{D2-1}, \eqref{D2-2}, \eqref{D2-4} 
and \eqref{D2-4-1}.

\item {\bf Proof of \eqref{D2-1}, \eqref{D2-2}, \eqref{D2-4} 
and \eqref{D2-4-1}. } 
At first, we show that two relations \eqref{D2-4} and \eqref{D2-4-1} are casted into 
the following unified form, 
\begin{align}
[h_{i,r+1},x^\pm_{j,s}]-[h_{i,r},x^\pm_{j,s+1}]
= \pm\frac{a_{ij}}{2}\{h_{i,r},x^\pm_{j,s}\}
\quad \text{for}\quad  i,j\in \{1,2,3\} \,. 
\label{D2-4gen}
\end{align}
In fact, if $(i,j)\neq (2,2)$\,, then  \eqref{D2-4gen} is nothing but \eqref{D2-4}\,.
If $(i,j)= (2,2)$\,, it 
reduces to
\begin{align}
[h_{2,r+1},x^\pm_{2,s}]-[h_{2,r},x^\pm_{2,s+1}]=0
\end{align}
for any $r,s\in \ZZ_+$ since $a_{22}=0$\,. This implies that 
\begin{align}
[h_{2,r},x^\pm_{2,s}]=[h_{2,0},x^\pm_{2,r+s}]\,. 
\end{align}

\ms 

Hence, it is enough for us to prove the relations \eqref{D2-1}, \eqref{D2-2}, and \eqref{D2-4gen}\,. 
We show these relations simultaneously by induction with respect to 
a number $K\in \ZZ_{\geq0}$ such as $r+s\leq K$\,. 

\ms

%
%Before proceeding the proof of induction, 
%it is noted that the relation \eqref{D2-4-1} is equivalent to \eqref{D2-4} 
%with $i=j=2$\,. Indeed, setting $i=j=2$ in \eqref{D2-4}, it becomes 
%for any $r,s\in \ZZ_+$
%\begin{align}
%[h_{2,r+1},x^\pm_{2,s}]-[h_{2,r},x^\pm_{2,s+1}]=0
%\end{align}
%because of $a_{22}=0$\,. By using these relations repeatedly, 
%we get 
%\begin{align}
%[h_{2,r},x^\pm_{2,s}]=[h_{2,0},x^\pm_{2,r+s}]\,. 
%\end{align}
%Since the generator $x^\pm_{2,r+s}$ is given by 
%\begin{align}
%x^\pm_{2,r+s}=\mp [\tilh_{1,1},x^\pm_{2,r+s-1}]
%=\mp {\rm ad}_{\tilh_{1,1}} (x^\pm_{2,r+s-1}) 
%=\cdots 
%=(\mp)^{r+s} {\rm ad}^{r+s}_{\tilh_{1,1}} (x^\pm_{2,0})
%\end{align}
%by the definition \eqref{iso-map} and the generator $h_{2,0}$
%commutes with $\tilh_{1,1}$\,, we have 
%\begin{align}
%[h_{2,r},x^\pm_{2,s}]=[h_{2,0},x^\pm_{2,r+s}]
%=(\mp)^{r+s} {\rm ad}^{r+s}_{\tilh_{1,1}} ([h_{2,0}, x^\pm_{2,0}]) 
%=0\,. 
%\end{align}
%Thus, the relation \eqref{D2-4-1} follows from \eqref{D2-4} with $i=j=2$. 
%Inversely, the relation \eqref{D2-4} with $i=j=2$ is automatically satisfies 
%by \eqref{D2-4}. 
%Therefore, it is sufficient to prove the following relations for any $i,j=1,2,3$\,;
%\begin{align}
%[h_{i,r+1},x^\pm_{j,s}]-[h_{i,r},x^\pm_{j,s+1}]
%= \pm\frac{a_{ij}}{2}\{h_{i,r},x^\pm_{j,s}\}\,. 
%%\label{D2-4gen}
%\end{align}
%%
%Now let us prove the relations \eqref{D2-1}, \eqref{D2-2} and \eqref{D2-4gen}
%by induction. 

When $K=0$\,, they hold by Definition \ref{def:lie} and \ref{def:yan}\,. 
When $K=1$\,, they are derived from Definition \ref{def:yan}
and \eqref{eq5} in Lemma \ref{lem:deg-2}. 

\ms

For a fixed $K\geq1$ and any $r,s\in \ZZ_{\geq0}$ such that $r+s\leq K$\,, 
suppose that the relations \eqref{D2-1}, \eqref{D2-2}, and \eqref{D2-4gen} hold. 
Then, we will show that these relations also hold for 
any $r,s\in \ZZ_+$ satisfying $r+s\leq K+1$\,. 

\ms

Firstly,  the relation \eqref{D2-1} holds for any $r,s$ such that $r+s\leq K+1$
by Lemma \ref{lem:hh}. 

\ms

Secondly, let us prove that the relation \eqref{D2-2} holds for $r+s\leq K+1$\,. 
Commuting $\tilh_{k,1}$ with the relations 
\begin{align}
[x^+_{i,l},x^-_{j,K-l}]=\de_{ij} h_{i,K} 
\qquad \text{for} \qquad 0\leq l \leq K\,, 
\end{align}
we have 
\begin{align}
a_{ki}[x^+_{i,l+1},x^-_{j,K-l}]-a_{kj}[x^+_{i,l},x^-_{j,K-l+1}] =0
\label{aaxxpm}
\end{align}
because of $[\tilh_{k,1}, h_{i,K}]=0$\,. 
When $i\neq j$\,, since $a_{ii}a_{jj}-a_{ij}a_{ji}\neq 0$\,, 
we get 
\begin{align}
[x^+_{i,l},x^-_{j,K-l+1}]=0
\qquad \text{for} \qquad 0\leq l \leq K+1\,. 
\end{align}
When $i=j$\,, it is always possible to chose $k$ such as $a_{ik}\neq0$\,. 
For this $k$\,, the relation \eqref{aaxxpm} implies that 
\begin{align}
[x^+_{i,l+1},x^-_{i,K-l}]=[x^+_{i,l},x^-_{i,K-l+1}]
\qquad \text{for} \qquad 0\leq l \leq K\,. 
\end{align}
In particular, setting $l=K$\,, we see 
by definition $h_{i,K+1}$ in \eqref{iso-map} that  
%on the left hand side, we have by definition \eqref{iso-map}, 
%\begin{align}
%[x^+_{i,K+1},x^-_{i,0}]=h_{i,K+1}\,. 
%\end{align}
%Thus, we obtain 
\begin{align}
[x^+_{i,l},x^-_{i,K+1-l}]=h_{i,K+1}
\qquad \text{for} \qquad 0\leq l \leq K+1\,. 
\end{align}
Therefore, we have proved that, for $i,j\in \{1,2,3\}$\,, it holds that  
\begin{align}
[x^+_{i,l},x^-_{j,K+1-l}]=\de_{ij} h_{i,K+1} 
\qquad \text{for} \qquad 0\leq l \leq K+1\,. 
\label{xpm-h}
\end{align}

\ms 

Finally, we are going to show \eqref{D2-4gen}\,. 
By the induction assumption, it holds that 
\begin{align}
0&={\rm S}_{ij}(l,K-l) 
\qquad \text{for} \qquad 0\leq l \leq K-1 
\nln
&=[h_{i,l+1},x^\pm_{j,K-l}]-[h_{i,l},x^\pm_{j,K-l+1}]
\mp\frac{a_{ij}}{2}\{h_{i,l},x^\pm_{j,K-l}\}\,.
\end{align}
Commuting $\tilh_{k,1}$ with $k$ such that $a_{kj}\neq 0$ with this relation, 
we have 
\begin{align}
0={\rm S}_{ij}(l,K-l+1) 
\qquad \text{for} \qquad 0\leq l \leq K-1\,, 
\end{align}
because we have already proved that $[\tilh_{i,1},h_{j,m}]=0$ for $0\leq m\leq K$. 
Then, the remaining task is to prove the cases of $l=K$ and $K+1$, 
\begin{align}
{\rm S}_{ij}(K,1)&= [h_{i,K+1},x^\pm_{j,1}]-[h_{i,K},x^\pm_{j,2}]
\mp\frac{a_{ij}}{2} \{h_{i,K},x^\pm_{j,1}\} =0 \,, 
\nln 
{\rm S}_{ij}(K+1,0)&= [h_{i,K+2},x^\pm_{j,0}]-[h_{i,K+1},x^\pm_{j,1}]
\mp\frac{a_{ij}}{2} \{h_{i,K+1},x^\pm_{j,0}\} =0 \,. 
\label{S-K1}
\end{align}
To prove these relations, we shall start with 
\begin{align}
{\rm R}_{ij}(r,s)=
[x^\pm_{i,r+1}, x^\pm_{j,s}]-[x^\pm_{i,r}, x^\pm_{j,s+1}]
\mp\frac{a_{ij}}{2} \{x^\pm_{i,r}, x^\pm_{j,s}\} =0\,, 
\end{align}
for $r,s\in \ZZ_{\geq0}$\,.
% where ${\rm R}_{ij}(r,s)$ is given in \eqref{R-rs}. 
Then, commuting $x^\mp_{i,0}$ with ${\rm R}_{ij}(r,s)=0$ and using 
the result \eqref{xpm-h}, we have 
\begin{align}
{\rm S}_{ij}(r,s)+\de_{ij}(-1)^{\de_{i,2}}\,{\rm S}_{ii}(r,s)=0
\qquad \text{for} \qquad 
r, s\in \ZZ_{\geq0}\,. 
\end{align}
When $i\neq j$\,, %the second term on the left hand side vanishes
we obtain
\begin{align}
{\rm S}_{ij}(r,s)=0
\qquad \text{for} \qquad 
r, s\in \ZZ_{\geq0}\,. 
\end{align}
This proves the relations \eqref{S-K1} by setting 
$(r,s)=(K,1)$ and $(K+1,0)$, respectively.  
When $i=j$\,, we have 
\begin{align}
{\rm S}_{ii}(r,s)+(-1)^{\de_{i,2}}\,{\rm S}_{ii}(r,s)=0
\qquad \text{for} \qquad 
r, s\in \ZZ_{\geq0}\,. 
\end{align}
In particular, by setting $(r,s)=(K,1)$ and $(K+1,0)$\,, the above relations reduce to 
\begin{align}
{\rm S}_{ii}(K,1)+(-1)^{\de_{i,2}}\,{\rm S}_{ii}(1,K)&=0\,, 
\nln
{\rm S}_{ii}(K+1,0)+(-1)^{\de_{i,2}}\,{\rm S}_{ii}(0,K+1)&=0\,. 
\end{align}
Here, we notice that, by the relations \eqref{eq5} and \eqref{eq1}  
in Lemma \ref{lem:deg-2}, the following quantities vanish
\begin{align}
{\rm S}_{ii}(1,K)
&=[h_{i,2},x^\pm_{i,K}]-[h_{i,1},x^\pm_{i,K+1}]
\mp\frac{a_{ii}}{2} \{h_{i,1},x^\pm_{i,K}\} =0\,, \nln
{\rm S}_{ii}(0,K+1)
&=[h_{i,1},x^\pm_{i,K+1}]-[h_{i,0},x^\pm_{i,K+2}]
\mp\frac{a_{ii}}{2} \{h_{i,0},x^\pm_{i,K+1}\} =0\,.  
\end{align}
Thus, the relations \eqref{S-K1} with $i=j$ are proved. 

\ms

Therefore, we have completed the induction procedure 
and proved 
\eqref{D2-1}, \eqref{D2-2} and \eqref{D2-4gen}, which includes the relations   
\eqref{D2-4} and 
\eqref{D2-4-1}\,. 
% as we mentions at the beginning of this proof. 
\qed  

\ms 

We are now left with proving the Serre relations \eqref{D2-8} and \eqref{D2-9}. 

\item {\bf Proof of \eqref{D2-8}.}
We start to prove the following lemma. 
\begin{lem}
\label{lem:h2} 
The operators defined by 
\begin{align}
\tilh_{i, 2} 
=h_{i,2} -h_{i,0} \tilh_{i,1} -\frac{1}{6}h_{i,0}^3
%-\frac{1}{12}a_{ij}^2\tss h_{i,0}
\label{h2} 
\end{align}
satisfy the following relations, for $i,j=1,2,3$ and $r\in \ZZ_{\geq 0}$\,,
\begin{align}
[\tilh_{i,2},x_{j,r}^\pm]=\pm a_{ij} x_{j,r+2}^\pm
\pm \frac{1}{12}a_{ij}^3\tss x_{j,r}^\pm\,. 
%\qquad (\tss r\geq0 \tss)\,. 
\end{align}
\end{lem}

\begin{proof}
By Proposition \ref{prop:deg2-boost} and the definition of 
$\tilh_{i,2}=B_{ij}+\dfrac{1}{12}a_{ij}^2h_{i,0}$\,, 
we immediately obtain
\begin{align}
[\tilh_{i,2},x_{j,r}^\pm]
=[B_{ij}+\dfrac{1}{12}a_{ij}^2h_{i,0}, x_{j,r}^\pm]
=\pm x_{j,r+2}^\pm \pm \dfrac{1}{12}a_{ij}^3 x_{j,r}^\pm\,. 
\end{align}
This proves the lemma. 
\end{proof}

We set, for $r, s, t\in \ZZ_{\geq0}$ and $j=1, 3$\,,  
\begin{align}
X(r,s\,;t)= [x^\pm_{j,r},[x^\pm_{j,s},x^\pm_{2,t}]]
+[x^\pm_{j,s},[x^\pm_{j,r},x^\pm_{2,t}]]\,. 
\end{align}
We prove \eqref{D2-8} by induction with respect to the number $K=r+s+t$\,. 
We may assume that $r\geq s$ without loss of generality. 

\ms 

When $K=0$\,, it is nothing but the Serre relations of the Lie algebra \eqref{lie-5}, 
\begin{align}
X(0,0\,;0)= 2 [x^\pm_{j,0},[x^\pm_{j,0},x^\pm_{2,0}]]=0\,, 
\qquad (\,j=1,3\,) \,. 
\end{align}
When $K=1$\,,
commuting $\tilh_{k, 1}$ with $X(0,0; 0)=0$\,, we get 
\begin{align}
a_{kj}X(1,0\,;0)+a_{k2}X(0,0\,;1)=0 \,. 
\end{align}
Setting $k=j$ and $2$ respectively, we obtain 
\begin{align}
X(1,0\,;0)=X(0,0\,;1)=0 
\end{align}
since $a_{jj}a_{22}-a_{j2}a_{2j}\neq 0$ ($j=1,3$)\,. 

\ms 

Next, we suppose that for a fixed $K\geq 0$ it holds that 
\begin{align}
&X(r,s\,;t)=0 \,, 
\label{ih1}
\\ 
&X(r+1,s\,;t)=X(r,s+1\,;t)=X(r,s\,;t+1)=0 \,,  
\label{ih2}
\end{align}
for $r, s, t \in \ZZ_{\geq 0}$ such as $r\geq s$\,, $r+s+t=K$\,. 
Then, we show that 
\begin{align}
&X(r+2,s\,;t)=X(r+1,s+1\,;t)=X(r,s+2\,;t)=0\,, 
\nln
&X(r+1,s\,;t+1)=X(r,s+1\,;t+1)=0\,,
\nln
&X(r,s\,;t+2)=0 \,.  
\label{eq:x1}
\end{align}
Commuting $\tilh_{k, 1}$ ($k=j, 2$) with three relations in \eqref{ih2}\,, we get 
\begin{align}
&X(r+2,s\,;t)+X(r+1,s+1\,;t)=0\,, \qquad 
X(r+1,s\,;t+1) =0\,, 
\nln
&X(r+1,s+1\,;t)+X(r,s+2\,;t)=0\,, \qquad 
X(r,s+1\,;t+1) =0\,, 
\nln
&X(r+1,s\,;t+1)+X(r,s+1\,;t+1)=0\,, \qquad 
X(r,s\,;t+2) =0\,. 
\label{eq:x2}
\end{align}
While, by Lemma \ref{lem:h2},  
commuting $\tilh_{k, 2}$  ($k=j, 2$)  with the relation in \eqref{ih1}\,, we get 
\begin{align}
&X(r+2,s\,;t)+X(r,s+2\,;t)=0\,, \qquad 
X(r,s\,;t+2) =0\,. 
\label{eq:x3}
\end{align}
Thus, \eqref{eq:x2} and \eqref{eq:x3} give \eqref{eq:x1}\,.  
This completes the induction procedure. 
Hence, we have proved \eqref{D2-8}\,.

\item {\bf Proof of \eqref{D2-9}.}
%\begin{align}
%\bigl[[x^\pm_{1,r},x^\pm_{2,0}],[x^\pm_{3,s},x^\pm_{2,0}]\bigr]&=P^\pm_{r+s} 
%%\label{D2-9}
%\end{align} 
For convenience, we introduce a notation, 
\begin{align}
Z_{(k,l,m,n)}=[[x^\pm_{1,k},x^\pm_{2,l}],[x^\pm_{3,m},x^\pm_{2,n}]]\,. 
\end{align}
By definition \eqref{iso-map}, $P_r^\pm=Z_{(r,0,0,0)}$\,.  
Then, the relations \eqref{D2-9} claim the following two relations,  
for $K\in \ZZ_{\geq0}$\,, %it holds that  
\begin{align}
&Z_{(K,0,0,0)} = Z_{(K-l,0,l,0)} \qquad \text{for} \qquad 0\leq l \leq K\,, 
\label{eq:p1}
\\
&[\,Z_{(K,0,0,0)}\,, J\,]=0 \qquad \text{for all} \qquad J\in \Y_L(\g)\,. 
\label{eq:p2}
\end{align}

\ms 

We show \eqref{eq:p1} and \eqref{eq:p2} by induction with respect to 
$K\in \ZZ_{\geq0}$\,.
When $K=0$\,, these relations hold by Definition \ref{def:lie} and \ref{def:yan}.  
%When $K=1$\,, 
%reduce to \eqref{lie-last} and \eqref{yan-serre}\,, respectively. 

\ms 

Suppose that for a fixed $K\in \ZZ_{\geq0}$ the relations 
\eqref{eq:p1} and \eqref{eq:p2} hold. 

\ms 

Firstly, we consider the relation \eqref{eq:p1}. 
Commuting $\tilh_{2,1}$ with \eqref{eq:p1}\,, 
it yields for $0\leq l \leq K$ 
\begin{align}
0&=[\tilh_{2,1}\,, Z_{(K-l,0,l,0)}]
\nln 
&=\mp (Z_{(K-l+1,0,l,0)}-Z_{(K-l,0,l+1,0)}) \,. 
%\qquad \text{for} \qquad 0\leq l \leq K\,.  
\end{align}
Hence, this proves that the relation \eqref{eq:p1} holds for $K+1$\,, 
\begin{align}
Z_{(K+1,0,0,0)} = Z_{(K+1-l,0,l,0)} \qquad \text{for} \qquad 0\leq l \leq K+1\,. 
\label{eq:p1K+1}
\end{align}

\ms 

Next, we show the relation \eqref{eq:p2} for $K+1$\,. 
For this purpose, it is sufficient to prove  
\begin{align}
[\,Z_{(K,0,0,0)}\,, J\,]=0 \qquad \text{for any} \qquad J\in 
\{~x_{j,0}^\pm\,, x_{1,1}^+~|~j=1,2,3~\}\,,  
\label{eq:p2a}
\end{align}
because the other generators in $\Y_L(\g)$ are generated by 
the above generators 
$x_{j,0}^\pm$ with $j=1,2,3$\,, and $x_{1,1}^+$\,. 
%$\{~x_{j,0}^\pm\,, x_{1,1}^+~|~j=1,2,3~\}$\,.  
The calculations are straightforward with the help of the following two lemmas. 

\begin{lem}
\label{lem:x13}
For $r,s\in \ZZ_{\geq0}$\,, it holds that 
\begin{align}
[x^\pm_{1,r}, x^\pm_{3,s}]=0\,. 
\end{align}
\end{lem} 
\begin{proof}
Commuting $\tilh_{1,1}$ by $r$ times and $\tilh_{3,1}$ by $s$ times respectively 
with 
\begin{align}
[x^\pm_{1,0}, x^\pm_{3,0}]=0\,,  
\notag 
\end{align}
we obtain the desired relations. 
\end{proof}

\begin{lem}
\label{lem:ffx}
For an odd generator $F\in \Y_L(\g)$\,, if $[F,F]=0$\,, 
then it holds that 
\begin{align}
[F, [F, x]]=0 \qquad \text{for any} \qquad x\in \Y_L(\g)\,. 
\end{align}
\end{lem} 
\begin{proof}
By the super-Jacobi identity, we have 
\begin{align}
[F, [F, x]]&=[[F, F] , x]-[F, [F, x]]\,. 
%0 \qquad \text{for any} \qquad x\in \Y_L(\g)\,. 
\end{align}
Since $[F, F]=0$\,, it gives us $[F, [F, x]]=0$\,. 
\end{proof}

\ms 
The commutativity with the Lie algebraic generators can be shown as follows.  

\ms 

With the generators $x_{1,0}^\pm$\,, it is calculated as 
\begin{align}
[\,Z_{(K+1,0,0,0)}\,, x_{1,0}^\pm\,]
&=[\,Z_{(0,0,K+1,0)}\,, x_{1,0}^\pm\,]
\nln
&=[[[x^\pm_{1,0},x^\pm_{2,0}],[x^\pm_{3,K+1},x^\pm_{2,0}]]\,, x_{1,0}^\pm]
\nln
&=[[x^\pm_{1,0},x^\pm_{2,0}],[x^\pm_{3,K+1},[x^\pm_{2,0},  x_{1,0}^\pm]]]
\nln
&=0\,. 
\end{align}
Here, the first equality is due to \eqref{eq:p1K+1}\,. 
The third equality is obtained from the Serre relation 
$[x^\pm_{1,0},[x^\pm_{1,0},x^\pm_{2,0}]]=0$ and Lemma \ref{lem:x13}. 
The last equality is owing to Lemma \ref{lem:ffx} because of 
${\rm deg}([x^\pm_{1,0},x^\pm_{2,0}])=1$ and 
$[[x^\pm_{1,0},x^\pm_{2,0}], [x^\pm_{1,0},  x_{2,0}^\pm]]=0$\,. 

\ms 
With the generators $x_{1,0}^\mp$\,, the commutation relation is computed as 
\begin{align}
[\,Z_{(K+1,0,0,0)}\,, x_{1,0}^\mp\,]
&=[\,Z_{(0,0,K+1,0)}\,, x_{1,0}^\mp\,]
\nln
&=[\,[[x^\pm_{1,0},x^\pm_{2,0}],[x^\pm_{3,K+1},x^\pm_{2,0}]]\,, x_{1,0}^\mp\,]
\nln
&=\pm [[h_{1,0},x^\pm_{2,0}],[x^\pm_{3,K+1},x^\pm_{2,0}]]
\nln
&=- [x^\pm_{2,0}, [x^\pm_{3,K+1},x^\pm_{2,0}]]
\nln
&=0\,,  
\end{align}
where the last equality is due to Lemma \eqref{lem:ffx}. 
The commutativity of $P^\pm_{K+1}=Z_{(K+1,0,0,0)}$ with the other 
Lie algebraic generators $x^\pm_{2,0}\,, x^\pm_{3,0}$ are proved by 
the similar computations.

\ms 

Finally, let us confirm the commutativity with the degree one generator $x^\pm_{1,1}$\,, 
\begin{align}
[\,Z_{(K+1,0,0,0)}\,, x_{1,1}^\pm\,]
&=[\,Z_{(1,0,K,0)}\,, x_{1,1}^\pm\,]
\nln
&=[[[x^\pm_{1,1},x^\pm_{2,0}],[x^\pm_{3,K},x^\pm_{2,0}]]\,, x_{1,1}^\pm]
\nln
&=[[x^\pm_{1,1},x^\pm_{2,0}],[x^\pm_{3,K},[x^\pm_{2,0},  x_{1,1}^\pm]]]
\nln
&=0\,. 
\end{align}
Here, the first equality is due to \eqref{eq:p1K+1} with $l=K$\,. 
The third equality is obtained from the Serre relation 
$[x^\pm_{1,1},[x^\pm_{1,1},x^\pm_{2,0}]]=0$ and Lemma \ref{lem:x13}. 
The last equality is owing to Lemma \ref{lem:ffx} because of 
${\rm deg}([x^\pm_{1,1},x^\pm_{2,0}])=1$ and 
\begin{align}
[[x^\pm_{1,1},x^\pm_{2,0}], [x^\pm_{1,1},  x_{2,0}^\pm]]
= -[[x^\pm_{1,1},[x^\pm_{1,1},  x_{2,0}^\pm]],x^\pm_{2,0}]
+[x^\pm_{1,1}, [x^\pm_{2,0},[x^\pm_{1,1},  x_{2,0}^\pm]]]
=0\,. 
\end{align}
Hence, we have proved that 
\begin{align}
[\,Z_{(K+1,0,0,0)}\,, J\,]=0 \qquad \text{for any} \qquad J\in \Y_L(\g)\,. 
%\label{eq:p2}
\end{align}
This completes the induction procedure, 
and the relation \eqref{D2-9} is proved. 

\ms 

Now, we have shown that 
the map $\phi : \Y_D(\g)\to \Y_L(\g)$ defined in \eqref{iso-phi}
is a homomorphism. 
Therefore, Theorem \ref{thm:D2} is proved. 
\end{proof}

%\section{Poincar\'e-Birkhoff-Witt  theorem}
%\label{sec:PBW}
%\setcounter{equation}{0}

\section*{Acknowledgment}

We would like to thank the hospitality of the School of Mathematics and Statistics, the University of Sydney.
Most of this work has been done during his stay at the university.
We thank Prof.\ Hiroyuki Yamane for the variable discussions.  
We also appreciate our colleagues, Prof.\ Yoshiyuki Koga, and 
Prof.\ Yuji Satoh, for their communication with the related subjects. 
We also thank Prof.\ Alessandro Torrielli for his correspondence. 
Finally, we are very grateful to Prof.\ Alex Molev for his variable comments 
and encouragement for this project. 
This work was supported by JSPS KAKENHI Grant Number JP19K03421. 
%Grant-in-Aid for Scientific Research(C)

%%%%%%%%%%%%%%%%%%%%%%%%%%%%%%%%%%%%%%%%%%%%%
\appendix 
\renewcommand{\theequation}{\Alph{section}.\arabic{equation} }

\section{Relation to Drinfeld's first realization }
\label{app:D1}
\setcounter{equation}{0}

In this Appendix, we propose the relation of our Yangian $\Y_L(\g)$ to 
the Drinfeld first realization \cite{b:yan, dri:qg}. 
For this purpose, we introduce the Lie algebraic generators 
\begin{align}
x_{i}^\pm&= x_{i,0}^\pm\,, \quad 
h_{i}=   h_{i,0}\,,  \quad 
P^\pm=  P_{0}^\pm 
\quad \text{with}\quad i=1,2,3, 
\label{gen:D1lie}
\end{align}
and the following {\it hatted} degree one generators, 
%
%It is interesting to see how the Yangian $\Y_L(\g)$ in Definition \ref{def:yan} is 
%relating the Drinfeld first relization \cite{b:yan}. 
%For this purpose, we introduce the generators 
%$x^\pm_{i}\tss,\,h_i\tss,\,P^\pm$ 
%and the hatted generators 
%$\hat{x}^\pm_{i}\tss,\,\hat{h}_i\tss,\,\hat{P}^\pm$ 
%($i=1,2,3$) as follows, 
\begin{align}
%x_{i}^\pm&= x_{i,0}^\pm\,, \quad  h_{i}=   h_{i,0}\,,  \quad P^\pm=  P_{0}^\pm 
%\el \\  
\hat{x}_{1}^+&=  x_{1,1}^+ -\frac{1}{2}\bigl( h_{1,0} x_{1,0}^+ 
- \sum_{\mu=3,4} E_{2\mu} E_{\mu1} \bigr)
\el\\ 
\hat{x}_{2}^+&=  x_{2,1}^+ -\frac{1}{2}\bigl( h_{2,0} x_{2,0}^+ 
+ E_{12} E_{31} + E_{34} E_{42}- E_{14} P_0^+  \bigr)
\el\\ 
\hat{x}_{3}^+&=  x_{3,1}^+ - \frac{1}{2}\bigl(h_{3,0} x_{3,0}^+ 
-\sum_{l=1,2}E_{l3} E_{4l} \bigr)
\el\\
\hat{h}_{1}&=  \tilh_{1,1}  + x_{1,0}^- x_{1,0}^+ 
- \frac{1}{2}\sum_{\mu=3,4} \bigl(E_{1\mu} E_{\mu1}- E_{2\mu} E_{\mu2}\bigr) 
\el\\ 
\hat{h}_{2}&=  \tilh_{2,1}  
-\frac{1}{2}\bigl(E_{12} E_{21}- E_{13} E_{31}+E_{34} E_{43}+ E_{24} E_{42} -P_0^- P_0^+\bigr) 
\el\\ 
\hat{h}_{3}&=  \tilh_{3,1} - x_{3,0}^- x_{3,0}^+ 
-\frac{1}{2} \sum_{l=1,2} \bigl( E_{l3} E_{3l}- E_{l4} E_{4l} \bigr)
\el\\ 
\hat{x}_{1}^- &=  x_{1,1}^-  - \frac{1}{2}\bigl(x_{1,0}^- h_{1,0} 
-\sum_{\mu=3,4}  E_{1\mu} E_{\mu2} \bigr) 
\el\\ 
\hat{x}_{2}^-&=  x_{2,1}^- - \frac{1}{2}\bigl(x_{2,0}^- h_{2,0} 
- E_{13} E_{21} - E_{24} E_{43} -P_0^- E_{41}   \bigr)
\el\\ 
\hat{x}_{3}^-&=  x_{3,1}^- - \frac{1}{2}\bigl(x_{3,0}^-  h_{3,0} 
+\sum_{l=1,2}  E_{l4} E_{3l} \bigr)
\el \\
\hat{P}^+&=  P_{1}^+ + C_0 P_0^+ 
\el \\
\hat{P}^-&= P_{1}^- + P_0^- C_0 \,.     
\label{gen-D1} 
\end{align}
Here, the generators $E_{ij}$ of the matrix algebra in \eqref{lie-mat} 
via the identification \eqref{lie-iso}. 

\ms

We expect that the generators in \eqref{gen:D1lie} and \eqref{gen-D1} satisfy 
the defining relations of Drinfeld's first realization.  
Though we have not checked all commutation relations, 
the following proposition supports our proposal. 
\begin{prop}
The generators in \eqref{gen-D1} satisfy the relations, 
\begin{align}
[\hat{h}_i,x^\pm_j]&=[h_i,\hat{x}^\pm_j]=\pm a_{ij} \hat{x}^\pm_j \\
[\hat{x}^+_i,x^-_j]&=[x^+_i,\hat{x}^-_j]=\de_{ij}\hat{h}_i \\
[\hat{x}^\pm_i,x^\pm_j] &= [x^\pm_i,\hat{x}^\pm_j] \\
[[\hat{x}^\pm_{1},x^\pm_{2}],[x^\pm_{3},x^\pm_{2}]]&=\hat{P}^\pm \,. 
\end{align}
\end{prop}

\begin{proof}
Using the relations \eqref{lie-1}--\eqref{lie-last}, \eqref{lie-mat} and \eqref{yan-1}--\eqref{yan-serre}, 
we can verify these relations by direct computations.  
\end{proof}

\begin{remark}
It is noted that the antipodes for the generators listed in \eqref{gen-D1}  do not 
contains the Lie algebraic generators \eqref{gen:D1lie}, {\it i.e.}
\begin{align}
{\rm S}(\hat{J})=-\hat{J} \qquad \text{for}\qquad 
\hat{J}\in \{\tss \hat{x}^\pm_{i}\tss,\,\hat{h}_i\tss,\,\hat{P}^\pm \tss\}\,. 
\end{align}
This is consistent with the fact that the Lie superalgebra $\g$ has a vanishing Killing form.  
\end{remark}

%%%%%%%%%%%%%%%%%%%%%%%%%%%%%%%%%%%%%%%%%%%%%

\end{document}